\documentclass[a4paper]{article}

\usepackage[utf8]{inputenc}
\usepackage{etex}
\usepackage{amssymb}
\usepackage{color}
\usepackage{amsmath,amscd,amsthm}
\usepackage[color]{xypic}
\usepackage{mathrsfs}
\usepackage{grffile}
\usepackage{ifthen}
\usepackage{multicol}
\usepackage{tikz,tikz-cd}
\usetikzlibrary{matrix,arrows,positioning,calc,decorations.markings,scopes}
\usetikzlibrary{arrows.meta}
\usepackage{enumerate}
\usepackage{todonotes}
\usepackage{mathtools}
\usepackage{blkarray}

\usepackage{chngcntr}
\counterwithin{figure}{section}

\usepackage{subcaption}
\usepackage{caption}

\DeclareMathOperator{\im}{im}

\DeclareMathOperator{\Ker}{Ker}
\DeclareMathOperator{\coker}{coker}
\DeclareMathOperator{\cost}{cost}

\DeclareMathOperator{\rank}{rank}

\DeclareMathOperator{\op}{op}
\DeclareMathOperator{\Id}{Id}
\DeclareMathOperator{\Offset}{\mathcal O}

\DeclareMathOperator{\bary}{bary}

\DeclareMathOperator{\gr}{gr}

\DeclareMathOperator{\Rk}{Rk}
\newcommand{\B}[1]{{\mathcal B}(#1)}

\newcommand{\barc}[1]{\mathrm \xi\ifthenelse{\equal{#1}{}}{}{(#1)}}

\newcommand{\F}{F}

\renewcommand{\dim}{\mathop{\bf dim}}

\newcommand{\R}{\mathbb R}

\newcommand{\N}{\mathbb N}

\def\int{{I}}
\newcommand{\Int}{\ensuremath{\mathcal{I}}}

\newcommand{\grid}{\mathcal{G}}

 \newcommand{\kVect}{\mathbf{Vect}}
  \newcommand{\Set}{\mathbf{Set}}
  
\renewcommand{\Im}{\mathrm{Im}}

\newcommand{\Rec}{\ensuremath{\mathcal{R}}}
\newcommand{\Sec}{\ensuremath{\mathcal{S}}}

\newcommand{\Rank}{\Rk}

\newcommand{\Rips}{\mathrm{Rips}}

\newcommand{\SRips}{\mathrm{SRips}}
\newcommand{\DRips}{\mathrm{DRips}}

\newcommand{\SupVR}{\mathrm{Rips^{\downarrow}}}
\newcommand{\IntVR}{\mathrm{Rips}}

\newcommand{\Z}{\mathbb Z}

\newcommand{\Top}{\mathbf{Top}}
\newcommand{\Simp}{\textup{\bf{Simp}}}
\newcommand{\Vect}{\textup{\bf{Vect}}}
\newcommand{\vect}{\textup{\bf{vect}}}
\newcommand{\id}{\textup{Id}}

\newcommand{\dmod}{\text{Mod}}

\newcommand{\FI}[1]{\mathcal S(#1)}
\newcommand{\HH}{{H}}

\numberwithin{equation}{section}
\newtheorem{theorem}{Theorem}
\numberwithin{theorem}{section}
\newtheorem{proposition}[theorem]{Proposition}
 
\newtheorem{corollary}[theorem]{Corollary}

\theoremstyle{definition}
\newtheorem{example}[theorem]{Example}
\newtheorem{remark}[theorem]{Remark}
\newtheorem{remarks}[theorem]{Remarks}
\newtheorem{definition}[theorem]{Definition}
\usepackage
[
colorlinks=true,
linkcolor=black,
anchorcolor=blue,
citecolor=black,
urlcolor=blue,
plainpages=false,
pdfpagelabels
]{hyperref}
\usepackage[capitalise]{cleveref}

\title{An Introduction to Multiparameter Persistence}

\author{Magnus Bakke Botnan and Michael Lesnick}
\date{}

\usepackage{graphicx}

\begin{document}

\maketitle
\begin{abstract}
In topological data analysis (TDA), one often studies the shape of data by constructing a filtered topological space, whose structure is then examined using persistent homology.  However, a single filtered space often does not adequately capture the structure of interest in the data, and one is led to consider multiparameter persistence, which associates to the data a space equipped with a multiparameter filtration.  Multiparameter persistence has become one of the most active areas of research within TDA, with exciting progress on several fronts.  In this article, we introduce multiparameter persistence and survey some of this recent progress, with a  focus on ideas likely to lead to practical applications in the near future.
\end{abstract}
\setcounter{tocdepth}{2}
\tableofcontents
\section{Introduction}
\label{sec:intro}

\subsection{Overview}
Topological data analysis (TDA) is a branch of data science which applies topology to study the \emph{shape} of data, i.e., the coarse-scale, global, non-linear geometric features of data.  Examples of such features include clusters, loops, and tendrils in point cloud data, as well as modes and ridges in functional data.  While the history of TDA dates back to the 1990's, 
in recent years the field has advanced rapidly, leading to a rich theoretical foundation, highly efficient algorithms and software, and many applications \cite{otter2017roadmap,carlsson2009topology,rabadan2019topological,giunti2021TDA}.
 
On a high level, the TDA pipeline generally involves two steps: Given a data set $X$ (e.g., a finite set of points in $\R^n$), we
\begin{enumerate}
\item construct a commutative diagram of topological spaces $F(X)$ whose topological structure encodes information about the shape of $X$, and then
\item analyze the structure of $F(X)$ using tools from topology and abstract algebra.
\end{enumerate}
The most common instance of this pipeline is \emph{(1-parameter) persistent homology}.  Here the diagram $F(X)$ is assumed to be indexed by a totally ordered set $A$ (e.g., $A=\R$, or $A=\N$), and one applies homology with field coefficients, along with an algebraic structure theorem, to obtain signatures of data called \emph{barcodes}. We give a detailed overview in \cref{sec:1D}.

However, there are many settings in TDA where one wishes to construct a diagram of spaces $F(X)$ indexed by a poset that is not totally ordered.  In particular, it is often  natural to consider a diagram indexed by a product $A$ of $n$ totally ordered sets, e.g. $A=\R^n$ or $A=\mathbb N^{\op}\times \R$; we call these \emph{($n$-parameter) filtrations}, or \emph{multifiltrations}.  The branch of TDA which studies data via multifiltrations is called \emph{multiparameter persistence} or (when homology is used to study these diagrams) \emph{multiparameter persistent homology} (MPH).  The cases $n=2,3$ are of primary interest and the $n=2$ case has received the most attention thus far.  As we will see in \cref{sec:multi-d-pipe,Sec:Multifiltrations}, the use of 2- and 3-parameter persistent homology is natural when working with 
\begin{itemize}
\item data with outliers or variations in density (indeed, standard constructions of 1-parameter persistent homology are notoriously unstable to outliers, which motivates a 2-parameter approach; see \cref{Sec:Not_Robust}),
\item time-varying data, or more generally, data equipped with a real-valued function (e.g., the partial charge function on the atom centers of a protein),
\item data with \emph{tendrils} emanating from a central core,
\item functional data (e.g. image data) with noise that is large in magnitude but localized (``spike noise"). 
\end{itemize}

Multiparameter persistence was introduced by Frosini and Mullazani, who considered multiparameter persistent homotopy groups \cite{frosini1999size}.  Multiparameter persistent homology was first considered by Carlsson and Zomorodian \cite{carlsson2009theory}; several of the key ideas about MPH we will encounter in this introduction first appeared in \cite{carlsson2009theory}.

MPH provides algebraic invariants of data called \emph{($n$-parameter) persistence modules}, simply by applying homology with field coefficients to a multifiltration; in the case that $A=\mathbb Z^n$, these are $\Z^n$-graded modules over a polynomial ring in $n$ variables, objects classically studied in commutative algebra \cite{miller2004combinatorial}.  Persistence modules can also be viewed as commutative diagrams of vector spaces, which have been studied extensively in quiver representation theory \cite{derksen2005quiver}.  Thus, classical work in algebra offers a highly developed set of mathematical tools for working with MPH.  However, the sorts of mathematical questions about these objects that arise in TDA are often very different than the ones encountered in classical work in algebra; as we will see, these questions have stimulated the development of a rich theory for MPH, distinct from classical work.

When $n\geq 2$, there are difficulties with defining the barcode of a $n$-parameter persistence module; see  \cref{sec:diff-barc} for a detailed discussion.  Quiver representation theory provides an illuminating perspective on this: multiparameter persistence modules have \emph{wild representation type}.  In brief, this means that in contrast to the 1-parameter case, the natural generalization of a barcode in the $n$-parameter setting is enormously complicated, and in general far too complex to work with in practice.  We give an overview of this quiver representation theory and its connections to MPH in  \cref{sec:multidbarcode}. 
 
Because of these difficulties, the standard methodology for data analysis using 1-parameter persistence tends not to extend naively to the multiparameter setting, and new ideas are needed to make practical use of multiparameter persistence.  Recently, there has been substantial progress on the development of such ideas.  These advances suggest that MPH is an especially natural approach to the study of the shape of data, and one with potential to have substantial impact on data science.  Yet, for this potential to be reached, much critical work remains.

In this article, we introduce multiparameter persistent homology, focusing on some of the topics that we feel are most fundamental and most likely to lead to practical applications in the near future.  We hope that this article will serve as an invitation for others--mathematicians, computer scientists, statisticians, and application specialists---to join the effort to realize the potential of this approach to data analysis.  While we have tried to make this article as accessible as possible, we do assume some familiarity with a few elementary mathematical topics that arise in the study of 1-parameter persistence, such as basic category theory, homology, and abstract algebra.

We mention three other introductory resources on multiparameter persistence: First, Karag\"uler's recent M.S. thesis \cite{karaguler2021survey} provides an introduction to the theoretical aspects of multiparameter persistence, covering some of the same territory as our article, as well as some topics not discussed here.  Second, two chapters in Dey and Wang's textbook \cite{dey2022computational} focus on multiparameter persistence, particularly its computational aspects.  Third, Hal Schenck's recent book \cite{schenckmathematical} has some material on the algebraic aspects of multiparameter persistence.

\subsection{Main Themes} 
This article explores six main themes which underly much of the recent progress in MPH:

\paragraph{Invariants} Though defining multiparameter barcodes is problematic, many simple invariants of persistence modules are available to us, which can serve as a surrogate for the barcode in applications.  Three simple invariants which play a central role in the MPH literature are the \emph{Hilbert function}, the \emph{rank invariant}, and the \emph{multigraded Betti numbers}.  As we will explain, the rank invariant is equivalent to each of two different barcode-like invariants, the \emph{fibered barcode} \cite{lesnick2015interactive} and the \emph{signed barcode} \cite{botnan2021signed}.  

This signed barcode in fact belongs to a family of closely related generalized barcode constructions, which can be defined either in terms of M\"obius inversion or relative homological algebra \cite{betthauser2022graded,mccleary2020edit,patel2018generalized,botnan2022bottleneck}.  In particular, the rank invariant admits a natural extension, the \emph{generalized rank invariant} \cite{kim2018generalized}, which has an analogous signed barcode \cite{asashiba2019approximation,kim2018generalized,botnan2021signed}.  

We introduce invariants of multiparameter persistence modules in \cref{subsec:invariants}, and consider signed barcodes in detail in \cref{sec:rankinv}.  

\paragraph{Visualization} Visualization of barcodes has been critical to the practical success of 1-parameter persistence.   The problem of visualizing (invariants of) MPH in a practical, computationally efficient way seems to be of similar importance to the success of MPH.  RIVET \cite{lesnick2015interactive,rivet},  a practical software for working with 2-parameter persistent homology, introduces a novel interactive visualization of the Hilbert function, the bigraded Betti numbers, and the fibered barcode.  We illustrate RIVET's visualization on several examples in \cref{subsec:visual}.  Other ideas for visualizing MPH include the signed barcode \cite{botnan2021signed}, the \emph{multiparameter persistent landscape} \cite{vipond2021multiparameter} and (in an important special case) the \emph{staircase code} of \cite{cai2020elder}.

\paragraph{Metrics and Stability} Metrics on barcodes play an especially important role in both theory and applications of 1-parameter persistence.  Analogously, to develop the theory and applications of MPH, we need good metrics in the multiparameter setting.  Perhaps surprisingly, though defining barcodes for MPH is problematic, we have well-behaved extensions of the standard metrics on barcodes to the multiparameter persistence modules.  Of these, the best known is the \emph{interleaving distance}, which extends the \emph{bottleneck distance} on persistence barcodes.  The interleaving distance is useful in the theory of MPH, e.g., in establishing that several MPH constructions are robust to outliers \cite{blumberg2020stability}.  It has been shown that computing the interleaving distance on $n$-parameter persistence module is NP-hard for $n\geq 2$ \cite{bjerkevik2019computing}, which motivates the search for a more computable surrogate.  The \emph{matching distance} \cite{cerri2013betti} has emerged as one attractive choice.  

The interleaving distance is an $\ell^\infty$-type distance, i.e., it can be defined in terms of the $\ell^\infty$ distance on Euclidean spaces \cite{bjerkevik2021ell}.  As such, it is insensitive to small-scale algebraic structure that can be important in applications.  Several recent works explore the problem of defining $\ell^p$-type distances on multiparameter persistence modules. 

We discuss metrics on persistence modules and the stability theory for MPH in \cref{sec:metrics}.  

\paragraph{Computation and Software} Efficient algorithms are a critical prerequisite to practical applications of MPH, as is user-friendly software implementing such algorithms.  Activity in the computational aspects of MPH has accelerated in recent years, especially in the 2-parameter setting, leading to several important advances which have lowered the barrier to applications.  Arguably, the most fundamental computational problem in MPH is \emph{minimal presentation computation}.  Recent work shows that the standard algorithm for computing barcodes in the 1-parameter setting extends to a similarly efficient algorithm for computing minimal presentations of 2-parameter persistence modules, which has been implemented in \cite{lesnick2019computing} and (in a substantially improved form) in \cite{kerber2021fast}.  Other notable progress on the computational aspects of MPH concerns the computation of certain bifiltrations \cite{corbet2021computing,lesnick2021computing,edelsbrunner2020simple,edelsbrunner2021multi}, the computation of distances between multiparameter persistence modules \cite{kerber2020efficient,bjerkevik2019computing,kerber2019exact,dey2018computing}, and decomposition of MPH modules into indecomposables \cite{dey2019generalized}.  
We discuss computation in several places where related concepts arise, specifically in \cref{subsec:decomp-algo,Sec:Density-Sensitive Multifiltrations,Sec:Interleavings,Sec:Matching_Dist,Sec:Computing_Pres_Res,Sec:Compuitng_Signed_Barcode}.

\paragraph{Barcodes in Special Cases} In some special cases of interest, we do have well-defined barcodes that are simple enough to work with.  In particular, a fundamental 2-parameter persistence invariant of $\R$-valued functions called \emph{interlevel set persistence} has simple barcodes analogous to those in the 1-parameter setting \cite{cochoy2020decomposition,bendich2013homology,botnan2020decomposition}.  A part of the MPH literature focuses on developing a theoretical and computational understanding of such special cases, and some interesting discoveries have been made. This will be discussed in \cref{sec:block}.

\paragraph{Applications} There have been several efforts to develop practical applications of MPH, e.g., to image analysis, computational chemistry, and shape analysis, which hint at the promise of multiparameter persistence as a practical approach to data analysis \cite{keller2018persistent,adcock2014classification,vipond2018multiparameter,betancourt2018pseudo,biasotti2008multidimensional,schiff2021characterizing,carriere2020multiparameter,chung2021multi,beltramo2021euler,vipond2021multiparameter}. 
Nevertheless, in spite of the widespread interest in MPH among the TDA community and very encouraging recent progress in the field, applications of MPH are still in their infancy.  Arguably, this is in large part because the algorithms and software tools needed to study data at scale using MPH have been introduced only very recently, and are still very much under development.  We discuss two recent applications of MPH to cancer imaging in \cref{Sec:Applications}.  

\subsubsection{Other Key Themes}
The literature of multiparameter persistence is growing rapidly; we have made no attempt to discuss all of the interesting and promising recent research, and the choice of topics here is biased towards our own research interests. We draw the reader's attention to two interesting and important themes which are not covered in this article:

\paragraph{Sheaf-theoretic Viewpoints on Multiparameter Persistence}
There is a large literature on connections between sheaves and generalized persistence, which is very closely related to MPH, e.g., \cite{berkouk2021ephemeral, kashiwara2018persistent,de2016categorified,curry2014sheaves}.  The sheaf-theoretic framework naturally extends the multiparameter persistence framework, and offers powerful tools (e.g., sheafification) that promise to be useful in some of the settings considered in our paper.

\paragraph{Infinitely Presented but Tame Modules} The diagrams of spaces one deals with in applications of multiparameter persistence are often indexed by $\R^n$.  In many cases, the homology modules of these diagrams are finitely presented, in which case their study essentially reduces to the study of persistence modules indexed by $\mathbb N^n$.  But it is also common to encounter settings where the persistence modules are not finitely presented, yet satisfy a weaker finiteness property called \emph{tameness}, where generators and relations appear not only at single points in $\R^n$ but along continuous curves.  The algebraic properties of tame modules have been explored in depth in recent work of Miller \cite{miller2020homological}, whose work was motivated in part by the study of fly wing morphology \cite{miller2015fruit}; see also \cite{miller2020essential,miller2020primary}.  Two recent (independent) works \cite{baryshnikov2021biparametric,budney2021bi} draw on ideas from singularity theory \cite{whitney1955singularities} to show that for $M$ a compact manifold, the 2-parameter sublevel persistent homology (see \cref{sec:1D-pipeline}) of a generic smooth function $f:M\to \R^2$ is tame; in fact, a precise description of the topology of the sublevel filtration can be given in terms of the singularities of $f$.  We believe that these ideas are likely to be useful in both theory and applications of multiparameter persistent homology.  Closely related ideas also appear in \cite[Section 3]{cerri2019geometrical}.  

\subsection{Outline}
\cref{Sec:Filtrations_Persistence_Modues} introduces multiparameter filtrations and persistence modules, the main objects of study in multiparameter persistence.  \cref{sec:1D} reviews 1-parameter persistent homology and its stability properties.   \cref{sec:multi-d-pipe} introduces multiparameter persistent homology: \cref{sec:diff-barc} discusses the difficulty of defining barcodes of multiparameter persistence modules, \cref{subsec:invariants} discusses invariants of such modules, \cref{subsec:visual} gives an example of RIVET's visualization, and \cref{Sec:Applications} discusses two applications of MPH to cancer imaging.  \cref{Sec:Multifiltrations} introduces several constructions of multifiltrations from data, and discusses the computation of some of these.  \cref{sec:metrics} considers metrics on multiparameter persistence modules and stability results for MPH.  \cref{Sec:Min_Pres_Res}
discusses the problem of computing minimal presentations and resolutions of persistence modules.
\cref{sec:multidbarcode} provides a brief introduction to quiver representation theory and its connections to MPH.  \cref{sec:rankinv} introduces signed barcodes of multiparameter persistence modules.  \cref{sec:block} considers special settings where bipersistence modules decompose into \emph{interval modules}, and thus have nice barcodes, including the important case of \emph{interlevel persistent homology}.

\subsection*{Acknowledgements}
We have learned much of what we know about multiparameter persistence though conversations with friends and colleagues.  We especially thank our past and current collaborators in this area, who have had a major influence on this article.  We also thank the authors of \cite{vipond2021multiparameter} for helpful feedback on our discussion of their work in \cref{Sec:Applications} and for sharing \cref{fig:vipond} with us. We would also like to express our gratitude to Benjamin Blanchette and Sliem el Ela for bringing errors and typos in an earlier draft to our attention. Finally, we thank an anonymous reviewer for carefully reading the paper and providing valuable input.

\section{Filtrations and Persistence Modules}\label{Sec:Filtrations_Persistence_Modues}

\subsection{Posets}\label{Sec:Posets}
Given posets $P$ and $Q$, we let $P\times Q$ denote the \emph{product poset}, i.e. the partial order on $P\times Q$ is given by $(x_1, y_1)\leq (x_2, y_2)$ if and only if $x_1\leq x_2$ and $y_1\leq y_2$.  More generally, the product of $n$ posets \[P_1\times P_2\times \cdots \times P_n\] is defined analogously.

The \emph{opposite poset} of $P$, denoted $P^{\rm op}$, is the poset with the same underlying set, with $x\leq y$ in $P^{\rm op}$ if and only if $y\leq x$ in $P$. 

\begin{definition}
An \emph{interval} in a poset $P$ is a non-empty subset $I$ of $P$ satisfying the two following conditions:
\begin{enumerate}
\item If $s,t\in I$ and $s\leq u \leq t$, then $u\in I$, 
\item If $s,t\in I$, then there are $s=u_0, \ldots, u_{m}=t\in I$ such that $u_i$ and $u_{i+1}$ are comparable for all $0\leq i < m$.
\end{enumerate}
\end{definition}

A poset $P$ can be considered as a category with objects the elements of $P$ and morphisms:
\[
{\rm Hom}(x,y) = \begin{cases}
\{*\} & (x\le y), \\
\varnothing & (x \not\le y).
\end{cases}
\]

\subsection{Filtrations}\label{Sec:Filtrations}  Let $\Top$ denote the category of topological spaces and continuous maps.  A \emph{($P$-indexed) filtration} is functor $F\colon P\to \Top$ such that if $x\leq y$, then $F_x\subseteq F_y$ and the map $F_{x,y}:F_x\to F_y$ is the inclusion. In plainer language, $F$ is a choice of a topological space $F_x$ for each $x\in P$ such that $F_x\subseteq F_y$ whenever $x \leq y$.   If $P=T_1\times \cdots \times T_n$ where each $T_i$ is a totally ordered set, then we call $F$ a \emph{multiparameter or $n$-parameter filtration}. 
When $n=2$, we often call $F$ a \emph{bifiltration}. For example, a bifiltration $F\colon\N\times \N\to \Top$ is a diagram of topological spaces of the form
\[
\begin{tikzcd}[ampersand replacement=\&,row sep=3ex]  
\vdots                    \&\vdots                     \&\vdots       \\ 
\F_{0,2}\arrow[hook]{r}\arrow[hook]{u}  \&\F_{1,2}\arrow[hook]{r}\arrow[hook]{u}   \&\F_{2,2}\arrow[hook]{r}\arrow[hook]{u}  \&\cdots   \\
\F_{0,1}\arrow[hook]{r}\arrow[hook]{u}  \&\F_{1,1}\arrow[hook]{r}\arrow[hook]{u}   \&\F_{2,1}\arrow[hook]{r}\arrow[hook]{u}  \&\cdots   \\
\F_{0,0}\arrow[hook]{r}\arrow[hook]{u}   \&\F_{1,0}\arrow[hook]{r}\arrow[hook]{u}   \&\F_{2,0}\arrow[hook]{r}\arrow[hook]{u}  \&\cdots   
\end{tikzcd}
\]
We will often consider filtrations valued in the category of simplicial complexes $\Simp$, which we regard as a subcategory of $\Top$ via geometric realization.

For fixed $P$, the $P$-indexed filtrations form a category whose morphisms are natural transformations.  That is, a morphism $\eta\colon F\to G$ in this category is a choice of a continuous map
$\eta_x\colon F_x\to G_x$ for each $x\in P$, such that for all $x\leq y\in P$, the following diagram commutes:
\[
\begin{tikzcd}[ampersand replacement=\&]  
\F_x\arrow["\eta_x",swap]{d}\arrow[hook]{r}  \& F_y\arrow["\eta_y"]{d}    \\
G_x\arrow[hook]{r} \& G_y
\end{tikzcd}
\]
For example, given filtrations $F,G\colon \N\to \Top$, a morphism is an extension of $F$ and $G$ to a commutative diagram of the following shape:
\[
\begin{tikzcd}[ampersand replacement=\&]  
F_{0}\arrow[hook]{r}\arrow[dashed]{d}  \&F_{1}\arrow[hook]{r}\arrow[dashed]{d}   \&F_{2}\arrow[dashed]{d}\arrow[hook]{r} \&\cdots   \\
G_{0}\arrow[hook]{r}  \&G_{1}\arrow[hook]{r}  \&G_{2}\arrow[hook]{r}  \&\cdots  
\end{tikzcd}
\]

\paragraph{Sublevel Filtrations}

\begin{definition}\label{Def:Sublevel}
Given a poset $P$, a topological space $W$, and a (not necessarily continuous) function $\gamma\colon W\to P$, we define the \emph{sublevel filtration of $\gamma$} as
\begin{align*}
&\mathcal S^\uparrow(\gamma)\colon P\to \Top,\\
&\mathcal S^\uparrow(\gamma)_p=\{w\in W\mid \gamma(w)\leq p\}.
\end{align*}  
\end{definition}

In the case $P=\R^n$, this construction was introduced in \cite{frosini1999size}.  A filtration isomorphic to a sublevel filtration is said to be \emph{1-critical} \cite{carlsson2010computing,blumberg2017universality}, whereas a filtration that is not 1-critical is said to be \emph{multicritical}.  Note that if $F\colon P\to \Top$ is 1-critical, then each element of $\bigcup_{p\in P} F_p$ first appears in the filtration at some unique minimal index in $P$, whereas this need not be the case for multicritical filtrations. We will see later that 1-critical and multicritical bifiltrations both arise naturally from data.

\subsection{Persistence Modules}\label{Sec:Persistence_Modules} Fix a field $k$, for example, $k=\Z/2\Z$, and let $\Vect$ denote the category of $k$-vector spaces and linear maps.  For $P$ a poset, a \emph{($P$-indexed) persistence module}, or simply a \emph{$P$-module}, is a functor $M\colon P\to \Vect$. For $x\leq y$ we let $M_{x,y}$ denote the morphism $M_x \to M_y$.  If $P = T_1\times \cdots \times T_n$ where each $T_i$ is a totally ordered set, then $M$ is called a \emph{multiparameter or $n$-parameter persistence module}. When $n=2$, we also call $M$ a \emph{bipersistence module}.  

For example, a bipersistence module $M\colon\N\times \N\to \Vect$ is a commutative diagram of $k$-vector spaces of the form:
\[
\begin{tikzcd}[ampersand replacement=\&,row sep=2.5ex]  
\vdots                    \&\vdots                     \&\vdots       \\ 
M_{(0,2)}\arrow{r}\arrow{u}  \&M_{(1,2)}\arrow{r}\arrow{u}   \&M_{(2,2)}\arrow{r}\arrow{u}  \&\cdots   \\
M_{(0,1)}\arrow{r}\arrow{u}  \&M_{(1,1)}\arrow{r}\arrow{u}   \&M_{(2,1)}\arrow{r}\arrow{u}  \&\cdots   \\
M_{(0,0)}\arrow{r}\arrow{u}   \&M_{(1,0)}\arrow{r}\arrow{u}   \&M_{(2,0)}\arrow{r}\arrow{u}  \&\cdots   
\end{tikzcd}
\]
A persistence module is said to be \emph{pointwise finite-dimensional}, or \emph{p.f.d.} if $\dim(M_x)<\infty$ for all $x\in P$.  

Applying $i^{\mathrm{th}}$ homology with coefficients in $k$ to each space and each inclusion map in a filtration $F\colon P\to \Top$ yields a persistence module $H_i(F)\colon P\to \Vect$.

For fixed $P$, the $P$-modules form a category $\Vect^P$ whose morphisms are the natural transformations, just as for filtrations. By working pointwise, the category $\Vect^P$ inherits direct sums (i.e., coproducts), images, kernels, cokernels from $\Vect$.  In fact, $\Vect^P$ is an \emph{abelian category}.  For example, given a morphism $f\colon M\to N$ of $P$-modules, $\ker f$ is the $P$-module defined by $(\ker f)_x:=\ker (f_x)$ for all $x\in P$, with the internal maps of $\ker f$ given by restricting $M_{x,y}$ to $(\ker f)_x$. We say that $M\colon P\to \Vect$ is \emph{indecomposable} if $M\cong M'\oplus M''$ implies that $M'=0$ or $M''=0$. 

\begin{definition}
For an interval $I$ in a poset $P$, the \emph{interval module} $k_I$ is defined by
\begin{align*}
(k_I)_x&=
\begin{cases}
k &{\textup{if }} x\in I, \\
0 &{\textup{ otherwise},}
\end{cases}
& (k_I)_{x,y}=
\begin{cases}
\Id_k &{\textup{if }} x\leq y\in I,\\
0 &{\textup{ otherwise}.}
\end{cases}
\end{align*}
\end{definition}
It is not hard to verify that $k_I$ is indecomposable; see e.g. \cite[Proposition 2.2]{botnan2018algebraic}.  

\subsection{Persistence Modules as Multigraded Modules}
\label{sec:multigraded}
As the name suggests, persistence modules can be seen as modules over a ring, and sometimes this perspective is very useful. 
We explain this just in the case where $P=\Z^n$, though the story is much the same when $P=\R^n$ \cite{lesnick2015theory}.

\begin{definition}
Let $\mathbf e_i$ denote the $i^{\mathrm{th}}$ standard basis vector in $\Z^n$.  An \emph{$n$-grading} on a $k[x_1,\dots,x_n]$-module $M$ is a vector space decomposition \[M=\oplus_{z\in \Z^n} M_z\] such that $x_i M_z\subset M_{z+\mathbf e_i}$ for all $z\in \Z^n$ and $i\in \{1,\ldots,n\}$.   A $k[x_1,\dots,x_n]$-module $M$ is said to be \emph{$n$-graded} if it comes equipped with an $n$-grading.
\end{definition}

A \emph{morphism} $f:M\to N$ of $n$-graded modules is a module homomorphism (in the usual sense) such that $f(M_z)\subset N_z$ for all $z\in \N^n$.  With these morphisms, the $n$-graded modules form a category $\dmod$.

\begin{proposition}[Carlsson, Zomorodian \cite{carlsson2009theory,corbet2018representation}]\label{Prop:ModDiagramEquivalence}
The categories $\Vect^{\Z^n}$ and $\dmod$ are equivalent.
\end{proposition}

$k[x_1,\dots,x_n]$-modules are the basic objects of study in commutative algebra, and thus have a rich and highly developed theory \cite{eisenbud1995commutative}.
Proposition \ref{Prop:ModDiagramEquivalence} allows us to adapt standard language and constructions for $k[x_1,\dots,x_n]$-modules to the study of persistence modules, provided those constructions make sense in the $n$-graded setting.  Fortunately, as a rule of thumb, definitions and results about $k[x_1,\dots,x_n]$-modules do tend to carry over to the $n$-graded setting, and often become simpler there.

\section{1-Parameter Persistent Homology}
\label{sec:1D}
In this section, we briefly review 1-parameter persistent homology.   For a more thorough introduction, see, e.g., \cite{oudot2015persistence,edelsbrunner2017persistent}. 

As mentioned the introduction, 1-parameter persistent homology provides invariants of data called \emph{barcodes}.  The barcode is a collection of intervals lying in a fixed  totally ordered set, e.g., $\R$ or $\N$.  Intuitively, intervals in the barcode represent geometric features of the data, and the length of the interval is a measure of size of the corresponding feature.  For instance, the circular structure in \cref{fig:circle} is represented by the long interval in \cref{fig:barcode-circle}.

\subsection{The Persistent Homology Pipeline}
\label{sec:1D-pipeline}
Let $T$ be a totally ordered set. The standard pipeline for constructing barcodes proceeds in three steps, as illustrated by the following diagram:
\begin{center}
\framebox[1.1\width]{Data} $\Rightarrow$ \framebox[1.1\width]{$T$-Filtration} $\xRightarrow{{\rm Homology}}$ \framebox[1.1\width]{$T$-Module}
$\xRightarrow{{\rm Decomposition}}$ \framebox[1.1\width]{Barcode}
\end{center}
Depending on the type of data and the information one wants to capture, the construction of a filtration from data can take many forms.  The following examples are standard:

\begin{example}\label{Ex:Offset_Filtration}
For $X\subset \R^n$, define the  \emph{offset filtration} of $X$ to be $\Offset(X):=\mathcal S^\uparrow(d_X)$, where $d_X\colon\R^n\to [0,\infty)$ is the distance function to $X$, i.e.,
\[d_X(y)=\inf_{x\in X} \|y-x\|.\]
\end{example}
Note that $\Offset(X)_r$ is simply the union of the closed balls of radius $r$ centered at the points of $X$.

\begin{remark}
The offset filtration is topologically equivalent to two well-known simplicial filtrations, the \emph{\v Cech filtration}, and a sub-filtration thereof called the \emph{Delaunay filtration} (also known as the \emph{$\alpha$-filtration}) \cite{edelsbrunner2010computational}.  The Delaunay filtration is used frequently in persistent homology computations involving points in low-dimensional Euclidian space (particularly in $\R^3$).  
\end{remark}

\begin{definition}
For a finite metric space $X$ and $r\in [0,\infty)$, let $N(X)_r$ denote the neighborhood graph of $X$, i.e., the graph with vertex set $X$ containing the edge $[x,y]$ if and only if $d(x,y)\leq r$.  Let $\Rips(X)_r$ denote the \emph{clique complex} on $N(X)_r$, i.e., the largest simplicial complex with 1-skeleton $N(X)_r$.  These simplicial complexes assemble into a filtration $\Rips(X)\colon [0,\infty)\to \Simp$, called the \emph{(Vietoris-)Rips filtration}.
\label{def:Rips}
\end{definition}

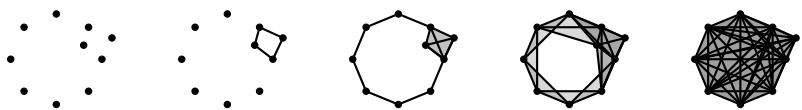
\begin{figure}[h!]
\begin{center}
\definecolor{aqaqaq}{rgb}{0.6274509803921569,0.6274509803921569,0.6274509803921569}

\begin{center}
\begin{tikzpicture}[scale=.6]
\coordinate (a) at (1,0);
\coordinate (b) at (0,1);
\coordinate (c) at (-1,0);
\coordinate (d) at (0,-1);
\coordinate (e) at (.707,.707);
\coordinate (f) at (-.707,.707);
\coordinate (g) at (-.707,-.707);
\coordinate (h) at (.707,-.707);
\coordinate (i) at  (1.22,.47);
\coordinate (j) at (.6,.31);
\coordinate (k) at (.6,.9);

\draw [fill] (a) circle [radius=0.07];
\draw [fill] (b) circle [radius=0.07];
\draw [fill] (c) circle [radius=0.07];
\draw [fill] (d) circle [radius=0.07];
\draw [fill] (e) circle [radius=0.07];
\draw [fill] (f) circle [radius=0.07];
\draw [fill] (g) circle [radius=0.07];
\draw [fill] (h) circle [radius=0.07];
\draw [fill] (i) circle [radius=0.07];
\draw [fill] (j) circle [radius=0.07];
\end{tikzpicture}
\hskip20pt
\begin{tikzpicture}[scale=.6]
\coordinate (a) at (1,0);
\coordinate (b) at (0,1);
\coordinate (c) at (-1,0);
\coordinate (d) at (0,-1);
\coordinate (e) at (.707,.707);
\coordinate (f) at (-.707,.707);
\coordinate (g) at (-.707,-.707);
\coordinate (h) at (.707,-.707);
\coordinate (i) at  (1.22,.47);
\coordinate (j) at (.6,.31);
\coordinate (k) at (.6,.9);

\draw[thick] (a) -- (i);
\draw[thick] (a) -- (j);
\draw[thick] (e) -- (i);
\draw[thick] (e) -- (j);

\draw [fill] (a) circle [radius=0.07];
\draw [fill] (b) circle [radius=0.07];
\draw [fill] (c) circle [radius=0.07];
\draw [fill] (d) circle [radius=0.07];
\draw [fill] (e) circle [radius=0.07];
\draw [fill] (f) circle [radius=0.07];
\draw [fill] (g) circle [radius=0.07];
\draw [fill] (h) circle [radius=0.07];
\draw [fill] (i) circle [radius=0.07];
\draw [fill] (j) circle [radius=0.07];
\end{tikzpicture}
\hskip20pt
\begin{tikzpicture}[scale=.6]
\coordinate (a) at (1,0);
\coordinate (b) at (0,1);
\coordinate (c) at (-1,0);
\coordinate (d) at (0,-1);
\coordinate (e) at (.707,.707);
\coordinate (f) at (-.707,.707);
\coordinate (g) at (-.707,-.707);
\coordinate (h) at (.707,-.707);
\coordinate (i) at  (1.22,.47);
\coordinate (j) at (.6,.31);
\coordinate (k) at (.6,.9);

\fill[color=aqaqaq, fill=aqaqaq, fill opacity=0.4] (a) -- (e) -- (i) -- cycle;
\fill[color=aqaqaq, fill=aqaqaq, fill opacity=0.4] (a) -- (e) -- (j) -- cycle;
\fill[color=aqaqaq, fill=aqaqaq, fill opacity=0.4] (a) -- (i) -- (j) -- cycle;
\fill[color=aqaqaq, fill=aqaqaq, fill opacity=0.4] (e) -- (i) -- (j) -- cycle;

\draw[thick] (a) -- (i);
\draw[thick] (a) -- (j);
\draw[thick] (e) -- (i);
\draw[thick] (e) -- (j);
\draw[thick] (a) -- (e);
\draw[thick] (i) -- (j);
\draw[thick] (b) -- (e);
\draw[thick] (b) -- (f);
\draw[thick] (c) -- (f);
\draw[thick] (c) -- (g);
\draw[thick] (d) -- (g);
\draw[thick] (d) -- (h);
\draw[thick] (a) -- (h);
\draw[thick] (a) -- (i);
\draw[thick] (a) -- (j);
\draw [fill] (a) circle [radius=0.07];
\draw [fill] (b) circle [radius=0.07];
\draw [fill] (c) circle [radius=0.07];
\draw [fill] (d) circle [radius=0.07];
\draw [fill] (e) circle [radius=0.07];
\draw [fill] (f) circle [radius=0.07];
\draw [fill] (g) circle [radius=0.07];
\draw [fill] (h) circle [radius=0.07];
\draw [fill] (i) circle [radius=0.07];
\draw [fill] (j) circle [radius=0.07];
\end{tikzpicture}
\hskip20pt
\begin{tikzpicture}[scale=.6]
%
\coordinate (a) at (1,0);
\coordinate (b) at (0,1);
\coordinate (c) at (-1,0);
\coordinate (d) at (0,-1);
\coordinate (e) at (.707,.707);
\coordinate (f) at (-.707,.707);
\coordinate (g) at (-.707,-.707);
\coordinate (h) at (.707,-.707);
\coordinate (i) at  (1.22,.47);
\coordinate (j) at (.6,.31);
\coordinate (k) at (.6,.9);

\fill[color=aqaqaq, fill=aqaqaq, fill opacity=0.4] (a) -- (e) -- (i) -- cycle;
\fill[color=aqaqaq, fill=aqaqaq, fill opacity=0.4] (a) -- (e) -- (j) -- cycle;
\fill[color=aqaqaq, fill=aqaqaq, fill opacity=0.4] (a) -- (i) -- (j) -- cycle;
\fill[color=aqaqaq, fill=aqaqaq, fill opacity=0.4] (e) -- (i) -- (j) -- cycle;

\fill[color=aqaqaq, fill=aqaqaq, fill opacity=0.4] (a) -- (b) -- (e) -- cycle;
\fill[color=aqaqaq, fill=aqaqaq, fill opacity=0.4] (a) -- (b) -- (i) -- cycle;
\fill[color=aqaqaq, fill=aqaqaq, fill opacity=0.4] (a) -- (b) -- (j) -- cycle;
\fill[color=aqaqaq, fill=aqaqaq, fill opacity=0.4] (a) -- (e) -- (h) -- cycle;

\fill[color=aqaqaq, fill=aqaqaq, fill opacity=0.4] (e) -- (b) -- (f) -- cycle;
\fill[color=aqaqaq, fill=aqaqaq, fill opacity=0.4] (b) -- (f) -- (c) -- cycle;
\fill[color=aqaqaq, fill=aqaqaq, fill opacity=0.4] (f) -- (c) -- (g) -- cycle;
\fill[color=aqaqaq, fill=aqaqaq, fill opacity=0.4] (c) -- (g) -- (d) -- cycle;
\fill[color=aqaqaq, fill=aqaqaq, fill opacity=0.4] (g) -- (d) -- (h) -- cycle;
\fill[color=aqaqaq, fill=aqaqaq, fill opacity=0.4] (d) -- (h) -- (a) -- cycle;

\fill[color=aqaqaq, fill=aqaqaq, fill opacity=0.4] (b) -- (e) -- (i) -- cycle;
\fill[color=aqaqaq, fill=aqaqaq, fill opacity=0.4] (b) -- (e) -- (j) -- cycle;
\fill[color=aqaqaq, fill=aqaqaq, fill opacity=0.4] (b) -- (j) -- (i) -- cycle;

\fill[color=aqaqaq, fill=aqaqaq, fill opacity=0.4] (h) -- (e) -- (i) -- cycle;
\fill[color=aqaqaq, fill=aqaqaq, fill opacity=0.4] (h) -- (e) -- (j) -- cycle;
\fill[color=aqaqaq, fill=aqaqaq, fill opacity=0.4] (h) -- (i) -- (j) -- cycle;

\fill[color=aqaqaq, fill=aqaqaq, fill opacity=0.4] (e)-- (f) -- (j) -- cycle;
\draw[thick] (a) -- (i);
\draw[thick] (a) -- (j);
\draw[thick] (e) -- (i);
\draw[thick] (e) -- (j);
\draw[thick] (a) -- (e);
\draw[thick] (i) -- (j);
\draw[thick] (b) -- (e);
\draw[thick] (b) -- (f);
\draw[thick] (c) -- (f);
\draw[thick] (c) -- (g);
\draw[thick] (d) -- (g);
\draw[thick] (d) -- (h);
\draw[thick] (a) -- (h);
\draw[thick] (a) -- (i);
\draw[thick] (a) -- (j);
\draw[thick] (a) -- (b);
\draw[thick] (b) -- (c);
\draw[thick] (c) -- (d);
\draw[thick] (d) -- (a);
\draw[thick] (e) -- (f);
\draw[thick] (f) -- (g);
\draw[thick] (g) -- (h);
\draw[thick] (h) -- (a);
\draw[thick] (e) -- (h);
\draw[thick] (h) -- (i);
\draw[thick] (h) -- (j);
\draw[thick] (i) -- (b);
\draw[thick] (j) -- (b);
\draw[thick] (f) -- (j);

\draw [fill] (a) circle [radius=0.07];
\draw [fill] (b) circle [radius=0.07];
\draw [fill] (c) circle [radius=0.07];
\draw [fill] (d) circle [radius=0.07];
\draw [fill] (e) circle [radius=0.07];
\draw [fill] (f) circle [radius=0.07];
\draw [fill] (g) circle [radius=0.07];
\draw [fill] (h) circle [radius=0.07];
\draw [fill] (i) circle [radius=0.07];
\draw [fill] (j) circle [radius=0.07];
\end{tikzpicture}
\hskip20pt
\begin{tikzpicture}[scale=.6]
%
\coordinate (a) at (1,0);
\coordinate (b) at (0,1);
\coordinate (c) at (-1,0);
\coordinate (d) at (0,-1);
\coordinate (e) at (.707,.707);
\coordinate (f) at (-.707,.707);
\coordinate (g) at (-.707,-.707);
\coordinate (h) at (.707,-.707);
\coordinate (i) at  (1.22,.47);
\coordinate (j) at (.6,.31);
\coordinate (k) at (.6,.9);

\fill[color=aqaqaq, fill=aqaqaq, fill opacity=1] (i) -- (e) -- (b) -- (f) -- (c) -- (g) -- (d) -- (h) -- cycle;
\draw[thick] (a) -- (b);
\draw[thick] (a) -- (c);
\draw[thick] (a) -- (d);
\draw[thick] (a) -- (e);
\draw[thick] (a) -- (f);
\draw[thick] (a) -- (g);
\draw[thick] (a) -- (h);
\draw[thick] (a) -- (i);
\draw[thick] (a) -- (j);

\draw[thick] (b) -- (c);
\draw[thick] (b) -- (d);
\draw[thick] (b) -- (e);
\draw[thick] (b) -- (f);
\draw[thick] (b) -- (g);
\draw[thick] (b) -- (h);
\draw[thick] (b) -- (i);
\draw[thick] (b) -- (j);

\draw[thick] (c) -- (d);
\draw[thick] (c) -- (e);
\draw[thick] (c) -- (f);
\draw[thick] (c) -- (g);
\draw[thick] (c) -- (h);
\draw[thick] (c) -- (i);
\draw[thick] (c) -- (j);

\draw[thick] (d) -- (e);
\draw[thick] (d) -- (f);
\draw[thick] (d) -- (g);
\draw[thick] (d) -- (h);
\draw[thick] (d) -- (i);
\draw[thick] (d) -- (j);

\draw[thick] (e) -- (f);
\draw[thick] (e) -- (g);
\draw[thick] (e) -- (h);
\draw[thick] (e) -- (i);
\draw[thick] (e) -- (j);

\draw[thick] (f) -- (g);
\draw[thick] (f) -- (h);
\draw[thick] (f) -- (i);
\draw[thick] (f) -- (j);

\draw[thick] (g) -- (h);
\draw[thick] (g) -- (i);
\draw[thick] (g) -- (j);

\draw[thick] (h) -- (i);
\draw[thick] (h) -- (j);

\draw[thick] (i) -- (j);

\draw [fill] (a) circle [radius=0.07];
\draw [fill] (b) circle [radius=0.07];
\draw [fill] (c) circle [radius=0.07];
\draw [fill] (d) circle [radius=0.07];
\draw [fill] (e) circle [radius=0.07];
\draw [fill] (f) circle [radius=0.07];
\draw [fill] (g) circle [radius=0.07];
\draw [fill] (h) circle [radius=0.07];
\draw [fill] (i) circle [radius=0.07];
\draw [fill] (j) circle [radius=0.07];
\end{tikzpicture}
\end{center}
\end{center}
\caption{The Rips filtration (Definition \ref{def:Rips}).}
\label{Fig:RipsEx}
\end{figure}

In view of computational considerations, Rips filtrations are a natural construction when working with high-dimensional or non-Euclidean data.  \cref{Fig:RipsEx} illustrates the Rips filtration of a set of points in the plane.

By applying the $i$-th homology functor with coefficients in the field $k$ to either an offset or Rips filtration, one obtains a $[0,\infty)$-persistence module.  The following theorem tells us that this module has a well-defined barcode.

\begin{theorem}[Structure of Persistence Modules \cite{botnan2020decomposition,crawley2015decomposition}]
\label{Thm:ordinary_Structure}
If $M$ is p.f.d. persistence module indexed by a totally ordered set $T$, then there exists a unique multiset $\B{M}$ of intervals in $T$, such that 
\[M\cong \bigoplus_{I\in \B{M}} k_I.\]  
\end{theorem}

We call $\B{M}$ the \emph{barcode of $M$}.  For $\F$ a filtration such that $H_i(\F)$ is p.f.d., we write $\B{H_i(\F)}$ simply as $\mathcal{B}_i(\F)$.  

\begin{remark}[History]
\cref{Thm:ordinary_Structure} was proven in the case $T=\Z$ by Webb in 1985 \cite{webb1985decomposition}, for all $T$ with a countable dense subset by Crawley-Boevey in 2012 \cite{crawley2015decomposition}, and for all totally ordered sets $T$ by Botnan and Crawley-Boevey in 2018 \cite{botnan2020decomposition}. The important special case of finitely presented modules \cite{zomorodian2005computing} is a slight variant of the standard structure theorem for finitely generated modules over a PID, which can be found in many undergraduate abstract algebra textbooks. 
\end{remark}

\paragraph{Computation}
The barcode of a filtered simplicial complex containing $n$ simplices can be computed in $O(n^3)$ using a variant of Gaussian elimination \cite{edelsbrunner2002topological}. There has been a great deal of work in the past two decades on optimizing this basic approach, leading to dramatic performance gains for important classes of filtrations \cite{bauer2021ripser,otter2017roadmap}.

\paragraph{Applications}
Persistent homology has been applied in many areas, including computational chemistry \cite{nguyen2019mathematical}, materials science \cite{lee2017quantifying,hiraoka2016hierarchical}, neuroscience \cite{sizemore2018cliques,giusti2016two,rybakken2019decoding, gardner2022toroidal}, and bioinformatics \cite{rabadan2019topological}.  It is useful for exploratory data analysis \cite{hiraoka2016hierarchical,chan2013topology,carlsson2008local}, and is commonly integrated into pipelines for supervised learning \cite{hensel2021survey}.

\paragraph{Vectorization of Barcodes}
In applications to machine learning or statistics, one often wants to work with features of data taking values in a vector space, since many standard machine learning methods and statistical tests require their input to be of this form.  Yet the space of barcodes does not naturally have the structure of a vector space, and indeed its geometry is rather complicated \cite[Theorem 2.5]{turner2014frechet}.   A simple and practical solution is to fix a map from barcode space to a vector space, and work with the images of barcodes under this map.  Many such maps have been proposed and applied to problems in machine learning.  See \cite{hensel2021survey} for an overview.

\subsection{The Bottleneck Distance}
The \emph{bottleneck distance} $d_B$ is the standard metric on barcodes in the persistence theory. Intuitively, for barcodes $\mathcal{C}$ and $\mathcal{D}$, $d_B(\mathcal{C},\mathcal{D})$ is the maximum distance we must perturb an endpoint on of an interval in $\mathcal{C}$ to transform $\mathcal{C}$ into $\mathcal{D}$.  While $d_B$ can be defined between any two  barcodes whose intervals lie in $\R$ \cite{bauer2016persistence}, for simplicity's sake we give the definition only for barcodes with finitely many intervals, each of the form $[a,b)$, where $a<b\in \R$.

Let $\mathcal{C}$ and $\mathcal{D}$ be two such barcodes.  We define a \emph{matching between $\mathcal{C}$ and  $\mathcal{D}$} to be a collection of pairs $\chi = \{(I,J)\in \mathcal{C}\times \mathcal{D}\}$ where each $I$ and each $J$ appears in at most one pair.  
 If $I\in \mathcal \mathcal{C}\cup \mathcal{D}$ does not appear in a pair, then we say that $I$ is \emph{unmatched}. 

Given intervals $I=[a,b)$ and $J=[c,d)$, let
\[c(I,J) = \max(|c-a|,|d-b|),\qquad c(I) = (b-a)/2.\]
For a matching $\chi$ between $\mathcal{C}$ and $\mathcal{D}$, we define 
\begin{equation}
\cost(\chi):=\max\left( \max_{(I,J)\in \chi} c(I,J),\ \max_{ I \in \mathcal{C}\cup \mathcal{D} \text{ unmatched}}\ c(I)\right).
\label{eq:cost-matching}
\end{equation}

\begin{definition}
\label{Def:bottleneck}
The \emph{bottleneck distance} between $\mathcal{C}$ and $\mathcal{D}$ is 
\[d_B(\mathcal{C}, \mathcal{D}) = \inf\, \{ \cost(\chi): \text{ $\chi$ is a matching between $\mathcal{C}$ and $\mathcal{D}$}\}.\]
\end{definition}

\begin{remark}
The bottleneck distance $d_B(\mathcal{C}, \mathcal{D})$ can be computed in time $O(n^{1.5}\log n)$ where $n=|\mathcal{C}| + |\mathcal{D}|$ \cite{kerber2017geometry,efrat2001geometry}.  An efficient implementation is available in the Hera software \cite{kerber2017geometry}.
\end{remark}

\subsection{Stability of Persistent Homology}
The \emph{stability theorem} for persistent homology asserts that the map from data to barcodes is 1-Lipschitz continuous with respect to suitable choices of metrics.  This plays an important role in both the statistical foundations of PH \cite{fasy2014confidence} and the computational theory \cite{sheehy2013linear}. And as we will see, it is also the starting point for some of the central ideas in multiparameter persistence.  In fact, there are many variants of the stability theorem in the literature \cite{cohen2007stability,chazal2012structure,bauer2015induced,botnan2018algebraic,bjerkevik2021ell}.  Here we will state a version of stability for sublevel filtrations, as well as corollaries of this concerning the offset and Rips filtrations.  In \cref{Sec:Interleavings}, we will also discuss a more general algebraic formulation of stability called the \emph{isometry theorem}.

\begin{definition}\label{Def:Distances_on_Point_Clouds}
~
\begin{itemize}
\item[(i)] Given $X,Y\subset \R^n$, the \emph{Hausdorff distance} between $X$ and $Y$ is given by \[d_H(X,Y)=\inf\,\{r\geq 0\mid X\subset \cup_{p\in Y} B(p,r),\ Y\subset \cup_{p\in X} B(p,r)\},\]
where $B(p,r)$ denotes the closed ball of radius $r$ centered at $p$.
\item[(ii)]
Given metric spaces $X,Y$, the\emph{ Gromov-Hausdorff} distance between $X$ and $Y$ is given by \[d_{GH}(X,Y)=\inf_{\gamma,\kappa} d_H(\gamma(X),\kappa(Y)),\]
where $\gamma:X\to Z$ and $\kappa:Y\to Z$ range over all isometric embeddings of $X$ and $Y$ into a common metric space $Z$. 
\end{itemize}
\end{definition}

\begin{theorem}[Stability of Persistent Homology \cite{cohen2007stability,chazal2009proximity,chazal2009gromov,chazal2014persistence}]
\label{Thm:Stability}
~
\begin{itemize}
\item[(i)] For any topological space $W$, functions $\gamma,\kappa:W\to \R$, and $i\in \N$ such that $H_i(S^\uparrow(\gamma))$ and $H_i(S^\uparrow(\kappa))$ are p.f.d., we have 
\[d_B(\mathcal{B}_i(S^\uparrow(\gamma)),\mathcal{B}_i(S^\uparrow(\kappa)))\leq 
 \sup_{w\in W} |\gamma(w) - \kappa(w)|.\]
\item[(ii)] For finite $X,Y\subset \R^n$ and $i\in \mathbb N$,
\[d_B(\mathcal B_i(\Offset(X)),\mathcal B_i(\Offset(Y)))\leq d_{H}(X,Y).\]
\item[(iii)] For finite metric spaces $X,Y$ and $i\in \mathbb N$,
\[d_B(\mathcal B_i(\Rips(X)),\mathcal B_i(\Rips(Y)))\leq d_{GH}(X,Y).\]
\end{itemize}
\end{theorem}

We will see multiparameter analogues of each of these three stability statements in \cref{sec:metrics}.

\subsection{Ordinary Persistent Homology is Not Robust}\label{Sec:Not_Robust}
\begin{figure}[h!]
\centering
\begin{subfigure}{.45\textwidth}
  \centering
  \includegraphics[width=.9\linewidth]{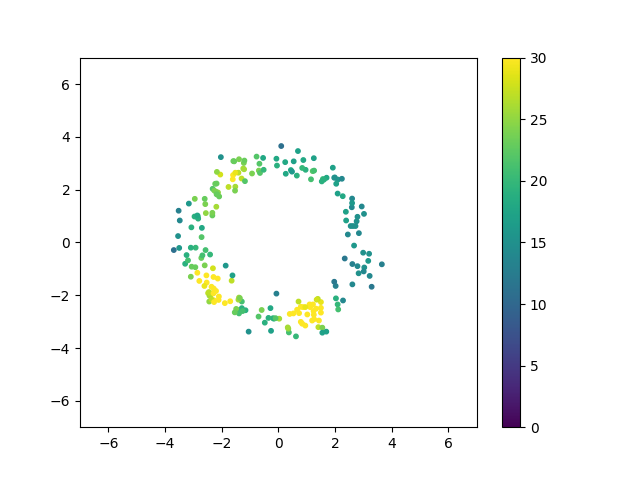}
  \caption{}
  \label{fig:circle}
\end{subfigure}%
\begin{subfigure}{.45\textwidth}
  \centering
  \includegraphics[width=.9\linewidth]{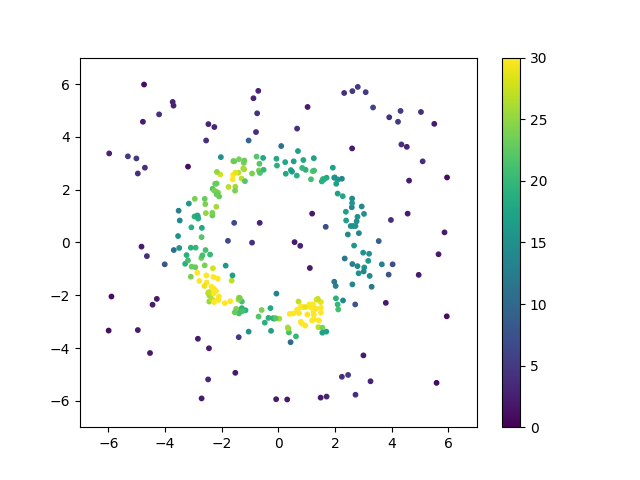}
  \caption{}
  \label{fig:circle-all-points}
\end{subfigure}
\begin{subfigure}{.45\textwidth}
  \centering
  \includegraphics[width=.9\linewidth]{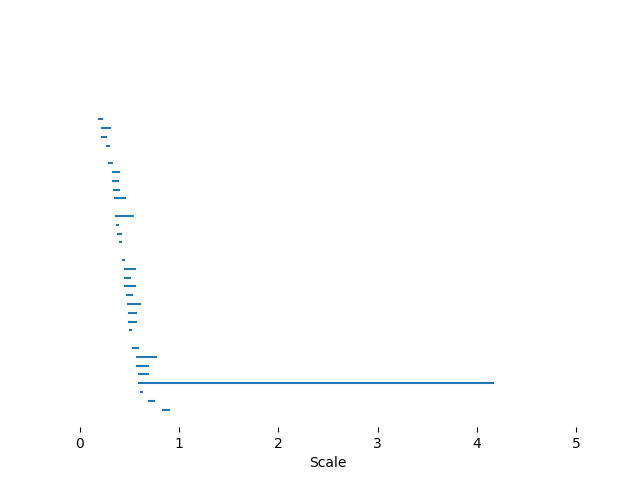}
  \caption{}
  \label{fig:barcode-circle}
\end{subfigure}%
\begin{subfigure}{.45\textwidth}
  \centering
  \includegraphics[width=.9\linewidth]{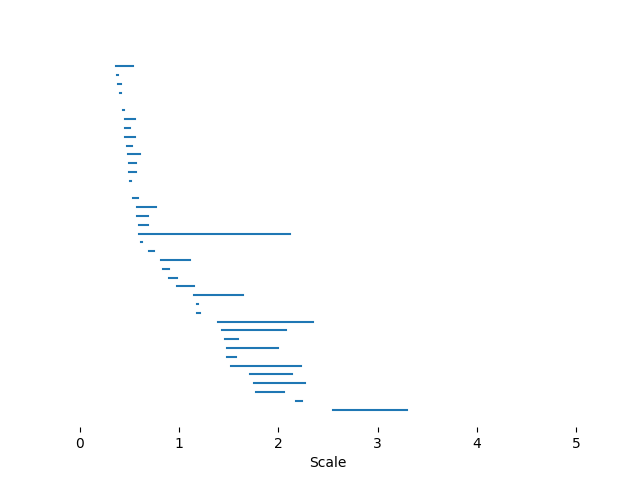}
  \caption{}
  \label{fig:barcode-all}
\end{subfigure}
\caption{\textbf{(a)} A data set with a circular shape colorized by a local density estimate. \textbf{(b)} The data from (a) with added noise and colorized by a local density estimate.  \textbf{(c)} The barcode of the data in (a). \textbf{(d)} The barcode of the data in (b). The barcodes were computed using Ripser \cite{bauer2021ripser}. }
\label{fig:circles-full}
\end{figure}

\cref{Thm:Stability}\,(ii) and (iii) notwithstanding, both the offset and Rips constructions of persistent homology are highly unstable to outliers \cite[Section  4]{blumberg2014robust}, and both constructions can be insensitive to geometric structure in high-density regions of the data.  For example, consider the two data sets shown in \cref{fig:circles-full}, and note that the second differs from the first only by the addition of a few points.  While the $H_1$ barcode of the first data set has a long interval, representing a loop in the data, this signal is lost in the barcode of the second data set.

There have been many proposals for dealing with such issues within the framework of 1-parameter persistence; see \cite{blumberg2020stability} for a recent overview. However, these approaches share important disadvantages, which we now explain, closely following \cite{blumberg2020stability}.  First, they require the user to fix one or more parameters, typically either a spatial scale or a density threshold.  The suitability of a particular parameter choice depends on the data, and it can be unclear how to best choose such a parameter. Moreover, for data with interesting features across a range of scales and densities, there may be no single parameter choice that fully captures the structure of interest in the data.  In addition, strategies that fix a scale parameter do not distinguish small features in the data from large features, and strategies that fix a density threshold do not distinguish features appearing at high densities from those appearing at low densities.

A natural solution to these problems is to consider 2-parameter persistent homology, where one of the parameters is a scale parameter, as in the
Rips filtration, and the other parameter is a density threshold.  In fact, as we will see in \cref{Sec:Density-Sensitive Multifiltrations}, there are several natural constructions in this spirit with different advantages and disadvantages.  

\section{Multiparameter Persistent Homology}
\label{sec:multi-d-pipe}
To develop multiparameter persistent homology, one would like to generalize the persistent homology pipeline of \cref{sec:1D-pipeline}.  The first two arrows of the pipeline generalize without difficulty, as follows: 

\medskip
\begin{center}
\framebox[1.1\width]{Data} $\Rightarrow$ \framebox[1.1\width]{Multifiltration} $\xRightarrow{\text{Homology}}$ \framebox[1.1\width]{Multiparameter Persistence Module}
\end{center}
\medskip

Indeed, applying homology to each space and each map in a multifiltration yields a multiparameter persistence module, exactly as in the 1-parameter case. Furthermore, there are many  natural ways of constructing multifiltrations from data; several of these will be discussed in \cref{Sec:Multifiltrations}.  As one example, we have the following density-sensitive extension of the Rips filtration: 
\begin{definition}\label{Def:Degree_Rips}
For $X$ a metric space, $r\geq 0$, and $d>0$, let $\DRips(X)_{d,r}$ be the maximal subcomplex of $\Rips(X)_r$ whose vertices have degree at least $d-1$ in the 1-skeleton of $\Rips(X)_r$.  Varying $r$ and $d$, we obtain a bifiltration $\DRips(X)$, the \emph{degree-Rips bifiltration} \cite{lesnick2015interactive}.  
\end{definition}
\cref{fig:Degree-Rips} illustrates part of the degree-Rips bifiltration of the example from \cref{Fig:RipsEx}.  

\subsection{The Difficulty of Defining Barcodes of Multiparameter Persistence Modules} \label{sec:diff-barc}
We now turn to the salient question: Is there any good way to define the barcode of a multiparameter persistence module?  As we will now see, there are several illuminating ways one can approach this question.

We begin with some good news: As in the 1-parameter case, p.f.d. multiparameter persistence modules decompose into indecomposable summands in an essentially unique way.  In fact, this is true for persistence modules indexed by any poset:  

\label{sec:decomp}
\begin{theorem}\label{Thm:KS}
For a poset $P$ and persistence module $M\colon P\to \Vect$, 
\begin{itemize}
\item[(i)] If $M$ is pointwise finite-dimensional, then there exists a collection of indecomposables $\{M^\lambda\}_{\lambda \in \Lambda}$ such that \[M\cong \bigoplus_{\lambda \in \Lambda} M^\lambda.\]  
\item[(ii)] If \[M\cong \bigoplus_{\lambda \in \Lambda} M^\lambda \cong\bigoplus_{\gamma\in \Gamma} M^\gamma\]
with each $M^\lambda$ and $M^\gamma$ indecomposable, then there exists a bijection $\sigma\colon  \Lambda\to \Gamma$ such that
$M^\lambda \cong M^{\sigma(\lambda)}$ for all $\lambda\in \Lambda$. 
\end{itemize}
\end{theorem}
In view of this result, understanding the algebraic structure of p.f.d. persistence modules boils down to understanding the structure of the indecomposable modules.  

\begin{remark}
For finite posets, the proof of \cref{Thm:KS}\,(i) is elementary and can be found in many introductory texts on quiver representations; see e.g.  \cite[Chapter 4.1]{barot2015introduction}. In the generality given here, \cref{Thm:KS}\,(i) was sketched in work of Gabriel and Roiter in 1992 \cite{gabriel1992representations}, follows from work of Crawley-Boevey from 1994 \cite{crawley1994locally}, and was given a short direct proof in recent work of Crawley-Boevey and Botnan from 2018 \cite{botnan2020decomposition}. \cref{Thm:KS}\,(ii) is the Azumaya--Krull--Remak--Schmidt theorem \cite{azumaya1950corrections}.
\end{remark}

\begin{definition}
A $P$-module $M$ is \emph{interval-decomposable} if there exists a multiset  $\B{M}$  of intervals in $P$ such that 
\[ M \cong \bigoplus_{I\in \B{M}} k_I.\]
We call $\B{M}$ the \emph{barcode} of $M$. 
\end{definition} 

By virtue of \cref{Thm:KS}\,(ii), the barcode of an interval-decomposable is well-defined whenever it exists.  However, it turns out that for $n\geq 2$, not all $n$-parameter persistence modules are interval-decomposable; drawing on classical ideas from quiver representation theory, we shall see in \cref{sec:multidbarcode} that the space of indecomposables that can arise in the multiparameter setting is enormously complex.  In particular, there is no way of parameterizing this space by collections of nice regions in the parameter space, as in the 1-parameter setting. Thus, while one could define the barcode of a multiparameter parameter persistence module to be the multiset of isomorphism types of its incomposable summands, this object is generally too complex to work with in practice. 

One might nevertheless hope that there is a good way to define the barcode of a multiparameter persistence module $M:P\to \Vect$ as a collection of regions in $P$.  But it turns out not to be possible, in the following sense.

\begin{definition}\label{Def:Good_Barcode}
A multiset $B$ of subsets of $P$ is a \emph{good barcode of $M$} if for all $x\leq y\in P$ we have \[\Rk (M_x\to M_y) = |\{S\in B : x,y\in S\}|,\] i.e., the rank of the map $M_x\to M_y$ is the number of elements of $B$ containing both $x$ and $y$.  
\end{definition}

Given how barcodes of 1-parameter persistence modules are usually interpreted and used in TDA, the goodness condition of Definition  \ref{Def:Good_Barcode} is quite natural. However, the next example shows that a good barcode of $M$ need not exist.

\begin{example}\label{Ex:Not_Good}
Let  $P= \{0,1,2\}\times \{0,1,2\}$ and let $M$ be the following $P$-module, where $=$ denotes identity maps:
\begin{equation}\begin{tikzcd}
k\arrow{r}{=} & k\ar[r] & 0 \\
k\arrow{r}{[1,0]^T} \arrow{u}{=} & k^2\arrow{r}{[1,1]}\arrow{u}{[1,0]} & k\ar[u]\\
0\ar[r]\ar[u] & k\arrow{u}{[0,1]^T}\arrow{r}{=} & k\ar[u,"="]
\end{tikzcd}
\label{eq:not-rank2}
\end{equation}
A simple argument by contradiction shows that $M$ cannot have a good barcode: If $B$ is a good barcode of $M$, then since \[\Rk (M_{(0,1)} \to M_{(2,1)})=\Rk (M_{(0,1)} \to M_{(1,2)})=\Rk (M_{(1,0)} \to M_{(2,1)})=1,\] $B$ must contain intervals $I,J,K$ (not necessarily distinct) with 
\[(0,1),(2,1)\in I,\qquad (0,1),(1,2)\in J,\qquad (1,0),(2,1)\in K.\]  But since $\dim M_{0,1}=\dim M_{2,1}=1$, we then have $I=J=K$, implying that $(1,0),(1,2)\in I$.  This contradicts the fact that $\Rk (M_{(1,0)} \to M_{(1,2)}) = 0$.
\end{example}
Perhaps surprisingly, recent work \cite{botnan2021signed} has shown that if we allow the elements of the barcode to be \emph{signed} subsets of $P$ (i.e., to be labeled positive or negative) then it is possible to give a well-behaved definition of the barcode, which encodes ranks in a way analogous to Definition \ref{Def:Good_Barcode}.  We shall discuss this in detail in \cref{sec:rankinv}.

\subsection{Invariants}
\label{subsec:invariants}
Invariants of multiparameter persistence modules play a role analogous to barcodes in 1-parameter persistence: They can be used to visualize persistence modules and can be fed as input to machine learning methods and statistical tests.  Many invariants of persistence modules have been proposed in the TDA literature, and one can find yet more in the classical literature on commutative algebra and representation theory.  The main question for TDA is which such invariants can be useful in the development of data analysis methodology.  We are still in the early stages of understanding this.  

In this section, we will introduce several simple and well-known invariants of multiparameter persistence modules.  Signed barcodes, a class of potentially useful invariants which are equivalent to some of the invariants discussed in this section, will be introduced later, in \cref{sec:rankinv}.

\paragraph{Three Simple Invariants}
We begin by introducing three simple invariants which arise frequently in the MPH literature.  These are a good starting point for thinking about how to do practical work with MPH.  Let $P$ be a poset.  

\begin{enumerate} 
\item The \emph{Hilbert function} of a p.f.d. persistence module $M\colon P\to \Vect$ is the function $P\to \mathbb N$ sending $z$ to $\dim M_z$.
\item The \emph{rank invariant} of $M$ \cite{carlsson2009theory} records the rank of the linear map $M_a \to M_b$ for every $a\leq b\in P$.  Clearly, this is a refinement of the Hilbert function.  When $P=\N$ or $\R$, the rank invariant is \emph{complete}, i.e., it determines the isomorphism type of $M$; see \cite[Theorem 12]{carlsson2009theory} and \cite{crawley2012decomposition}.  However, for $P=\N^2$, the rank invariant is incomplete.

\item For $M\colon \R^n\to \Vect$ finitely presented (see \cref{SubSec:Pres_Res}) and $z\in \R^n$, the \emph{(multi-graded) Betti numbers of $M$} at $z$ are natural numbers \[\beta^M_0(z),\beta^M_1(z),\ldots,\beta^M_n(z).\]
Informally, $\beta^M_0(z)$ and $\beta^M_1(z)$ count the number of generators and relations in $M$ at $z$, respectively, while for $i>1$, $\beta^M_i(z)$ counts higher-order relations. 
 More formally, given a minimal free resolution of $M$ 
\[0\to F_n \to  \cdots \to F_1\to F_0\to M\to 0\]
(see Definition \ref{def:minimal-free}), $\beta^M_i(z)$ is the number of elements at bigrade $r$ in a basis for $F_i$ \cite{eisenbud2005geometry,lesnick2015interactive}.  These are standard invariants in commutative algebra \cite{miller2004combinatorial}.
\end{enumerate}

\paragraph{The Fibered Barcode}
We say an affine line $L\subset \R^n$ is \emph{admissible} if the product partial order on $\R^n$ restricts to a total order on $L$.  Note that the admissible lines in $\R^2$ are those with non-negative slope. The \emph{fibered barcode} of an $\R^n$-indexed module $M$ \cite{cerri2013betti} is the function which maps each admissible line $L$ to the barcode $\B{M\circ L}$.  It is equivalent to the rank invariant of $M$, but more convenient for visualization and stability analysis; see \cref{subsec:visual,Sec:Matching_Dist}.  

\paragraph{The Generalized Rank Invariant}
The following extension of the rank invariant was proposed in \cite{kim2018generalized}.

\begin{definition}\label{def:generalized-rank_invariant}
Let $M\colon P\to \Vect$ be p.f.d.  Given an interval $\int\subseteq P$, the {\em generalized rank} of $M$ over $\int$, denoted by $\Rank_\int M$, is defined by:
\[ \Rank_\int M = \rank\left[ \varprojlim M|_\int \to \varinjlim M|_\int \right]. \]
Given a collection $\Int$ of intervals, the {\em generalized rank invariant} of $M$ over $\Int$ is the function $\rank_\Int M \colon \Int \to \N \cup \{ \infty \}$ defined by $\rank_\Int M(\int) =  \Rk_\int M$.
\end{definition}

\begin{remark}
Kim and Memoli~\cite{kim2018generalized} proved that for interval-decomposable modules indexed by \emph{essentially finite} posets, the generalized rank invariant is complete.  This result was later extended slightly in \cite{botnan2021signed}.   Moreover, Kim and Moore~\cite{kim2021bigraded} showed that the generalized rank invariant determines the Betti numbers of $M$. 
\end{remark}

\paragraph{Vectorizations of Persistence Modules}
As noted in \cref{sec:1D-pipeline}, many applications of 1-parameter persistent homology involve mapping barcodes into linear spaces and then applying standard statistics and machine learning methods.  It is natural to pursue the same idea in the multiparameter setting.  To this end, several novel maps from the space of persistence modules into linear spaces have been proposed and applied to data \cite{vipond2018multiparameter,carriere2020multiparameter,corbet2019kernel,vipond2021multiparameter,beltramo2021euler}.   Most such maps proposed so far depend only on the rank invariant, and several depend only on a part of the rank invariant.  

\paragraph{Invariants from Metrics}
A choice of metric $d$ on persistence modules can be used to define invariants.  For example, fixing a persistence module $N$, we may consider the map $M\mapsto d(M,N)$.  

For a more interesting example, suppose we are given an $\R$-valued invariant $f$ of multiparameter persistence modules; for example, we may take $f(M)$ to be the sum of all $0^{\mathrm{th}}$ Betti numbers of $M$.  Scolamiero et. al \cite{scolamiero2017multidimensional} define an invariant $\hat f$, called the \emph{hierarchical stabilization} of $f$, taking values in the set of decreasing functions from $[0,\infty)$ to $\R$.   Specifically, they define \[\hat f(M)(r)\coloneqq\inf\, \{f(N)\colon  d(N,M)\leq r\}.\]  It is shown that the map $M\mapsto \hat f(M)$ is 1-Lipchitz continuous with respect to $d$ and an interleaving metric on the functions.  This construction in fact extends immediately from persistence modules to any metric space; a variant was previously considered in the context of mode detection for $\R$-valued functions \cite{bauer2017persistence}. 
While mathematically attractive, computation of hierarchical stabilization in the multiparameter setting can be hard \cite{gafvert2017stable}.

\paragraph{Other Work on Invariants of MPH}
We briefly mention a few other works which consider invariants of multiparameter persistence modules in the context of TDA; we make no attempt to be exhaustive.  Harrington et al. \cite{harrington2019stratifying} study several classical invariants from commutative algebra in the multigraded setting, namely, Hilbert series, associated primes, and $0^{\mathrm{th}}$ local cohomology. Extending an idea of Patel in the 1-parameter setting \cite{patel2018generalized}, several papers \cite{mccleary2020edit, asashiba2019approximation,kim2018generalized} have explored the idea of defining signed barcodes in the multiparameter setting via \emph{M\" obius inversion} \cite{rota1964foundations}.  This work sets the stage for the results of \cite{botnan2021signed} discussed above on signed decompositions of the rank invariant; see Remark \ref{Rem:Other_Signed_Barcodes}.  
  Cai et al. \cite{cai2020elder} introduce a well-behaved and computable barcode, called the \emph{staircase code}, for the $0^{\mathrm{th}}$ homology modules of superlevel-Rips bifiltrations (Definition \ref{def:superlevel-rips}).
\subsection{Visualization}
\label{subsec:visual}
\begin{figure}[h!] 
\begin{subfigure}{.4\textwidth}
\centering
\definecolor{aqaqaq}{rgb}{0,0.62,0.62}
\input{Rips_Ex_Tikz_2_Params}
\caption{}
\label{fig:Degree-Rips}
\end{subfigure}
\begin{subfigure}{.55\textwidth}
\centering
\includegraphics[scale=.14]{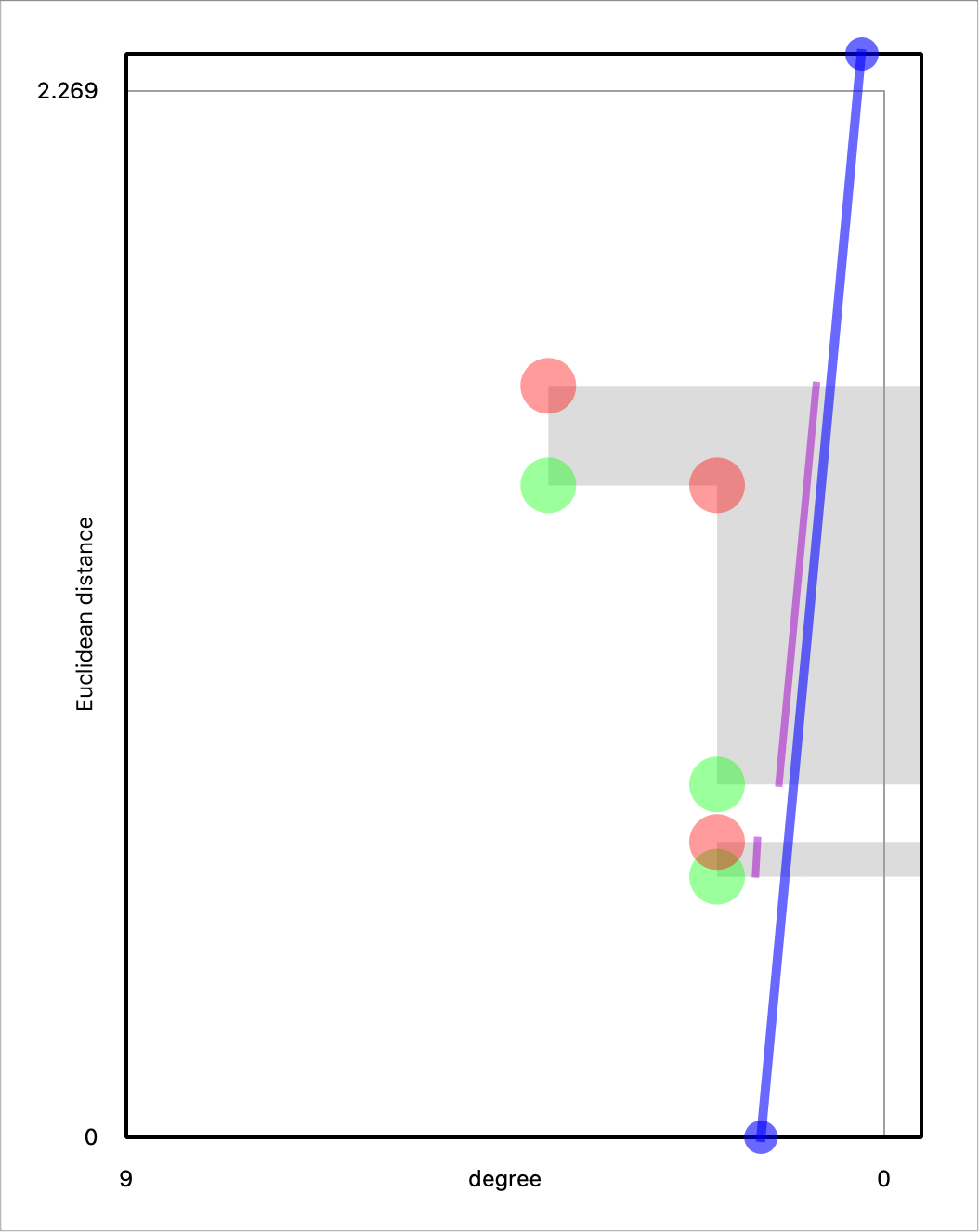}
\caption{}
\label{fig:degree-Rips_small}
\end{subfigure}
\caption{(a) The Degree-Rips Bifiltration of the same data set considered in \cref{Fig:RipsEx}, for a few choices of the degree and scale parameter: The left, middle, and right columns correspond to the degree parameters 4, 2, and 0, respectively. (b) The corresponding RIVET visualization in homology dimension 1. }
\end{figure}

\begin{figure}[h!]
\centering
\includegraphics[scale=.12]{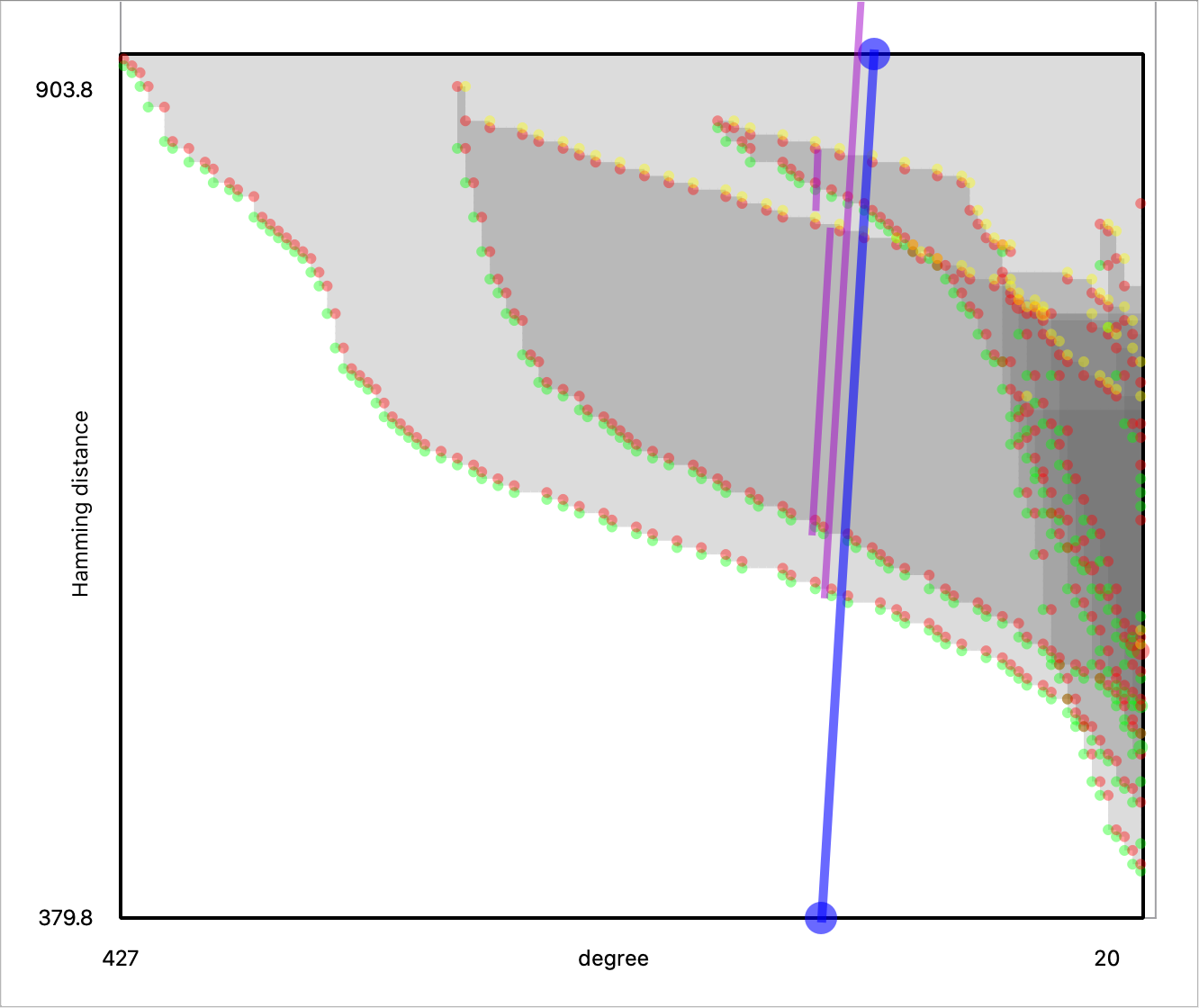}
\hskip10pt
\includegraphics[scale=.12]{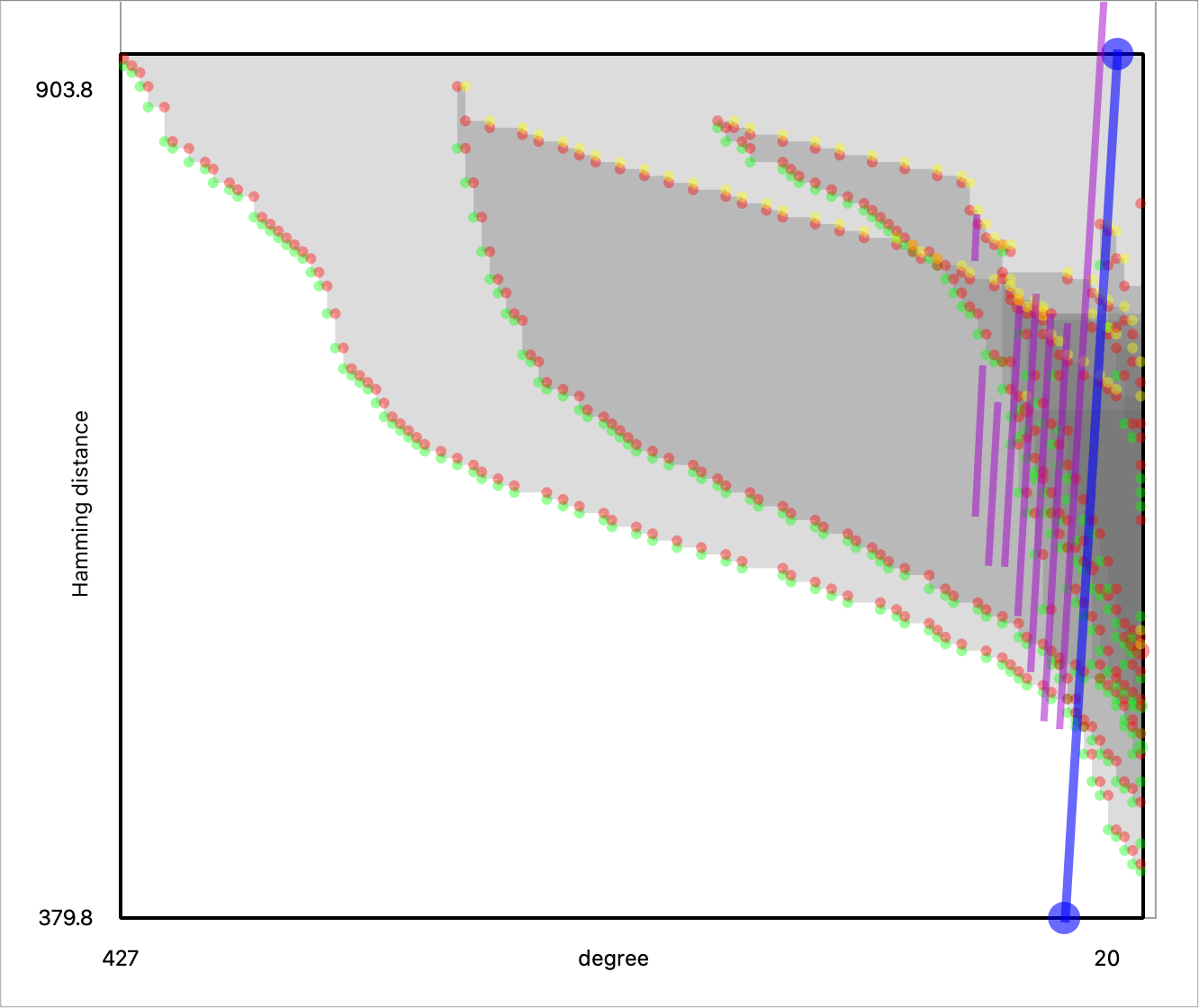}
\caption{RIVET's visualization of the 0th degree-Rips persistent homology of a data set of 1088 HIV-1 genomes, for two different choices of the line $L$ (shown in blue).  The visualization indicates the presence of two major clusters in the data (large grey regions to the left), each with several hundred points, as well as 5 smaller clusters of less than 40 points (darker grey regions to the right).  Beyond this, the plots of the Hilbert functions and bigraded Betti numbers exhibit interesting geometry which encodes subtle information about the size and shape of the clusters.}
\label{fig:HIV}
\end{figure}

We visualize the Hilbert function, fibered barcode, and bigraded Betti numbers in two examples using RIVET \cite{lesnick2015interactive,rivet}. Figures \ref{fig:degree-Rips_small} and \ref{fig:HIV} show (part of) RIVET's visualization of degree-Rips bifiltrations: As a simple example, \cref{fig:degree-Rips_small} visualizes the 1st PH (i.e., loop structure) of the degree-Rips bifiltration of the data of \cref{Fig:RipsEx}.  \cref{fig:HIV}  visualizes the 0th PH (i.e., cluster structure) of the degree-Rips bifiltration of 1088 pre-aligned HIV-1 genomes from \cite{chan2013topology}, metrized using the Hamming distance.\footnote{To control the size of this bifiltration, RIVET coarsens it slightly so that all simplices are born on a $250\times 250$ grid.  This coarsening is \emph{stable}, i.e., it changes the modules only slightly in the interleaving distance.}  

To explain the figures, first note that the $x$-axis is mirrored in each figure so that values decrease from left to right.  The Hilbert function is represented by greyscale shading:  In each figure, the darkness of shading is proportional to the vector space dimension, and the lightest non-white shade represents a value of $1$.  The bigraded Betti numbers are represented by translucent colored dots whose area is proportional to the value; the $0^{\mathrm{th}}$, $1^{\mathrm{st}}$, and $2^{\mathrm{nd}}$ Betti numbers are shown in green, red, and yellow.  For the fibered barcode visualization, the query line $L$ is shown in blue, and the barcode is shown in purple, with each interval offset perpendicularly from $L$.

Like other viruses, HIV has a rich and epidemiologically important subtype structure \cite{hemelaar2019global}; \cref{fig:HIV} indicates that the degree-Rips PH is able to see key aspects of this structure, without any data preprocessing or parameter choices that may bias the results.  In contrast, the 1-parameter Rips PH of this data (not shown) sees no cluster structure, because of the presence of low density outliers between the clusters.
\cite{eisenbud2005geometry,lesnick2015interactive}.

\begin{remark}\label{Rem:Aug_Arr_RIVET}
One of the main features of RIVET's visualization is an interactive scheme for visualizing the fibered barcode: The user selects the line $L$ by clicking and dragging the mouse, and the display of the barcode $\mathcal{B}(M\circ L)$ updates in realtime.  To support this real-time interactivity, RIVET precomputes a data structure called the \emph{augmented arrangement,} which can be queried for the barcode $\mathcal{B}(M\circ L)$ along a generic line $L$ in time $O(|\mathcal{B}(M\circ L)|+\log n)$, where $n$ is the size of a grid containing the supports $M$ \cite{lesnick2015interactive}.
\end{remark}

\begin{remark}
A different approach to visualizing the rank invariant arises from the idea of signed barcodes, discussed in \cref{sec:rankinv}.
\end{remark}

\subsection{Applications}\label{Sec:Applications}
We briefly discuss two recent applications of MPH to cancer imaging.

\paragraph{Spatial Patterns of Immune Cells in Tumors}
Vipond et al. \cite{vipond2021multiparameter} use a vectorization of MPH called the \emph{multiparameter persistence landscape} \cite{vipond2018multiparameter} to analyze the spatial patterns of immune cells in cancerous tumors.  This work uses RIVET's computational backend, as well as Vipond's code for multiparameter persistence landscapes \cite{vipondGitHub}.

\begin{figure}
\centering
\includegraphics[scale=1]{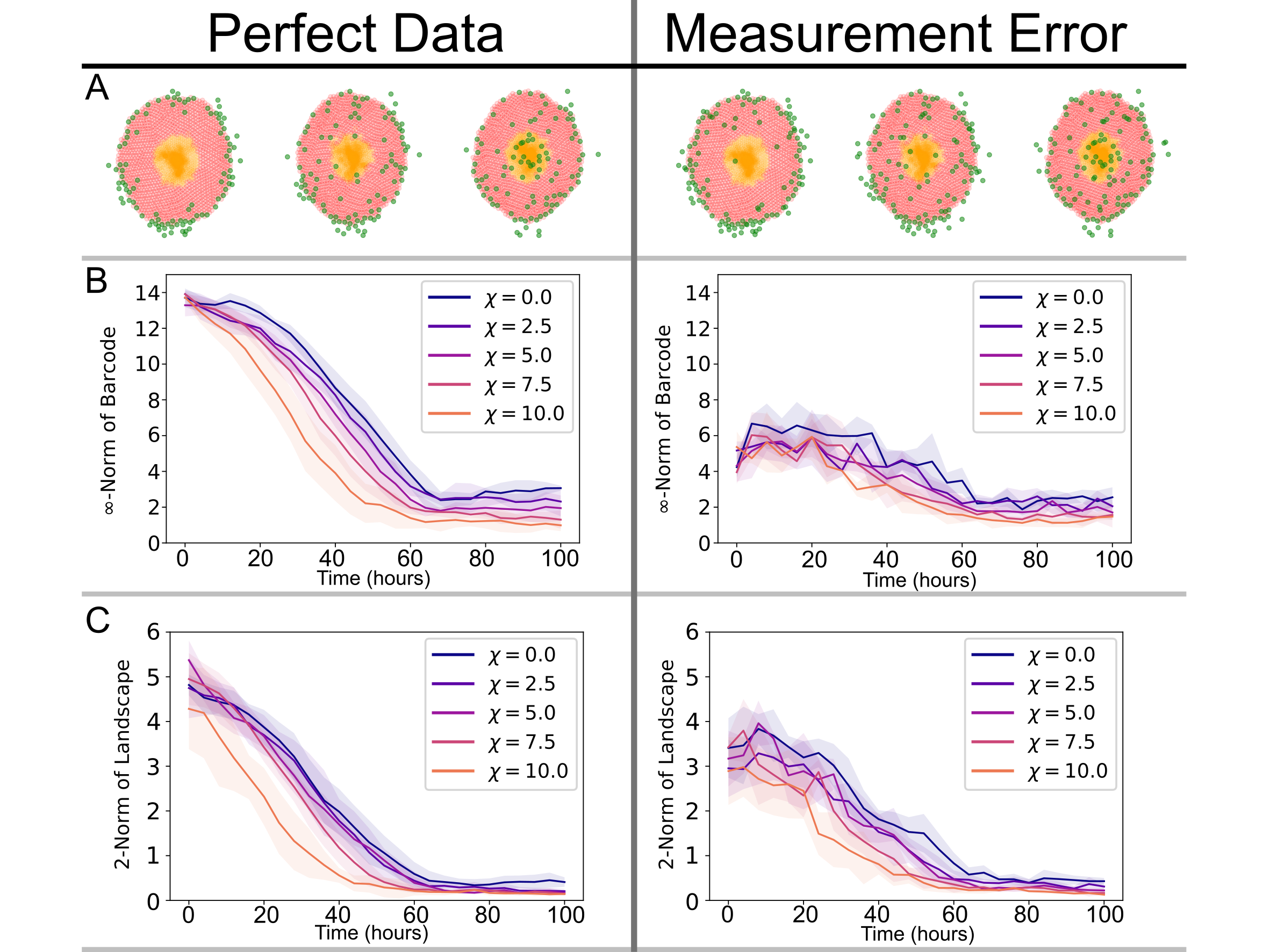}
\caption{(A) Cell distributions at three time points from a simulation of immune cells infiltrating a tumor.  Green, red, and orange represent macrophages, viable tumor cells, and necrotic tumor cells, respectively. There are ${\sim}100$ macrophages and ${\sim}20$ tumor cells falsely registered as macrophages. (B) Decay curves tracking the length of the longest bar in the $H_1$ barcode against time for simulations with different chemotaxis parameter $\chi$. Average curves with standard deviation  bands are shown, computed with five simulations for each value of $\chi$. (C) MPH decay curves tracking the 2-norm of the MPH landscape against time. An extended version of this figure can be found in \cite[Fig. 1]{vipond2021multiparameter}. }
\label{fig:vipond}
\end{figure}  
 
Information about the spatial distribution of immune cells in and around a tumor is useful both for prognosis and for guiding treatment \cite{fridman2017immune,galon2006type,kather2018topography,nawaz2015beyond,bull2020combining}.  The extent to which immune cells infiltrate the tumor is of particular importance.  In \cite{vipond2021multiparameter}, the positions of the immune cells are determined from \emph{chromogenic immunohistochemistry} images of the cancer tissue via a semiautomated method.

To quantify the spatial patterns of immune cells, prior approaches have compared the  densities of immune cells at the boundary and inner core of the tumor, or used established spatial statistics.  However, these approaches capture limited information about the spatial patterns of the cells; one hopes that the geometric information provided by TDA would be of clinical use. 
A key challenge is that the analysis of spatial distributions of immune cells is complicated by noise in the data: the classification of cells in a tumor image according to cell type is prone to errors, which is expected to cause problems for the standard 1-parameter constructions of persistent homology.   It is thus natural to consider MPH in this setting.

In \cite{vipond2021multiparameter}, Vipond et al. analyze both synthetic data and real histological data from head and neck tumors. The synthetic data is generated using a 2-dimensional  agent-based simulation of macrophages infiltrating a disc-shaped tumor.  Initially, the macrophages are distributed along the boundary of the tumor; see \cref{fig:vipond}\,(A).  As time evolves, the macrophages move inward towards the center of the tumor, driven by a gradient in the concentration of a chemical attractant.  The speed of movement is determined by a \emph{chemotaxis} parameter $\chi$.  The authors track the dynamics topologically, using Vietoris--Rips persistent homology on the point cloud determined by the locations of the macrophages: the long bar present in the $H_1$ barcode at time 0 gradually shrinks as time evolves. In order to quantify this, the authors compute a persistent homology \emph{decay curve}, which tracks the length of the longest bar in the barcode as a function of time.  It is observed that different choices of $\chi$ give rise to distinct decay curves, with larger values of $\chi$ yielding steeper curves. %
However, when a biologically realistic amount of noise is added to the simulation by misclassifying some tumor cells as macrophages, the decay curves change qualitatively, muddying the relationship between $\chi$ and the shape of the decay curve; see \cref{fig:vipond}\,(B).  Such noise causes similar problems for an analogue of the decay curve constructed using an established spatial statistic called the \emph{pair correlation function}.  
To rectify this, the authors use a \emph{codensity-Rips bifiltration}, a variant of the density-Rips bifiltration introduced in Example \ref{Ex:Density_Rips} below, where one filters by both radius and a nearest-neighbors codensity function. The resulting \emph{multiparameter decay curves}, constructed from the 2-norm of the $H_1$ multiparameter persistence landscape \cite{vipond2018multiparameter}, are observed to be more robust to noise than their 1-parameter counterparts; see \cref{fig:vipond}\,(C).  This finding serves as a proof of concept for the application of multiparameter persistence landscapes to immunohistological data. 

The authors next study the spatial patterns of three types of immune cells (namely, the T-cells CD8$^+$ and  FoxP3$^{+}$, and the macrophage CD68$^{+}$) in images of head and neck tumors.  16 tumors are considered in total, though the key analyses and findings concern smaller subsets of these tumors.  For each image and cell type, the cell locations determine a point cloud.  
The image is decomposed into regions of a suitable fixed size and a MPH landscape is computed for the restriction of the point cloud to each region.  This yields an ensemble of landscapes for each tumor and cell type.  Each landscape can be loosely interpreted as a measure of the extent to which the immune cells infiltrate the corresponding region; intuitively, regions which have few immune cells lead to features of large persistence.  
  
An analysis of the landscapes of five of the largest tumors provides insight into the patterns of tumor infiltration for each cell type, and into the relationship between these patterns.  For instance, the authors observe via PCA and linear discriminant analysis that the landscape ensembles of different cell types have qualitatively distinct distributions, reflecting differences in the extent and pattern of immune cell infiltration \cite[Figures 3, S10]{vipond2021multiparameter}.  In addition, statistics constructed from the landscapes indicate that FoxP3$^{+}$ cells usually infiltrate the tumors less than either CD8$^+$ or CD68$^{+}$ cells \cite[Tables S4, S7]{vipond2021multiparameter}.  Moreover, it is observed that a landscape-based measure of how FoxP3$^{+}$ cells infiltrate the tumor is correlated with oxygen levels of the tumor \cite[Supporting information, lines 75-85 and Tables S4, S7]{vipond2021multiparameter}.  
The sample sizes considered are too small to establish the statistical significance or biological reproducibility of these findings, but the results suggest that it may be fruitful to conduct similar analyses with a larger cohort.

While \cite{vipond2021multiparameter} does not make a direct connection between topological invariants and clinically important  variables such as disease prognosis, it demonstrates that MPH is a viable tool in the study of this type of data, capable of providing  biological insights.  

\begin{figure}[h!]
\center
\includegraphics[scale=0.25]{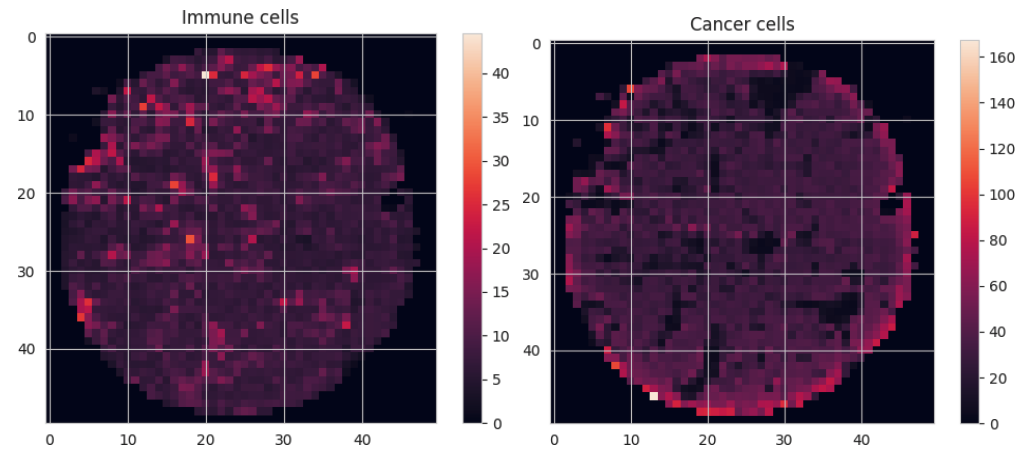}
\caption{Quantiative imunofluoresence images of a breast cancer tissue sample.  Left: pixel intensity corresponds to number of immune cells. Right: pixel intensity corresponds to number of cancer cells.  The figure is taken from \cite{carriere2020multiparameter}.}
\label{fig:immune}
\end{figure}

\paragraph{Classification of Breast Cancer Tissue}
As part of a larger investigation of feature maps for MPH in machine learning, Carri{\`e}re and Blumberg \cite{carriere2020multiparameter} use MPH in the study of \emph{quantitative immunofluorescence images} of breast cancer tissue.  For each tissue sample, one has two images with the same domain; in one, pixel intensity corresponds to the number of immune cells at that pixel and in the other, pixel intensity corresponds to the number of cancer cells; see \cref{fig:immune}. The objective is to train a classifier that predicts patient survival, using only these images.  

As explained in \cite{carriere2020multiparameter}, the usual approach to this problem begins by thresholding the intensities, in order to identify one subset of the pixels as ``immune pixels,'' and another subset as ``cancer pixels.''  For each cancer pixel the distance to its nearest immune pixel is computed, yielding a multiset of distance values.  Symmetrically, for each immune pixel the distance to its nearest cancer pixel is computed, yielding a second multiset of distance values.  The means and variances of these two multisets serve as features for classification.  However, this featurization method is unstable, requires manual thresholding, and is insensitive to geometric information at larger spatial scales.

To address these shortcomings, the authors construct a 2-parameter sublevel filtration from the two images associated to each tissue sample.  They then compute the $H_1$ persistent homology module of this filtration, and vectorize this using each of three different MPH feature maps: the \emph{multiparameter persistence image} introduced in the same paper; the multiparameter persistence kernel of \cite{corbet2019kernel} (using sliced Wasserstein kernels \cite{carriere2017sliced}); and the multiparameter persistence landscape \cite{vipond2018multiparameter}.  Several 1-parameter persistence feature maps are also considered; these are obtained by restricting the bipersistence module to a single diagonal line. Comparing the performance of these different feature maps in the classification task, the authors find that the three multiparameter approaches perform similarly to each other, and far better than either the 1-parameter persistence approaches or the usual approach based on nearest neighbor distributions.
\section{Multifiltrations}
\label{Sec:Multifiltrations}
The \emph{degree-Rips} bifiltration of Definition \ref{Def:Degree_Rips} is just one of many natural constructions of bifiltrations from data.  Here, we introduce several interesting constructions of bifiltrations and trifiltrations that have either already played a significant role in MPH or that we think might in the future.  

\subsection{Superlevel-Rips Bifiltrations}
\begin{definition}[\cite{carlsson2009theory}] 
 \label{def:superlevel-rips}
 Consider a finite metric space $X$ and function $\gamma:X\to \R$.
The \emph{superlevel-Rips bifiltration}  $\SupVR(\gamma)\colon\R^{\op}\times [0,\infty)\to \mathbf{Simp}$ is given by \[\SupVR(\gamma)_{a,r}=\Rips(\gamma^{-1}[a,\infty))_r.\]
 \end{definition}
The sublevel-Rips bifiltration can be defined analogously, and \v Cech complexes can also be used in place of Rips complexes.  In \cite{carlsson2009theory}, several choices of the function $\gamma$ are proposed; we mention two of these.
 
\begin{example}[Density-Rips bifiltrations]\label{Ex:Density_Rips}
Informally, a \emph{density function} on a metric space $X$ is a function $\gamma\colon X\to [0,\infty)$ whose value is high in dense regions of the data and low near sparse regions of the data.  Standard examples include kernel density functions, e.g., with a Gaussian or disk kernel, and nearest neighbor density functions.

Given a density function $\gamma\colon X\to [0,
\infty)$, we call the superlevel-Rips bifiltration $\SupVR(\gamma)$ a \emph{density-Rips} bifiltration.  This is a practical, computationally tractable choice of density-sensitive bifiltration on a finite metric space.  We have already discussed an application to cancer imaging in \cref{Sec:Applications}.
\end{example}

\begin{example}[Eccentricity-Rips bifiltrations]\label{Ex:Eccentricity_Rips}
Define $\gamma\colon X\to [0,
\infty)$ by \[\gamma(x)=\frac{1}{|X|}\sum_{y\in X} d(x,y),\] i.e., $\gamma(x)$ is the average distance of $x$ to points in $X$.  We call $\gamma$ an \emph{eccentricity} function, and $\SupVR(\gamma)$ the \emph{eccentricity-Rips} bifiltration of $X$.  When $X$ has the structure of multiple tendrils emanating radially from a central core, the superlevel sets of $\gamma$ decompose into clusters in corresponding to the tendrils of $X$, which can be detected as persistent features in $H_0(\SupVR(\gamma))$. 
\end{example}

\begin{example}\label{Ex:Application}
In some applications, the data set $X$ comes equipped with a function $\gamma\colon X\to \R$, and we can consider the superlevel-Rips bifiltration of this.  For one example, in the study of time-varying data, $\gamma(x)$ may be the time of appearance of a data point.  As a second example, in applications to TDA computational chemistry, $X$ is taken to be the atom centers of a ligand (drug candidate), with the Euclidean distance and $\gamma$ is taken to be the partial charge function \cite{keller2018persistent}.  
\end{example}

\subsection{Interlevel-Rips Trifiltrations}\label{Sec:Interlevel_Rips}
In defining the superlevel-Rips bifiltration, we have made a choice of filtration direction: we have chosen to filter $X$ by the superlevel sets of $\gamma$ rather than by sublevel sets.  This is natural for some examples of  $\gamma$, e.g., those of Example \ref{Ex:Density_Rips} and Example \ref{Ex:Eccentricity_Rips}.  But for other examples of $\gamma$, such as those of Example \ref{Ex:Application}, imposing such a choice of direction can be unnatural.  For such examples, it is arguably more natural to consider a direction-agnostic variant of the construction, which we will call the \emph{interlevel-Rips trifiltration} \cite{gafvert2021topological}.  The idea is to filter $X$ not by superlevel sets, but by \emph{interlevel-sets}, i.e., the inverse images of intervals. In the following $\overline{U}$ denotes the subposet of $\R^{\op}\times \R$ given by $\overline{U}=\{(a,b) : a\leq b\}$. 

\begin{definition}\label{Def:Function_Rips_Trifiltration}
Consider a finite metric space $X$ and function $\gamma\colon X\to \R$.  The \emph{interlevel-Rips trifiltration} 
 $\overline{U}\times [0,\infty)\to \mathbf{Simp}$ is defined by  \[\IntVR(\gamma)_{(a,b),r}=\Rips(\gamma^{-1}[a,b])_r.\]
 \end{definition}
 
Observe that for $\gamma$ bounded, $\IntVR(\gamma)$ is a refinement of $\SupVR(\gamma)$.

In \cref{Sec:Interlevel_Sets}, we will consider a similar construction, the \emph{interlevel bifiltration}, which is defined for any real-valued function on a topological space.  In the case where one has a 1-parameter family of metrics $d=\{d_t\}_{t\in [0,1]}$ on a fixed finite set $S$, M\' emoli and Kim give a different construction of an interlevel-Rips trifiltration \cite{kim2021spatiotemporal}, which is perhaps more natural in that setting.  Such 1-parameter families arise in the study of collective motion of animals, e.g., flocking and swarming \cite{doi:10.1063/1.5125493,olfati2006flocking,bhaskar2019analyzing,toner1998flocks}, and also in the topological study of time series via sliding window embeddings \cite{perea2015sliding}.

\subsection{Multifiltrations from Images}
Many applications of 1-parameter persistent homology concern image analysis, where sublevel filtrations are often used.  There is not yet a consensus on what the most natural or useful multifiltrations are for image analysis, but one promising idea is that a second persistence parameter can be used to thicken sublevel or superlevel sets, thereby introducing some sensitivity to the \emph{width} of features that the ordinary sublevel and superlevel filtrations lack.

One construction along these lines is due to Chung, Day, and Hu \cite{chung2021multi}:  For $(X,d)$ a metric space and a function $f\colon X\to \R$, define a bifiltration $F\colon \R^2\to \Top$ by 
\[
F_{a,t}=
\begin{cases}
 \{y\in X :  d\left(y,f^{-1}[a,\infty)    \right)  \leq t\} & t\geq 0,\\
 \{y\in f^{-1}[a,\infty) : d\left(y,X\setminus f^{-1}[a,\infty) \right) \geq -t\} & t < 0.
 \end{cases}
\]
By considering interlevel sets rather than superlevel sets, one can also define a more symmetric 3-parameter analogue of this filtration.

Intuitively, these multifiltrations should, in some sense, exhibit better robustness to ``spike noise'' (i.e., noise of large magnitude but small width) than the ordinary sublevel or interlevel filtrations, but to the best of our knowledge there is no theoretical result along these lines.

\subsection{Parameter-Free, Density-Sensitive Multifiltrations}\label{Sec:Density-Sensitive Multifiltrations}
One important drawback of the density-Rips bifiltration introduced above is that, as a rule, the density function $\gamma$ depends on choice of \emph{bandwidth parameter}, e.g., the variance of the Gaussian kernel, or the number of nearest neighbors.  As suggested in \cref{Sec:Not_Robust}, it is generally preferable to work with constructions which do not require us to choose parameters. 

The degree-Rips bifiltration of Definition \ref{Def:Degree_Rips} is one natural choice of density-sensitive bifiltration whose construction requires no parameters.  Its simplicity is appealing, and as we discuss below, it has the advantage of being computable.  However, there are other natural alternatives, some of which have superior robustness properties (\cref{Thm:Prohorov_Stability}).  We now discuss several of these.

\paragraph{The Multicover Bifiltration}
Given $X\subset \R^n$, let $\mu_X$ denote the counting measure of $X$. 
\begin{definition}
Given any metric space $X$ and measure $\mu$ on $X$, the \emph{measure bifiltration} of $\mu$ is the $(0,\infty)^{\op}\times [0,\infty)$-indexed bifiltration $\mathcal M(\mu)$ defined as 
\[
\mathcal M(\mu)_{m,r}=
\{x\in X : \mu(B(x,r))\geq m\},
\]
where as in Definition \ref{Def:Distances_on_Point_Clouds}, $B(y,r)$ denotes the closed ball in $X$ of radius $r$ centered at $x$.  For $X\subset \R^n$, the bifiltration $\mathcal M (\mu_X)$ is called the \emph{multicover bifiltration} \cite{chazal2011geometric,sheehy2012multicover}.
\end{definition}
 Fixing $m=1$, we recover the offset filtration of Example \ref{Ex:Offset_Filtration}; the multicover bifiltration can thus be seen as a density-sensitive refinement of this.

\paragraph{Subdivision Bifiltrations}
Let $\bary(T)$ denote the barycentric subdivision of a simplicial complex $T$, i.e., the abstract simplicial complex whose simplices are sets $\{\sigma_1,\ldots,\sigma_l\}$ of simplices in $T$ such that $\sigma_i\subset \sigma_{i+1}$ for $\in 1,\ldots, l-1$.  Note that vertices of $\bary(T)$ correspond to simplices of $T$.  For $m\in (0,\infty)$, let $\mathcal S(T)_m$ denote the maximal subcomplex of $\bary(T)$ whose vertices correspond to simplices in $T$ of dimension at least $m-1 $.  These subcomplexes assemble into a $(0,\infty)^{\op}$-filtration $\mathcal S(T)$ of $\bary(T)$.  Moreover, for any simplicial filtration $F$, applying this construction at each index yields a bifiltration.
\begin{definition}[\cite{sheehy2012multicover}]
  When the above construction is applied to the filtration $F= \Rips(X)$, we call the resulting $(0,\infty)^{\op}\times [0,\infty)$-indexed bifiltration  the \emph{subdivision-Rips bifiltration} of $X$, and denote it as $\SRips(X)$. A subdivision-\v Cech bifiltration can be defined analogously.
\end{definition}
Vertices of $\SRips(X)_{m,r}$ correspond to cliques in $\Rips(X)_{r}$ of size at least $m$, which for large $m$  represent dense regions in $X$ at scale $r$.  In this sense, this bifiltration is indeed density sensitive.  

It follows from work of Sheehy and Cavanna \cite{cavanna2017and,sheehy2012multicover} that slight variants of the subdivision-\v Cech and the multicover bifiltrations are topologically equivalent  (more precisely, weakly equivalent \cite{blumberg2017universality}) and therefore have isomorphic persistent homology modules; see also \cite[Section 4]{blumberg2020stability}.  However, both the Rips and \v Cech subdivision filtrations have exponentially many vertices with respect to $|X|$, and are therefore too large to handle directly in practical computations.  

\paragraph{The Rhomboid Bifiltration}
For $X\subset \R^n$ finite, Corbet et al.  \cite{corbet2021computing} showed that a polyhedral bifiltration in $\R^{n+1}$ called the \emph{rhomboid bifiltration}, introduced by Edelsbrunner and Osang \cite{edelsbrunner2021multi}, is topologically equivalent to the multicover bifiltration of $X$.   The rhomboid bifiltration is a density-sensitive extension of the Delaunay filtration.  We refer the reader to \cite{corbet2021computing,edelsbrunner2021multi} for the definition.  The rhomboid bifiltration has size $O(|X|^{n+1})$.   

\paragraph{The Kernel Bifiltration}
Given a \emph{kernel function} $K\colon [0,\infty)\to [0,\infty)$ satisifying mild conditions, Rolle and Scoccola extend the definition of the associated density-Rips bifiltration to obtain a 3-parameter filtration, essentially by considering all possible values of the kernel bandwidth \cite{rolle2020stable}.  The construction thus no longer depends on a persistence parameter.  For a suitable choice of $K$, this trifiltration also extends the degree-Rips bifiltration.

\subsection{Computation of Density-Sensitive Bifiltrations} 
For a simplicial multifiltration $F$ indexed by $P$ and $a\in P$, let $\mathrm{births}(a)$ denote the number of simplices in $F_a$ not  present in $F_b$ for any $b<a$.  Define the \emph{size} of $F$ to be $\sum_{a\in P} \mathrm{births}(a)$.

For fixed $K$, the $K$-dimensional skeleton of the degree-Rips bifiltration of a metric space $X$ has size $O(|X|^{K+2})$.  However, if the bifiltration is coarsened to lie on a grid of constant size, then the $K$-skeleton has size $O(|X|^{K+1})$, which agrees (asymptotically) with that of the ordinary Rips filtration.  Using a line sweep algorithm designed by Roy Zhao and implemented in RIVET, these low-dimensional skeleta are readily computed in practice \cite{lesnick2021computing}.

Edelsbrunner and Osang \cite{edelsbrunner2020simple} recently introduced a clever and relatively simple algorithm for computing the rhomboid bifiltration.  The algorithm and its complexity analysis depend on a choice of algorithm for computing weighted Delaunay bifiltrations.  For a suitable such choice, the algorithm computes the rhomboid bifiltration of $X\subset \R^3$ in time $O(|X|^5)$.  An implementation of the algorithm is available  \cite{osang2020rhomboid}, and can compute the full rhomboid bifiltration of at least 200 points in $\R^3$; examples and timing results appeared in \cite{corbet2021computing} and \cite{edelsbrunner2020simple}.
This brings us close to practical computation of multicover PH for low-dimensional data sets, though it seems that further algorithmic work is needed to fully realize the practical potential of these ideas.  

The problem of computing the homology of the subdivision-Rips bifiltration is not yet understood.  As indicated above, the subdivision-Rips bifiltration is too large to be directly computed, but it may be that a smaller bifiltration has (exactly or approximately) the same homotopy type.

\section{Metrics and Stability}
\label{sec:metrics}
There have been many proposals for metrics on multiparameter persistence modules.  In this section we discuss a few of these and use them to formulate stability theorems. We focus primarily on the interleaving and matching distances, but also briefly discuss several other distances.

\subsection{The Interleaving Distance}\label{Sec:Interleavings}
The interleaving distance \cite{chazal2009proximity,lesnick2015theory,chazal2012structure}, an algebraic generalization of the bottleneck distance, is the best known and most thoroughly studied distance on multiparameter persistence modules.  As we will see, it is a natural choice of distance for formulating stability and inference results for multiparameter persistent homology.

For intuition, let us first consider a simple case:  A \emph{1-interleaving} between $\Z$-persistence modules $M$, $N$ is a commutative diagram of vector spaces of the following form, extending $M$ and $N$:
\[
\begin{tikzcd}
\cdots \ar[r]\arrow[dashed]{dr} & M_{-2}\ar[r]\arrow[dashed]{dr}  &M_{-1}\ar[r]\arrow[dashed]{dr} & M_0\ar[r]\arrow[dashed]{dr} & M_{1}\ar[r]\arrow[dashed]{dr} &M_{2}\ar[r]\arrow[dashed]{dr}& \cdots \\
\cdots \ar[r]\arrow[dashed]{ur} & N_{-2}\ar[r]\arrow[dashed]{ur}  &N_{-1}\arrow[dashed]{ur}\ar[r] & N_0\arrow[dashed]{ur}\ar[r] & N_{1}\arrow[dashed]{ur}\ar[r] &N_{2}\arrow[dashed]{ur}\ar[r] & \cdots
\end{tikzcd}
\]
Similarly, a \emph{2-interleaving} between $M$ and $N$ is a commutative diagram 
\[
\begin{tikzcd}
\cdots \ar[r]\arrow[dashed]{drr} & M_{-2}\ar[r]\arrow[dashed]{drr}  &M_{-1}\ar[r]\arrow[dashed]{drr} & M_0\ar[r]\arrow[dashed]{drr} & M_{1}\ar[r]\arrow[dashed]{drr} &M_{2}\ar[r]& \cdots \\
\cdots \ar[r]\arrow[dashed]{urr} & N_{-2}\ar[r]\arrow[dashed]{urr}  &N_{-1}\arrow[dashed]{urr}\ar[r] & N_0\arrow[dashed]{urr}\ar[r] & N_{1}\arrow[dashed]{urr}\ar[r] &N_{2}\ar[r] & \cdots
\end{tikzcd}
\]
For any integer $p\geq 0$, we can define a $p$-interleaving between $M$ and $N$ analogously, taking the arrows between $M$ and $N$ to increase indices by $p$.    
Moreover, the definition extends to $\Z^n$, by taking the arrows to increase indices by the vector $\vec p:=(p, \ldots, p)$.  

We can also define interleavings between $\R^n$-indexed modules in essentially the same way.  We now give the formal definition.
For $u\in [0, \infty)^n$, define the \emph{$u$-shift} $(-)(u)$ to be the endofunctor on $\R^n$-modules given on objects by 
 $M(u)_a= M_{u+a}$ and $M(u)_{a,b} = M_{a+u, b+u}$, and on morphisms by $f(u)_a = f_{u+a}$. Furthermore, let $\phi_M^u\colon M\to M(u)$ be the morphism whose restriction to each $M_a$ is the linear map $M_{a,a+u}\colon M_a \to M_{a+u}$. 

\begin{definition}\label{Def:Interleavings}
Given $\epsilon\in [0, \infty)$, we say that persistence modules $M,N\colon \R^n \to \kVect$ are \emph{$\epsilon$-interleaved} if there exist morphisms 
\[f\colon M\to N(\vec\epsilon) \qquad  \quad g\colon N\to M(\vec\epsilon)\] such that 
\[g(\vec\epsilon)\circ f = \phi_M^{2\vec\epsilon}\qquad \qquad f(\vec\epsilon)\circ g = \phi_N^{2\vec\epsilon}.\]
The \emph{interleaving distance} is defined by
\[d_I(M,N) = \inf\,\{\epsilon\geq 0 \colon M \text{ and } N \text { are $\epsilon$-interleaved.}\}.\]
\end{definition}

\begin{remark}
Clearly, Definition \ref{Def:Interleavings} generalizes to functors $\R^n\to \mathbf{C}$, for any category $\mathbf{C}$.  The cases $\mathbf{C}=\Top$ and $\mathbf{C}=\mathbf{Set}$ play important roles in TDA.
\end{remark}

The following theorem tells us that the interleaving distance is an extension of the bottleneck distance to multiparameter persistence modules. 

\begin{theorem}[The isometry theorem]
For any p.f.d. persistence modules $M,N\colon \R\to \Vect$, 
\[d_I(M,N) = d_B(\B{M}, \B{N}).\]
\label{thm:IST}
\end{theorem}
The inequality $d_I(M,N) \geq d_B(\B{M}, \B{N})$, known as the \emph{algebraic stability theorem}, was first proven by Chazal et al. \cite{chazal2009proximity,chazal2012structure}.  Other proofs have appeared in \cite{bauer2015induced,bauer2016persistence,bjerkevik2016stability}.  The converse inequality first appeared in \cite{lesnick2015theory}.

\paragraph{Multiparameter Algebraic Stability}

The bottleneck distance $d_B$ on intervals in $\R$ admits a natural extension to collections of intervals in $\R^n$.  It is convenient to define this via the interleaving distance:
\begin{definition}[Bottleneck Distance on Barcodes in $\R^n$ \cite{botnan2018algebraic}]\label{Def:bottleneck_general}
For a matching $\chi$ between barcodes $\mathcal{C}$ and $\mathcal{D}$ in $\R^n$, let
\begin{equation}
\cost(\chi):=\max\left( \max_{(X,Y)\in \chi} d_I(k_X,k_Y),\ \max_{ X \in \mathcal{C}\cup \mathcal{D} \text{ unmatched}}\ d_I(k_X,0)\right).
\label{eq:cost-matching}
\end{equation}
We then define the $d_B(\mathcal{C},\mathcal{D})$ to be the infimal cost of a matching between $\mathcal{C}$ and $\mathcal{D}$, exactly as in Definition \ref{Def:bottleneck}.  
\end{definition}

For certain types of intervals arising naturally in TDA, e.g., those considered in \cref{Thm:Generalized:Algebraic_Stability} below, one has simple, explicit formulae for the interleaving distance between interval modules.  On barcodes consisting of such intervals, these formulae yield simpler formulations of Definition \ref{Def:bottleneck_general}  \cite{botnan2018algebraic, bjerkevik2016stability,botnan2022bottleneck}.

Given that $d_I$ and $d_B$ are both defined the $n$-parameter setting, one might hope that the isometry theorem extends to interval-decomposable $\R^n$-persistence modules.  The following summarizes what is known about this:

\begin{theorem}[\cite{botnan2018algebraic,bjerkevik2016stability,botnan2022bottleneck}]\label{Thm:Generalized:Algebraic_Stability}
Fix $n\geq 2$, and assume that $M$ and $N$ are finitely presented, interval-decomposable $\R^n$-modules.
\begin{itemize}
\item[(i)] There is no constant $c$ such that \[d_B(\B{M}, \B{N}) \leq c\, d_I(M,N)\] for all such $M$ and $N$.
\item[(ii)] If all intervals in $\B{M}$ and $\B{N}$ are \emph{blocks} (Definition \ref{Def:Block}), then \[d_I(M,N)=d_B(\B{M}, \B{N}).\]
\item[(iii)] If $M$ and $N$ are free (Definition \ref{Def:Free}), then \[d_I(M,N)\leq d_B(\B{M}, \B{N})\leq (n-1)d_I(M,N).\] These bounds are tight for $n=2$ and $n=4$.  
\item[(iv)] If all intervals in $\B{M}$ and $\B{N}$ are \emph{rectangles} (Definition \ref{Def:Rectangle}) or all are \emph{hooks} (Definition \ref{Def:Hook}), then \[d_I(M,N)\leq  d_B(\B{M}, \B{N}) \leq (2n-1)d_I(M,N).\] 
These bounds are tight for $n=2$.
\end{itemize}
\end{theorem}

\paragraph{Universality}
\begin{definition}
We will say that a distance $d$ (i.e., extended pseudometric) on multiparameter persistence modules is \emph{stable} if for all topological spaces $W$, functions $\gamma,\kappa:W\to \R^n$, and $i\geq 0$, we have 
\[d(H_i(\mathcal S^\uparrow(\gamma)),H_i(\mathcal S^\uparrow(\kappa))) \leq \sup_{w\in W} \|\gamma(w)-\kappa(w)\|_\infty.\]

\end{definition}

\begin{theorem}[Universality of the Interleaving Distance \cite{lesnick2015theory}]\label{Thm:Universality}
\mbox{}
\begin{enumerate}[(i)]
\item The interleaving distance $d_I$ is stable.
\item Assume that the field $k$ is prime, i.e., $k=\mathbb Q$ or $k=\Z/p\Z$ for some prime $p$.  Then for any other stable distance $d$ on multiparameter persistence modules, we have $d\leq d_I$.
\end{enumerate}
\end{theorem}

\begin{remarks}
~
\begin{enumerate}
\item A version of \cref{Thm:Universality} for the special case of 1-parameter persistence and $0^\mathrm{th}$ homology first appeared in work by d'Amico et al. \cite{d2010natural}. 
\item The proof of \cref{Thm:Universality}\,(i) turns out to be trivial, but the proof of \cref{Thm:Universality}\,(ii) is not.
\item In view of \cref{thm:IST}, \cref{Thm:Universality}\,(i) extends \cref{Thm:Stability}\,(i).  
\item The generalization of \cref{Thm:Universality} to arbitrary fields is an open question.  
\item A similar universality result is given for filtrations in \cite{blumberg2017universality}, using the \emph{homotopy interleaving distance}.  
\end{enumerate}
\end{remarks}

\subsubsection{Computation} While the interleaving distance has good theoretical properties, it turns out to be hard to compute, even for very simple bipersistence modules. 
\begin{theorem}[\cite{bjerkevik2019computing}]\label{Thm:Computing_Interleaving_Distance_is_Hard}
For $k$ a finite field, it is NP-hard to approximate the interleaving distance on $2$-parameter persistence modules within a factor of 3. 
\label{thm:np-hard-int}
\end{theorem}

\cref{thm:np-hard-int} is proved by considering interleavings of \emph{staircase}-decomposable modules, where a staircase is a particularly simple type of interval.  That said, the hardness result still holds if all modules are assumed to be indecomposable.

The proof of \cref{Thm:Computing_Interleaving_Distance_is_Hard} is a reduction from the following problem \cite{bjerkevik2017computational}, which is shown to be NP-complete in \cite{bjerkevik2019computing}: let $A$ and $B$ be $n\times n$ matrices of distinct variables and suppose that we assign the value 0 to each variable in a given subset of the variables.  Is it possible to assign a value in $k$ to each of remaining variables such that the resulting product $AB$ is the identity matrix $I_n$?
\begin{example}
Let $A$ and $B$ be $3\times 3$ matrices with entries set to 0 as follows 
\[A= \begin{bmatrix} * & * & * \\ * & 0 & * \\ * & * & 0\end{bmatrix} \qquad \qquad  B =\begin{bmatrix} * & * & * \\ * & * & 0 \\ * & 0 & *\end{bmatrix} \]
In this case the decision problem has a positive answer:
\[\begin{bmatrix} 1 & 1 & 1 \\ 1 & 0 & 1 \\ 1 & 1 & 0\end{bmatrix}\cdot\begin{bmatrix} -1 & 1 & 1 \\ 1 & -1 & 0 \\ 1 & 0 & -1\end{bmatrix} = I_3\]
\end{example}

\subsection{The Matching Distance}\label{Sec:Matching_Dist}
\cref{thm:np-hard-int} motivates the search for a more computable surrogate for the interleaving distance.  The matching distance \cite{cerri2013betti} has emerged as a popular choice.  
It is defined by restricting the modules to (suitably parameterized) affine lines with positive slope, and taking bottleneck distances:

\begin{definition}\label{Matching Distance}
The \textbf{matching distance} between p.f.d. $\R^n$-indexed modules $M$ and $N$ is given by: \[ d_{\text{match}}(M,N) = \sup_L\ d_B(\B{M\circ L}, \B{N\circ L}),\]
where $L:\R\to \R^n$ ranges over parameterized lines of the form $L(t) = vt + b$ where $v\in [1,\infty)^n$ and $b\in \R^n$.
\end{definition}

Note that $d_M(M,N)$ depends only on the fibered barcodes of $M$ and $N$.  The following is an easy consequence of the algebraic stability theorem:
\begin{proposition}[\cite{landi2018rank,cerri2013betti}]\label{Prop:Matching_Bounds_Interleaving}
For all p.f.d. $\R^n$-persistence modules $M$ and $N$, \[d_{\textup{match}}(M,N)\leq d_I(M,N).\]
\end{proposition}

\subsubsection{Computation}
In 2011, Biasotti et al. \cite{biasotti2011new}  gave a quad tree-based algorithm to approximate the matching distance between bipersistence modules.  More recently, Kerber and Nigmetov \cite{kerber2020efficient} revisited their approach, arriving at a more efficient version, and released an implementation as part of the Hera software package.

The matching distance on bipersistence modules can in fact be exactly computed in polynomial time \cite{kerber2019exact,bjerkevik2021asymptotic}; the following describes the state of the art.

\begin{theorem}[\cite{bjerkevik2021asymptotic}]\label{Thm:ExactMatchingDist}
Given minimal free presentations of $\R^2$-indexed persistence modules $M$ and $N$ (see \cref{SubSec:Pres_Res}), $d_{\rm match}(M,N)$ can be computed deterministically in time $O(n^{6}\log n)$ and via randomization in time $O(n^{5}\log^3 n)$, where $n$ denotes the total number of generators and relations in the two presentations.
\end{theorem}

There is no literature on computing the matching distance between higher-parameter persistence modules, but it seems likely that the known approaches in the two-parameter case can be extended.

\subsection{Stability of Density Sensitive Bifiltrations}\label{Sec:Robustness}

The stability and robustness of the multicover, subdivision, and degree-Rips bifiltrations were recently studied in \cite{blumberg2020stability}, using the \emph{Prohorov} and \emph{Gromov-Prohorov} distances on measures and interleavings on persistence modules.  The Prohorov distance $d_{Pr}$ is a classical distance between measures on the same metric space, which can be thought of as a measure-theoretic analogue of the Hausdorff distance; the \emph{Gromov-Prohorov distance} $d_{GPr}$ is an extension of $d_{Pr}$ to measures on different metric spaces~\cite{greven2009convergence,janson2020gromov}:

\begin{definition}
The \emph{Prohorov distance} between measures $\mu$ and $\eta$ on a metric space $(Z,d_Z)$ is given by 
\[
d_{Pr}(\mu,\eta) = \sup_A \inf\{\delta \geq 0 :  \mu(A) \leq \eta(A^\delta) + \delta \textup{ and } \eta(A) \leq \mu(A^\delta) + \delta\},
\] 
where $A \subset  Z$ ranges over all closed sets and 
\[A^\delta=\{y\in  Z : d_{Z}(y,a) < \delta\textup{ for some }a \in A.\}\]
\end{definition}

\begin{definition}[\cite{greven2009convergence}]
The \emph{Gromov-Prohorov} distance between measures $\mu_X$ and $\mu_Y$ on metric spaces $X$ and $Y$ is given by 
\[
d_{GPr}(\mu_X,\mu_Y) = \inf_{\varphi, \psi} d_{Pr}(\varphi_*(\mu_X), \psi_*(\mu_Y)),
\]
where $\varphi \colon X \to Z$ and $\psi \colon Y \to Z$ range over all
isometric embeddings into a common metric space $Z$. 
\end{definition}
A important property of these distances is that they are stable to the addition of outliers, in a strong sense \cite[Remark 2.17]{blumberg2020stability}.  

Recall that for $X\subset \R^n$, we let $\mu_X$ denote the counting measure of $X$, seen as a measure on $\R^n$.  Similarly, for $X$ a finite metric space, we let $\eta_X$ denote the counting measure of $X$, seen as a measure on $X$.
\setlength{\maxdepth}{1pt}
\begin{theorem}[Stability of the Multicover and Subdivision Bifiltrations \cite{blumberg2020stability}]\mbox{}\label{Thm:Prohorov_Stability}
\begin{itemize}
\item[(i)] For all $X,Y\subset \R^n$ finite and $i\in \mathbb N$, \[d_{I}(H_i(\mathcal M(\mu_X)),H_i(\mathcal M(\mu_Y)))\leq d_{Pr}(\mu_X,\mu_Y).\]
\item[(ii)] For all finite metric spaces $X,Y$ and $i\in \mathbb N$, \[d_{I}(H_i(\SRips(X)),H_i(\SRips(Y)))\leq d_{GPr}(\eta_X,\eta_Y).\]
\end{itemize}
\end{theorem}
\cref{Thm:Prohorov_Stability} is a measure-theoretic analogue of items (ii) and (iii) of \cref{Thm:Stability}, and is proven using a closely analogous argument.  A variant of this result using normalized counting measures also holds, and implies a version which uses the Wasserstein distance on the measures in place of the Prohorov distance.  

In the same paper, an analogue of these results for degree-Rips bifiltrations is also given \cite[Theorem 1.6]{blumberg2020stability}, using the Gromov-Prohorov distance and a generalized definition of interleavings.  This turns out to be a qualitatively weaker kind of stability than the 1-Lipschitz stability of \cref{Thm:Prohorov_Stability}, but is tight (in a reasonable sense).  The upshot is that, as measured by interleavings, the degree-Rips bifiltration is qualitatively less robust to outliers than the subdivision-Rips bifiltration, but qualitatively more robust to outliers than the ordinary Rips filtration.

In a similar vein, Rolle and Scocolla have given a Lipschitz stability result for the kernel trifiltration (see \cref{Sec:Density-Sensitive Multifiltrations}) using the \emph{Gromov-Hausdorff-Prohorov distance,} an upper bound for the  Gromov-Prohorov which is density-sensitive but unstable to the addition of outliers \cite{rolle2020stable}.  This result implies a corresponding Lipschitz stability result for the degree-Rips bifiltration.  These ideas are applied in \cite{rolle2020stable} to obtain a stable and consistent clustering scheme based on 3-parameter persistence.

\begin{remark}[The Tension Between Robustness and Computability]
Together, the stability results discussed above and the discussion of computation in \cref{Sec:Density-Sensitive Multifiltrations} reveal a tension between robustness and computability for the 2-parameter PH of metric data: The multicover and subdivision-Rips bifiltrations have excellent robustness properties, but we do not yet have a practical computational framework for handling  either of them.  (However, for the multicover bifiltration of low dimensional data, we are already close.)  The degree-Rips bifiltrations, on the other hand, are computable, but the present theory guarantees robustness only in a far weaker sense.

This raises the following critical question: Can we develop a framework for MPH of metric data that is both computationally efficient and provably robust? 
\end{remark}

\subsection{$\ell^p$-Metrics on Multiparameter Persistence Modules}
The interleaving and matching distances are both $\ell^\infty$-distances, i.e., they can be defined in terms of $\ell^\infty$-metrics on Euclidean spaces \cite{bjerkevik2021ell}.  As such, they are insensistive to certain small scale differences between modules, which can be undesirable in both theory and applications.  To address this, it is natural to consider $\ell^p$-metrics on multiparameter persistence modules.  Several recent works have explored the question of defining such metrics \cite{skraba2020wasserstein,bjerkevik2021ell,bubenik2018wasserstein,giunti2021amplitudes,thomas2019invariants}.

 In the 1-parameter case, there is a standard $\ell^p$-generalization of the bottleneck distance, called the \emph{$p$-Wasserstein distance} and denoted $W_p^q$; here $q\in [0,1]$ is a parameter specifying an $\ell^q$-metric on $\R^2$ used in the definition of the distance.   $W^\infty_\infty$ is equal to the bottleneck distance, and varying $q$ changes the distance $W_p^q$ by at most a factor of 2.  The distances $W_p^q$ are used in many applications and in some TDA theory; see \cite[Section 1]{bjerkevik2021ell}.  
 
The problem of extending the Wasserstein distance to multiparameter persistence modules was first considered by Bubenik, Scott, and Stanley, who generalized $W_p^1$ to persistence modules indexed by an arbitrary abelian category \cite{bubenik2018wasserstein}.  However, in the case $p=\infty$, their distance is equal to neither the interleaving distance nor the matching distance, and has different qualitative properties.  

More recently, Bjerkevik and Lesnick \cite{bjerkevik2021ell} introduced $\ell^p$-extensions of the interleaving and matching distance to multiparameter persistence modules, called the $p$-\emph{presentation distance} $d_I^p$ and the \emph{$p$-matching distance} $d_{\mathrm{match}}^p$, respectively.  In the 1-parameter case, these distances are equal to $W_p^p$, and they share several of the desirable properties of the interleaving and matching distances discussed earlier in this section.  Most notably, on $n$-parameter persistence modules with $n\in \{1,2\}$, $d_I^p$ satisfies a universal property closely analogous to \cref{Thm:Universality}.  This result extends a fundamental  $\ell^p$-stability result for the 1-parameter persistent homology of cellular filtrations, due to Skraba and Turner \cite[Theorem 4.7]{skraba2020wasserstein}.  
 
\subsection{Other Metrics}
Besides the works referenced above, several other papers have introduced metrics for multiparameter persistence and developed theory for these; we briefly mention three: Scolamiero et. al \cite{scolamiero2017multidimensional} introduce the formalism of \emph{noise systems} and use this to define a family of metrics on multiparameter persistence modules.  Cerri et al. \cite{cerri2019geometrical,cerri2016coherent} study a variant of the matching distance requiring the matchings along different affine lines to be chosen coherently, and they prove a stability result for this.  McCleary and Patel \cite{mccleary2020edit} use M\" obius inversion to define a distance on the homology modules of certain lattice-indexed filtrations, and they establish functoriality and stability results for their approach.

\section{Minimal Presentations and Resolutions}\label{Sec:Min_Pres_Res}
In this section, we introduce minimal presentations and resolutions of persistence modules.  We then discuss the problem of computing minimal presentations, focusing on recent progress in the 2-parameter case. 

Minimal presentations  are particularly important in computation because they provide a memory-efficient representation of a persistence module, which then can be used as input to algorithms for computing invariants or distances.  Minimal presentations of persistent homology modules are often quite small in practice, in comparison to the filtrations from which they arise \cite{kerber2021fast}, which makes them convenient for computation.  However, minimal presentations are usually not unique, so they cannot be directly used in TDA the ways barcodes are used, e.g., as input to machine learning algorithms or statistical tests.  

\subsection{Free Modules}\label{Sec:Free_Modules}
Let $P=T_1\times \cdots \times T_n$, where each $T_i$ is a totally ordered set, and for $z\in P$, let $\langle z\rangle$ denote the interval $\{p\in P \colon p\geq z\}$. 

\begin{definition}\label{Def:Free}
A  $P$-persistence module $F$ is \emph{free} if it is interval-decomposable and all intervals in the barcode $\B{F}$ are of the form $\langle z\rangle$. 
\end{definition}
For $F$ free, we let $|F| = |\B{F}|$, where $|\cdot|$ denotes the cardinality of a multiset.  

Given a $P$-persistence module $M$ and $v\in M_z$, we write $\gr(v)=z$.  We say that $S\subset \bigcup_{z\in P} M_z$ is \emph{a set of generators} for $M$ if for any $v\in \bigcup_{z\in P} M_z$, \[v=\sum_{i=1}^m c_i M_{\gr(v_i),\gr(v)}(v_i)\] for some $v_1,v_2,\ldots, v_m\in S$ and scalars $c_1,\ldots, c_m\in k$.  We say $M$ is \emph{finitely generated} if there exists a finite set of generators for $M$.  

A \emph{basis} of a free module $F$ is a minimal generating set.  It can be shown via elementary linear algebra that if $B$ is any basis for a free module $F$, then for all $z\in P$, the number of elements in $B$ of grade $z$ is equal to the number of copies of $\langle z\rangle$ in $\B{F}$. 

\paragraph{Matrix Representation of Morphisms of Free Modules}
Let $\gamma\colon F'\to F$ be a morphism of finitely generated free modules, and let $B'=\{B'_i\}_{i=1}^{d'}$ and $B=\{B_i\}_{i=1}^{d}$ be ordered bases of $F'$ and $F$, respectively.  In analogy with ordinary linear algebra, we can represent $\gamma$ with respect to these bases via a matrix $[\gamma]$ with coefficients in the field $k$, together with a $P$-valued label for each row and each column of the matrix.  To explain the details, for $z\in P$, we represent $v\in F_z$ with respect to $B$ as a vector $[v]^{B}\in k^{|B|}$, by taking $[v]^B$ to be the unique vector such that $[v]_i^B=0$ if $\gr(B_i)\not\leq z$ and \[v=\sum_{i\colon \gr(B_i)\leq z} [v]^B_i F_{\gr(B_i),z}(B_i).\] Thus, $[v]^B$ records the field coefficients in the linear combination of $B$ giving $v$.   
We now define $[\gamma]$ as follows:
\begin{itemize}
\item the $j^{\mathrm{th}}$ column is $[\gamma(B'_j)]^{B}$,
\item the label of the $j^{\mathrm{th}}$ column is $\gr(B'_j)$,
\item the label of the $i^{\mathrm{th}}$ row is $\gr(B_i)$.
\end{itemize}

In the literature on multi-graded commutative algebra, the matrix $[\gamma]$ is called a \emph{monomial matrix} \cite{miller2004combinatorial}.

\subsection{Free Presentations and Resolutions}\label{SubSec:Pres_Res}
\begin{definition}~
\begin{itemize}
\item[(i)]
A \emph{(free) presentation}
$\mathcal F$ of a persistence module $M$ is a morphism of free modules
\begin{equation}
 F_1 \xrightarrow{\phi_1} F_0
\label{eq:minpres}
\end{equation}
with $\coker \phi_1 \cong M$.
\item[(ii)] A \emph{(free) resolution} $\mathcal{F}$ of $M$ is an exact sequence of free modules
\[\cdots \xrightarrow{\phi_3} F_2 \xrightarrow{\phi_2} F_1 \xrightarrow{\phi_1} F_0\]
with $\coker \phi_1 \cong M$.
\end{itemize}
\label{def:free-res}
\end{definition}
$M$ is said to be \emph{finitely presented} if there exists a presentation $\mathcal F$ of $M$ with $F_0$ and $F_1$ finitely generated.

\begin{remark}
In view of the discussion \cref{Sec:Free_Modules}, a free presentation $\gamma\colon F_0\to F_1$ with each $F_i$ finitely generated can be represented with respect to a choice of bases for $F_0$ and $F_1$ as a labeled matrix; we call this a \emph{presentation matrix}, or sometimes (by a slight abuse of terminology) simply a \emph{presentation}.  Similarly, a resolution can be represented as a sequence of matrices.
\end{remark}

Of particular interest are minimal presentations and resolutions. Several equivalent definitions can be given; the following approach is perhaps the most transparent.

\begin{definition}~
\begin{itemize}
\item[(i)] A \emph{trivial resolution} is a direct sum of free resolutions the form
\[\cdots 0\to 0 \to  F \xrightarrow{\id_F} F \to 0\to 0\to \cdots \to 0.\] 
where the two copies of $F$ may appear at any two consecutive indices.  
\item[(ii)] A \emph{trivial presentation} is a free presentation of the form
\[F\oplus F' \xrightarrow{\id_F\oplus\, 0} F.\]
\end{itemize} 
\end{definition}

\begin{definition}~
\begin{itemize}
\item[(i)] A resolution $\mathcal{F}$ of a persistence module $M$ is \emph{minimal} if any resolution $\mathcal{F'}$ of $M$ is isomorphic to the direct sum of $\mathcal{F}$ with a trivial resolution. 
\item[(ii)] A presentation $\mathcal{F}$ of $M$ is \emph{minimal} if any presentation $\mathcal{F'}$ of $M$ is isomorphic to the direct sum of $\mathcal{F}$ with a trivial presentation. 
\end{itemize}
\label{def:minimal-free}
\end{definition}

\begin{theorem}[Structure of Minimal Resolutions]\label{Thm:Structure_of_Resolutions}
Let $M$ be a finitely presented $n$-parameter persistence module. 
\begin{enumerate}
\item[(i)] A minimal free resolution $\mathcal{F}$ of $M$ exists and is unique up to natural isomorphism.
\item[(ii)] Each module $F_i$ in $\mathcal{F}$ is finitely generated and $F_i = 0$ for $i\geq n+1$. 
\end{enumerate}
\label{thm:minres-prop}
\end{theorem}

\begin{remarks}~
\begin{itemize}
\item[(i)] \cref{Thm:Structure_of_Resolutions}\,(i) is variant of a standard result in commutative algebra.    Using Azumaya's theorem (\cref{Thm:KS}\,(ii)), the proof given in \cite{peeva2011graded} for the case of finitely generated $\Z$-graded $\Z^n$-modules adapts readily to a proof of our version.
\item[(ii)]  \cref{Thm:Structure_of_Resolutions}\,(ii) is a variant of the well-known \emph{Hilbert's basis theorem}; the standard proofs adapt to our multigraded setting.
\item[(iii)] \cref{Thm:Structure_of_Resolutions} implies that a minimal presentation of $M$ also exists and is unique up to isomorphism.
\end{itemize}
\end{remarks}

The following easy corollary of \cref{Thm:Structure_of_Resolutions}, due to Scolamiero et al. \cite{scolamiero2017multidimensional}, turns out to be very useful in the study of 2-parameter persistent homology.
\begin{corollary}\label{Cor:Free_Kernels}
If $f\colon F_1\to F_0$ is a morphism of free, finitely presented bipersistence modules, then $\ker f$ is free.
\end{corollary}

\begin{example}
Consider the following indecomposable $\N^2$-module $M$
\[
\begin{tikzcd}
0\ar[r] & 0\ar[r] & 0\ar[r] & 0 \\
k\ar[u]\arrow{r}& 0\ar[u]\ar[r]  & 0\ar[u]\ar[r]  & 0\ar[u] \\
k^2\arrow{u}{[1,0]}\arrow{r}{[1,1]} & k\ar[u]\ar[r]  & 0\ar[r]\ar[u]  & 0\ar[u] \\
k^2\arrow{u}{=}\arrow{r}{=}  & k^2\arrow{u}{[1,1]}\arrow{r}{[0,1]} & k\ar[r]\ar[u]  & 0\ar[u] 
\end{tikzcd}
\]
A minimal free resolution of $M$ is given is given in matrix form by 
\[
A_1=
\begin{blockarray}{cccccc}
(0,2) & (0,3) & (1,1) & (2,0) & (3,0) \\
\begin{block}{(ccccc)c}
1 & 0 & 1 & 0 & 1 & (0,0) \\
0 & 1 & 1 & 1 & 0 & (0,0) \\
\end{block}
\end{blockarray}
\]
\[A_2=
\begin{blockarray}{cccc}
(1,3) & (2,2)  & (3,1) \\
\begin{block}{(ccc)c}
\textcolor{white}{-}1 & \textcolor{white}{-}1 &  \textcolor{white}{-}0 & (0,2) \\
\textcolor{white}{-}1 & \textcolor{white}{-}0 & \textcolor{white}{-}1 & (0,3) \\
-1 & -1 & -1 & (1,1) \\
\textcolor{white}{-}0 & \textcolor{white}{-}1 & \textcolor{white}{-}0 & (2,0) \\
\textcolor{white}{-}0 & \textcolor{white}{-}0 & \textcolor{white}{-}1 & (3,0) \\
\end{block}
\end{blockarray}
 \]
\end{example}

\subsection{Computing Minimal Presentations and Resolutions}\label{Sec:Computing_Pres_Res}
Considerations of efficiency aside, minimal presentations and resolutions can be computed using Gr\"obner basis techniques such as \emph{Schreyer's algorithm} and its variants \cite{la1998strategies,schreyer1980berechnung,erocal2016refined}, which in fact work in much greater generality.  The application of Schreyer's algorithm to MPH was first explored in \cite{carlsson2010computing}.  However, the problem instances one considers in TDA have very special structure: the modules are multi-graded (see \cref{sec:multigraded}), the number of persistence parameters is very small (usually at most three), and the problems are very large but very sparse.  Thus one would expect that on TDA problems, specialized algorithms and implementations would far outperform more classical approaches designed and optimized for other purposes.

\paragraph{Short Chain Complexes as Input to the Computations}
Let $P$ be a product of finite totally ordered sets with maximum element $P_{\max}$, and let $F\colon P\to \Simp$ be a simplicial filtration such that $F_{P_{\max}}$ is a finite simplicial complex.  We will consider the problem of computing a presentation of $H_i(F)$ for fixed $i$.

In analogy with ordinary simplicial homology, $F$ has an associated chain complex of free $P$-modules 
\[\cdots \xrightarrow{\partial_{i+1}} C_i(F)\xrightarrow{\partial_{i}} C_{i-1}(F)\xrightarrow{\partial_{i-1}}\cdots \xrightarrow{\partial_{1}} C_{0}(F)\to 0,\]
where for $z\in P$, $C_i(F)_z:=C_i(F_z;k)$ is the usual simplicial chain vector space with coefficients in $k$, and the internal maps in $C_i(F)$ are inclusions.  The $i^{\mathrm{th}}$ homology module of this chain complex is exactly $H_i(F)$.  Given the filtration $F$, constructing this chain complex (or portions thereof) is generally straightforward.  

It is not hard to verify that $C_i(F)$ is free if and only if $F$ is 1-critical (i.e., if each simplex has a unique birth index; see \cref{Sec:Filtrations}).  If $F$ is not 1-critical, then we may use a simple construction due to Scolamiero et al. \cite{chacholski2017combinatorial} to convert the short chain complex \[C_{i+1}(F) \xrightarrow{\partial_{i+1}} C_i(F)\xrightarrow{\partial_{i}} C_{i-1}(F)\]
into a short chain complex of free $P$-modules
\begin{equation}\label{eq:Short_Chain_Complex}
X\xrightarrow{f} Y \xrightarrow{g} Z
\end{equation}
such that $H_i(F)\cong \ker{g}/\operatorname{im}{f}$.  This construction is implemented in RIVET in the 2-parameter case, where it is  used to compute minimal presentations of degree-Rips bifiltrations.

Thus, whether or not $F$ is 1-critical, we may assume the input to our presentation computation is a short chain complex of free $P$-modules as in \cref{eq:Short_Chain_Complex} whose homology module is isomorphic to $H_i(F)$. 

\paragraph{Algorithms for Bipersistence Modules}
An efficient algorithm for computing a minimal presentation of a bipersistence module from a short chain complex was introduced in \cite{lesnick2019computing} and subsequently improved in \cite{kerber2021fast}.  Extensive computational experiments described  in \cite{kerber2021fast} demonstrate that this approach scales well with the size of the filtrations, and is efficient enough for practical use in TDA.  In fact, this algorithm underlies the MPH computations of \cite{vipond2021multiparameter} which we discussed  \cref{Sec:Applications}. 

The core computational engine behind this approach is a simple matrix reduction algorithm called the \emph{bigraded reduction}.  It is very similar to the standard matrix reduction algorithm for computing persistent homology \cite{zomorodian2005computing}, but instead of reducing an entire matrix in one pass, it proceeds by reducing submatrices of increasing size.  Slight variants of the bigraded reduction solve three fundamental algebra problems involving free bipersistence modules \cite{lesnick2019computing}: One variant computes a basis for the kernel of a map $\gamma\colon F\to G$ of free modules; a second variant minimizes a set of generators of a persistence module; and a third computes a minimal Gr\" obner basis for a submodule of a free module.  As explained in \cite{lesnick2019computing}, the first two variants are asymptotically more memory efficient than their more general classical counterparts, which first compute a Gr\"obner basis of $\im \gamma$.  

The first two variants are used in \cite{lesnick2019computing,kerber2021fast} to compute presentations.  Specifically, given a short chain complex $X\xrightarrow{f} Y \xrightarrow{g} Z$ of free modules, 
and letting $M=\ker{g}/\operatorname{im}{f}$, one computes a presentation matrix $Q$ for $M$ in three simple steps:
\begin{enumerate}
\item compute a minimal set of generators $S=\{s_1,\ldots, s_m\}$ for $\im f$,
\item compute an ordered basis $B$ for $\ker g$, and
\item using ordinary linear algebra, express each $s_i\in S$ as a linear combination of elements of $B$.  
\end{enumerate}
One then takes $Q$ to be the $|B|\times |S|$ matrix with $Q(*,i)=[s_i]^B$.  The row and column labels of $Q$ are taken to be the grades of the elements of $B$ and $S$, respectively.

This presentation will usually not be minimal, but it is \emph{semi-minimal}, in the sense that there is no presentation matrix with the same number of rows and fewer columns.  A standard procedure from commutative algebra for minimizing resolutions (essentially, a variant of Gaussian elimination) adapts to minimize such a presentation \cite[pages 127 and 166]{greuel2012singular}.  An efficient variant of this minimization algorithm in the setting of sparse matrices has been introduced in \cite{fugacci2018chunk}.  This algorithm in fact works for resolutions of multiparameter persistence modules with arbitrary numbers of parameters and, moreover, can be used more generally to minimize a chain complex of free modules.

Letting $l=|X|+|Y|+|Z|$, the above algorithm for computing a minimal presentation runs in time $\Theta(l^3)$ and memory $\Theta(l^2)$; perhaps surprisingly, this matches the asymptotic runtime and memory cost of the ordinary persistence algorithm in the 1-parameter setting \cite{zomorodian2005computing}.  

The approach adapts readily to compute a minimal resolution of a bipersistence module. In view of \cref{Thm:Structure_of_Resolutions}\,(ii), one needs only to perform a single additional kernel computation.

\paragraph{Further Directions}
The surprisingly good performance of these minimal presentation computations is cause for optimism about the prospects for applications of MPH.  Yet the speed and scalability of 2-parameter persistence computations is still well behind that of 1-parameter computations.  Recent approaches to 1-parameter computations use several clever optimizations that together make a huge difference in speed and scalability, as e.g., in \cite{bauer2021ripser}.  The incorporation of such optimizations into the multiparameter setting is an important direction for future research.

The approach to minimal presentation computation described above makes essential use of special structure in the 2-parameter setting. 
G\" afvert's recent Ph.D. thesis \cite{gafvert2021topological} provides a preliminary account of work with Bender and Lesnick on extending the approach of \cite{lesnick2019computing,kerber2021fast} to multiple parameters, drawing on Faug\'ere's well-known F5 algorithm for computing Gr\" obner bases \cite{faugere_F5_2002} and on ideas from the Ph.D. thesis of Skryzalin \cite{skryzalin2016numeric}.

\section{The Representation Theory of Multiparameter Persistence}
\label{sec:multidbarcode}

In this section, we use quiver representation theory to precisely describe the algebraic complexity of bipersistence modules indexed by the grid $[m]\times [n]$, for all $m, n\geq 1$. We begin in \cref{sec:finite-tame-wild} with a discussion of the relevant representation theory.  Then in \cref{Sec:Rep_Theory_of_MPMs}, we explain how this theory applies to the case of $[m]\times [n]$-modules.  \cref{sec:homreal,Sec:Hom_Dim_0} consider the homology modules of filtrations.  \cref{subsec:decomp-algo} briefly discusses the computation of direct sum decompositions.

\subsection{Posets of Finite, Tame, and Wild Type}
\label{sec:finite-tame-wild}

Drozd's theorem, a foundational result from quiver representation theory, provides a classification of finite posets into into three types (\emph{finite}, \emph{tame}, and \emph{wild}), according to the complexity of their spaces of persistence modules.  In what follows, we provide a brief introduction to these ideas.  We refer the reader to \cite{escolar2016persistence} for an introduction to quiver representation theory catered to a TDA audience, or to \cite{barot2015introduction,simson2007elements} for a thorough treatment.

We say that a finite poset $P$ is of \emph{finite representation type} if there exists a finite number of indecomposable $P$-modules up to isomorphism. 
\begin{example}
A poset of the form $\bullet \leftrightarrow \bullet \leftrightarrow  \cdots \leftrightarrow \bullet$, where each $\leftrightarrow$ denotes a single arrow pointing either to the left or to the right, is of finite representation type; the indecomposables are precisely the interval modules.  This a special case of \emph{Gabriel's theorem} \cite{gabriel1972unzerlegbare}; see \cite{ringel2016representation} for a direct and elementary proof.
\label{ex:AN}
\end{example}

\begin{example}
The poset $P=[2]\times [2]$ is also of finite type, and the indecomposables are precisely the interval modules.  
\label{ex:square}
\end{example}

\begin{example}
The poset $P=[3]\times [2]$ is of finite  type, but the indecomposables do not correspond to intervals.  For example, the following $P$-module is indecomposable:
\[\begin{tikzcd} k \arrow{r}{[0,1]^T} & k^2\arrow{r}{[1,1]} & k  \ \\ 0\ar[u]\ar[r] & k\ar[r, "="]\arrow{u}{[1,0]^T} & k  \ar[u, "="]\end{tikzcd}\]
See \cite{escolar2016persistence} for a discussion on how such indecomposables can be used to infer common topological structure in pairs of data sets, e.g. the atomic arrangements of silica glasses. 
\end{example}

The following gives an example of a poset $P$ that is \emph{not} of finite type. 

\begin{example}\label{Ex:4-star}
Let $P$ be the \emph{$4$-star} poset: 
\[
\begin{tikzcd} 
\bullet \ar[rr] & &  \bullet & & \bullet \ar[ll] \\
\bullet \ar[rru] & &  & & \bullet \ar[llu]  
\end{tikzcd}
\]
Let $n$ be a positive integer, $\lambda\in k$, and $J(\lambda)$ the $n\times n$ Jordan block over $k$:
\[ J(\lambda) = \begin{bmatrix} \lambda & 1 & 0 & \cdots & 0 \\ 0 & \lambda  & 1 & \cdots & 0 \\ \vdots & \vdots & \vdots & \ddots & \vdots \\ 0 & 0 & 0 & \cdots & \lambda\end{bmatrix}\]
Define $M^\lambda$ to be the $P$-module
\[
\begin{tikzcd}[ampersand replacement=\&]
k^n \arrow{rr}{A} \& \& k^n\oplus k^n \& \& k^n \arrow{ll}[swap]{C} \\
k^n   \arrow[swap]{urr}{B}\& \& \& \&  k^n  \ar{ull}{D}
\end{tikzcd}
\]
where the block matrices are $A=\begin{pmatrix} I_n \\ \textbf{0} \end{pmatrix}$, $B=\begin{pmatrix} \textbf{0}   \\I_n \end{pmatrix}$, $C=\begin{pmatrix} I_n\\ I_n  \end{pmatrix}$ and $D=\begin{pmatrix} I_n  \\ J(\lambda)  \end{pmatrix}$, and $I_n$ is the $n\times n$ identity matrix. One can verify that $M^\lambda$ is indecomposable and that $M^\lambda \not\cong M^{\lambda'}$ for $\lambda\neq \lambda'$.  In fact, for $k$ is algebraically closed and $n$ fixed, all but a finite number of indecomposables appear in one of finitely many such $1$-parameter families; see \cite{nazarova1973representations} for a complete classification of the indecomposables. Furthermore, if we work over $F_q$, the finite field with $q$ elements, then the number of isomorphism classes of indecomposables is a polynomial in $q$ of degree $n$ \cite{kac1982infinite}; see Remark \ref{rem:indec} below. 
\label{ex:tame-rep}
\end{example}

The poset in the previous example is of \emph{tame type}. Roughly speaking, this means that once the vector spaces dimensions have been fixed, all but finitely many of the indecomposables appear in one of finitely many ``$1$-parameter families'' of indecomposables. We refer to \cite[Section 9.5]{barot2015introduction} for a precise definition. 

We next observe that a small modification of the previous example can lead to a dramatic increase in complexity. 

\begin{example}\label{Ex:Five_Star}
Let $\lambda, \gamma \in k$, and consider the following two-parameter family of persistence modules indexed by the \emph{$5$-star poset}:
\begin{equation}
 M^{\lambda, \gamma} =\qquad \begin{tikzcd}[ampersand replacement=\&]
k^n \arrow{rr}{A} \& \&  k^n\oplus k^n \& \& k^n \arrow{ll}{C}  \\
k^n \arrow[swap]{rru}{B}\& \&   \& \& k^n \arrow{llu}{D}  \\
\& \& k^n \arrow{uu}{E}   \& \& 
\end{tikzcd}
\label{eq:wild}
\end{equation}
where $A$, $B$, $C$ and $D$ are given in Example \ref{ex:tame-rep} and  $E= \begin{pmatrix} I_n  \\ J(\gamma)  \end{pmatrix}$. Then, one can check that $M^{\lambda, \gamma}$ is indecomposable and that $M^{\lambda, \gamma} \not\cong M^{\lambda', \gamma'}$ for $(\lambda, \gamma) \neq (\lambda', \gamma')$. 

In contrast with Example \ref{ex:tame-rep}, however, the collection of such two-parameter families does not come close to describing all isomorphism classes of indecomposables. For instance, we can explicitly describe an $(n-1)$-parameter family of non-isomorphic indecomposables as follows: Fix $\lambda=0$ and replace $J(\gamma)$ in the matrix $E$ by the $n\times n$ matrix $K(\vec\alpha)$ whose entries on the super-diagonal are given by a vector $\vec{\alpha} \in  k^{n-1}$, with all other entries of $K(\vec\alpha)$ set to 0.  It can be shown that different choices of $\vec\alpha$ yield non-isomorphic indecomposables.

To give a further sense of the increase in complexity when passing from the 4-star poset to the 5-star poset, let $I(n,q)$ denote the number of isomorphism classes of indecomposables over $F_q$ for which the dimensions of the vector spaces are as in \cref{eq:wild}. Then $I(n,q)$ is a polynomial in $q$ of degree $1+n^2$ \cite{kac1982infinite}. 
\label{ex:wild}
\end{example}

\begin{remark}
More generally, Kac \cite{kac1982infinite} proved that if $\bar{d}$ is the \emph{dimension vector} of an indecomposable quiver representation over $F_q$, then the total number of indecomposables with dimension vector $\bar{d}$ is a polynomial $f\in \mathbb{Q}[q]$ with 
\[
\text{degree}(f) = \begin{cases}
1 - q_Q(\bar{d}) & \text{ if } \quad q_Q(\bar{d}) <0,\\
\gcd(\bar{d}) &\text{ if } \quad q_Q(\bar{d}) = 0,\\
0 & \text{ if } \quad q_Q(\bar{d}) >0, 
\end{cases}
\]
where $q_Q$ is the \emph{Tits form}. A simple computation of the Tits form recovers the degrees of the polynomials mentioned in Examples \ref{Ex:4-star} and \ref{Ex:Five_Star}. Computing the lower-order coefficients of the polynomial is difficult, in general; explicit examples for small quivers can be found in \cite{hua2000counting}. 

\label{rem:indec}
\end{remark}

The 5-star poset is an example of a poset of \emph{wild type}.  There are several equivalent definitions of wild type \cite[Chapter XIX, Theorem 1.11]{simson2007elements}.  We give three, treating the first two informally:

\begin{definition}\label{Def:Wild}
A finite poset $P$ is of \emph{wild type} if $P$ satisfies any of the following (equivalent) conditions:
\begin{itemize}
\item[(i)] There exists a 2-parameter family of indecomposable $P$-modules,
\item[(ii)] There exists an $n$-parameter family of indecomposable $P$-modules for all $n\in \mathbb N$,
\item[(iii)] For every finite poset $Q$, there exists an exact functor $F\colon \vect^Q\to \vect^{P}$ such that
\begin{enumerate}
\item if $M$ is indecomposable, then so is $F(M)$,
\item $F(M)\cong F(N)$ if and only if $M\cong N$, 
\item $F(cf + dg) = cF(f) + dF(g)$ for all $f, g\colon M\to N$ and $c, d\in k$. 
\end{enumerate}
\end{itemize}
\end{definition}

Thus, if $P$ is wild, a full classification of the indecomposable $P$-modules would yield a classification of the indecomposable $Q$-modules for any finite poset $Q$.  Obtaining such a classification is thought to be a hopelessly difficult task. 

The following ``trichotomy'' result is remarkable and surprising.

\begin{theorem}[Drozd \cite{drozd1980tame}]\label{Thm:Drozd_Trichotomy}
Over an algebraically closed field $k$, any finite poset is of finite, tame, or wild representation type. 
\end{theorem}

\subsection{The Representation Theory of $[m]\times [n]$}\label{Sec:Rep_Theory_of_MPMs} It turns out that the representation type of a commutative grid
is wild even for small grid sizes.  To see this, observe that any module over the 5-star poset embeds as multiparameter persistence
module over $[5]\times [5]$:
\[
\begin{tikzcd}
k^n \arrow{r}{A} & k^n\oplus k^n \arrow{r}{=} & k^n\oplus k^n \arrow{r}{=} & k^n\oplus k^n \arrow{r}{=} & k^n\oplus k^n \\
0 \ar[u]\arrow{r} & k^n \arrow{r}{B} \arrow{u}{B}& k^n\oplus k^n \arrow{u}{=}\arrow{r}{=} & k^n\oplus k^n \arrow{u}{=}\arrow{r}{=} & k^n\oplus k^n \arrow{u}{=}\\
0\ar[u]\arrow{r} &  0\ar[u]\arrow{r} \arrow{r}&k^n \arrow{u}{C}\arrow{r}{C} \ar[u] & k^n\oplus k^n \arrow{u}{=} \arrow{r}{=} & k^n\oplus k^n\arrow{u}{=} \\
0\ar[u]\arrow{r}&  0\ar[u]\ar[r] &  0\ar[u] \ar[r] & k^n \arrow{r}{D} \arrow{u}{D} & k^n\oplus k^n \arrow{u}{=}\\
0\ar[u]\arrow{r}& 0\ar[u] \ar[r] & 0\ar[u] \ar[r] & 0\ar[u]\arrow{r} \arrow{r} & k^n\arrow{u}{E}
\end{tikzcd}
\]
In particular, $[m]\times [n]$ must be of wild type for $m,n\geq 5$. The following theorem gives a complete description of the representation type of finite grids of varying sizes. 

\begin{theorem}[{\cite[Theorem 5]{MR1808681}, \cite[Theorem 2.5]{MR1273693}}]
    \label{thm.types_for_grids}
    The poset $[m]\times [n]$ is of finite type if
    \begin{itemize}
        \item $m = 1$ or $n = 1$,
        \item $(m,n) \in \{ (2,2), (2,3), (2,4), (3,2), (4,2)\}$.
    \end{itemize}
    It is of tame  type if
    \begin{itemize}
        \item $(m,n) \in \{ (2,5), (3,3), (5,2)\}$.
    \end{itemize}
    It is of wild type in all other cases.
\end{theorem}
\subsection{Homology in Dimension $\geq 1$}
\label{sec:homreal}
The previous discussion shows that the representation theory becomes very complicated even for small grids.  Moreover, Carlsson and Zomorodian have observed that, when working over a prime field $k$, one can realize any persistence module over a finite grid by applying $H_1(-; k)$ to a cellular or simplicial sublevel filtration.  As an example, \cref{fig:bifiltered-circle}  illustrates how one can realize the following indecomposable persistence module, previously considered in Example \ref{Ex:Not_Good}:
\begin{equation}\begin{tikzcd}
k\arrow{r}{=} & k\ar[r] & 0 \\
k\arrow{r}{[1,0]^T} \arrow{u}{=} & k^2\arrow{r}{[1,1]}\arrow{u}{[1,0]} & k\ar[u]\\
0\ar[r]\ar[u] & k\arrow{u}{[0,1]^T}\arrow{r}{=} & k\ar[u,"="]
\end{tikzcd}
\label{eq:not-rank}
\end{equation}
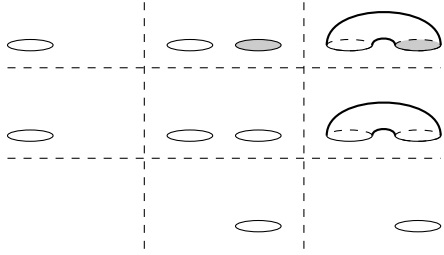
\begin{figure}
\centering
\begin{tikzpicture}[scale=0.3]
\begin{scope}
\end{scope}

\begin{scope}[yshift=4cm]
\draw (0,0) arc[x radius=1, y radius=0.25, start angle=0, end angle=-360];
\end{scope}

\begin{scope}[yshift=8cm]
\draw (0,0) arc[x radius=1, y radius=0.25, start angle=0, end angle=-360];
\end{scope}

\begin{scope}[xshift=7cm]
\draw (3,0) arc[x radius=1, y radius=0.25, start angle=0, end angle=-360];
\end{scope}

\begin{scope}[xshift=7cm, yshift=4cm]
%

\draw (0,0) arc[x radius=1, y radius=0.25, start angle=0, end angle=-360];
\draw (3,0) arc[x radius=1, y radius=0.25, start angle=0, end angle=-360];

\end{scope}

\begin{scope}[xshift=7cm, yshift=8cm]
%
%
\draw[fill=black!20] (3,0) arc[x radius=1, y radius=0.25, start angle=0, end angle=-360];
\draw (0,0) arc[x radius=1, y radius=0.25, start angle=0, end angle=-360];
\end{scope}

\begin{scope}[xshift=14cm]
\draw (3,0) arc[x radius=1, y radius=0.25, start angle=0, end angle=-360];
\end{scope}

\begin{scope}[xshift=14cm, yshift=4cm]
\draw[dashed] (3,0) arc[x radius=1, y radius=0.25, start angle=0, end angle=-360];
\draw (3,0) arc[x radius=1, y radius=0.25, start angle=0, end angle=-180];

\draw (0,0) arc[x radius=1, y radius=0.25, start angle=0, end angle=-180];
\draw [dashed] (0,0) arc[x radius=1, y radius=0.25, start angle=0, end angle=180];

\draw[-,thick]  (0,0) to [out=90,in=90] (1,0);
\draw[-,thick]  (-2,0) to [out=90,in=90] (3,0);
\end{scope}

\begin{scope}[xshift=14cm, yshift=8cm]
\draw[fill=black!20, dashed] (3,0) arc[x radius=1, y radius=0.25, start angle=0, end angle=-360];
\draw (3,0) arc[x radius=1, y radius=0.25, start angle=0, end angle=-180];

\draw[dashed] (0,0) arc[x radius=1, y radius=0.25, start angle=0, end angle=-360];
\draw (0,0) arc[x radius=1, y radius=0.25, start angle=0, end angle=-180];

\draw[-,thick]  (0,0) to [out=90,in=90] (1,0);
\draw[-,thick]  (-2,0) to [out=90,in=90] (3,0);
\end{scope}
\draw[dashed] (4,-1) -- (4,10);
\draw[dashed] (11,-1) -- (11,10);
\draw[dashed](-2, 3) -- (17,3);
\draw[dashed](-2, 7) -- (17,7);

\end{tikzpicture}
\caption{Applying $H_1(-;k)$ to the above bifiltration yields the indecomposable bipersistence module in \cref{eq:not-rank}.}
\label{fig:bifiltered-circle}
\end{figure}

It is shown in \cite{buchet2018realizations} that infinite families of indecomposables on a finite grid can arise as the homology modules of function-Rips bifiltrations.

\subsection{Homology in Dimension 0}\label{Sec:Hom_Dim_0}

Several bifiltrations that appear naturally in multiparameter persistence have the property that their $0$-th homology modules are surjective in at least one parameter direction. As an example, for $\gamma\colon P\to \R$ as in Definition \ref{def:superlevel-rips},
\[H_0(\SupVR(\gamma)_{(a,r), (a,s)})\colon H_0(\SupVR(\gamma)_{a,r})\to H_0(\SupVR(\gamma)_{a, s})\] is surjective for all $r\leq s$. 

\begin{remark}
$H_0$ factors as $H_0(-;k)\cong {\rm Free}\circ \pi_0$, where $\pi_0$ denotes path-components, and ${\rm Free}\colon \Set\to \Vect_k$ is the free functor sending a set $X$ to the free $k$-vector space with basis $X$. The fact that $H_0(\SupVR(\gamma))_{(a,r), (a,s)}$ is an epimorphism for $r\leq s$ is a direct consequence of this factorization. See \cite{chacholski2017combinatorial, brodzki2020complexity} for a further discussion of persistence modules factoring through the category of sets. 
\end{remark}

A natural question is to what extent such additional structure reduces the algebraic complexity. For an instance of a result in this direction, one can show that that if $M$ is an $[m]\times [n]$-module with surjections in one direction and injections in the other, then $M$ is interval-decomposable; see \cref{thm:block}.  On the other hand, the two-parameter families of persistence modules from Example \ref{ex:wild} all embed as $[5]\times [5]$-modules with all morphisms surjective. This suggests that the reduction in complexity can only be rather limited. The following theorem confirms that. 

\begin{theorem}[\cite{bauer2020cotorsion}]\label{Thm:Horizontal_Surjections}
The category of $[m]\times [n]$-modules whose horizontal maps are surjective is equivalent to the category of $[m]\times [n-1]$-modules, up to a finite, explicit list of indecomposables. 
\end{theorem}
An analogue of \cref{Thm:Horizontal_Surjections} for the case where all morphisms are surjective is also given in \cite{bauer2020cotorsion}. Furthermore,  \cref{Thm:Horizontal_Surjections} and \cref{thm.types_for_grids} together immediately imply the following:

\begin{corollary}
For $[m]\times [n]$-modules whose horizontal maps are surjective, the cases of finite type are:
\begin{itemize}
\item $m=1$ or $n\leq 2$ (all modules are interval-decomposable)
\item $(m,n)\in \{(2,3), (2,4), (2,5), (3,3), (4,3)\}$,
\end{itemize}
and the tame cases are: $(m,n)\in \{(2,6), (3,4), (5,3)\}$. 
\end{corollary}

\subsection{Computing Decompositions}
\label{subsec:decomp-algo}
While explicitly classifying the indecomposable multiparameter persistence modules is hopeless, one can still compute the indecomposable summands of any given persistence module.  General purpose polynomial-time algorithms for decomposing modules over finite-dimensional algebras exist in the literature and implementations are readily available; see e.g,. the MeatAxe algorithm \cite{holt1998meataxe}, which computes a decomposition of a $P$-module in time $\tilde{O}\left(\left(\sum_{p\in P} \dim M_p\right)^6\right)$, and is implemented in GAP \cite{GAP4}.  However, such algorithms are impractical in the context of TDA.

Dey and Xin \cite{dey2019generalized} have recently provided a decomposition algorithm tailored to multiparameter persistence. Their algorithm takes as input a presentation matrix of a $d$-parameter persistence module $M$.  Under the assumption that no two rows or columns have the same label, the algorithm outputs a decomposition of $M$ into indecomposables.  When this assumption does not hold, the algorithm still provides a decomposition of $M$ into summands, but the summands may not be indecomposabale.  For an input presentation matrix with a total of $N$ rows and columns, the algorithm runs in time $O(N^{(d-1)(2\omega +1)})$, where $\omega < 2.373$ is the exponent for matrix multiplication. 
While this work represents significant progress on the efficient decomposition of MPH modules, new ideas are needed for how to use such decompositions in the context of data analysis.

\section{Signed Barcodes}
\label{sec:rankinv}
Recall that the rank invariant $\Rk M$ of a persistence module $M$ is the function $(s,t)\mapsto\Rk(M_s\to M_t)$, where $s\leq t$.
 As noted in \cref{subsec:invariants}, the rank invariant is a complete invariant of a 1-parameter persistence module.  Thus, the barcode and rank invariant of a 1-parameter persistence module determine each other.  Though the rank invariant is incomplete outside of the 1-parameter setting, one may hope that the correspondence between barcodes and rank invariants can be extended to the multiparameter setting.  As we will explain, this is indeed possible if one considers \emph{signed barcodes}.  

In this section we will consider three signed barcode constructions appearing in the recent TDA literature, which we will call the \emph{rectangle}, \emph{interval}, and \emph{hook} barcodes.  The rectangle barcode \cite{botnan2021signed} and hook barcodes \cite{botnan2021signed,botnan2022bottleneck} are both equivalent to the rank invariant; the former can be constructed via M\"obius inversion, while the latter is constructed via relative homological algebra.  The hook barcode has the advantage of being stable in a stronger sense than the rectangle barcode.  The interval barcode (also known as the \emph{generalized persistence diagram}) \cite{botnan2021signed,kim2018generalized,asashiba2019approximation,kim2018generalized} is an analogue of the rectangle barcode which is equivalent to the generalized rank invariant (Definition \ref{def:generalized-rank_invariant}).  

\subsection{Rectangle Barcodes via Decomposition of the Rank Invariant}\label{Sec:Rectangle_Barcode}

In Example \ref{Ex:Not_Good}, we saw a module $M$ whose rank invariant is not equal to the rank invariant of any interval-decomposable module.  On the other hand, the rank invariant of $M$ can be expressed as the \emph{difference} between the rank invariants of two interval-decomposable modules:

\begin{equation}
\centering
\includegraphics[scale=0.8]{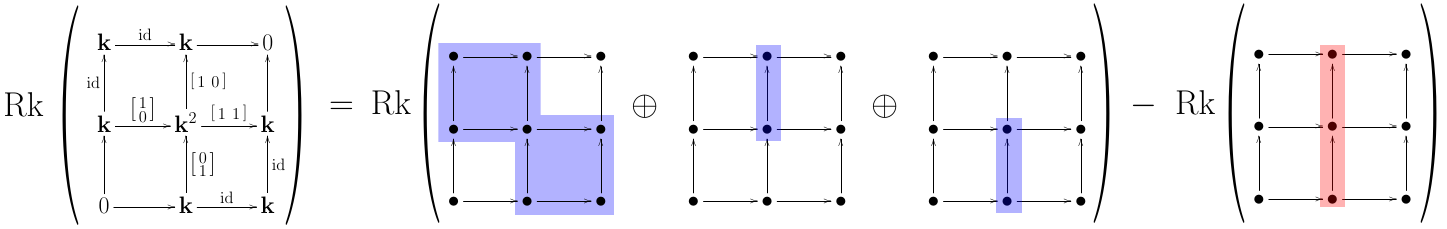}
\label{eq:rank-decomp-first}
\end{equation}

Such a decomposition need not be unique: 
\begin{equation}
\centering
\includegraphics[scale=0.8]{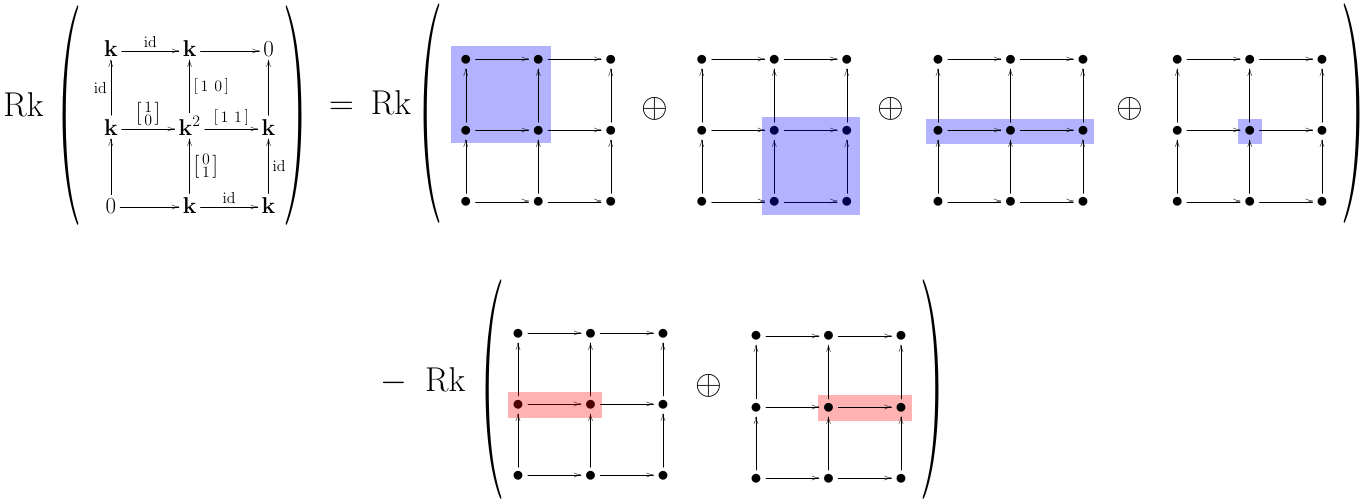}
\label{eq:rank-decomp}
\end{equation}

It turns out, however, that the decomposition is unique if we restrict to rectangles. 

\begin{definition}\label{Def:Rectangle}
A \emph{rectangle} in a product of totally ordered sets $T_1\times \cdots\times T_n$ is a poset of the form $I_1\times \cdots\times I_n$, where each $I_j$ is an interval in $T_j$.  
\end{definition}

\begin{theorem}[\cite{botnan2021signed}]\label{thm:R-S_multi-param}
Let $P=T_1\times \cdots \times T_n$, where each $T_i\subseteq \R$, and let $M$ be a finitely presented $P$-module. There exists a \emph{minimal} pair $(\Rec, \Sec)$ of finite multisets of rectangles in $P$ satisfying 
\begin{equation} \Rank M = \Rank\left(\bigoplus_{I\in \Rec} k_I\right)-\Rank\left(\bigoplus_{J\in \Sec} k_J\right);
\label{eq:RS}\end{equation}
here minimal means that if $(\Rec', \Sec')$ is any other such pair satisfying \cref{eq:RS}, then $\Rec\subseteq \Rec'$ and $\Sec\subseteq \Sec'$. 
\end{theorem}
We shall refer to the pair $(\Rec, \Sec)$ as the \emph{rectangle barcode} of $\Rank M$.  (It is called a \emph{minimal rank decomposition} in \cite{botnan2021signed}.)   Note that if $d=1$, then by \cref{Thm:ordinary_Structure} and the uniqueness of the signed barcode, we must have that $\Sec = \emptyset$.  

The rectangle barcode can be constructed by taking a M\" obius inversion of the restriction of the rank invariant to a finite grid \cite{botnan2021signed}.

\begin{remark}\label{Rem:Other_Signed_Barcodes}
The idea of using M\"obius inversion to associate a signed barcode to the rank invariant was introduced by Patel~\cite{patel2018generalized}, who worked in the setting of 1-parameter persistence modules taking values in certain symmetric monoidal categories. Subsequently, M\"obius inversions have been applied to persistence modules in several  settings \cite{mccleary2020edit, asashiba2019approximation,kim2018generalized,betthauser2022graded}. See \cite[Section 1]{botnan2021signed} for a detailed discussion of the connection between the various approaches.
\end{remark}

The 1-parameter slices of rectangle barcodes are stable, in the sense of the following theorem: 
\begin{theorem}[\cite{botnan2021signed}]\label{th:matching-distance_stability}
Given signed barcodes $(\Rec, \Sec)$ and $(\Rec', \Sec')$ of finitely presented  $\R^d$-persistence modules $M$ and $M'$, we have:
  \[
  d_{\mathrm{match}}\left(\bigoplus_{I\in \Rec} k_I \oplus \bigoplus_{J'\in \Sec'} k_{J'}, \bigoplus_{I'\in \Rec'} k_{I'} \oplus \bigoplus_{J\in \Sec} k_J\right) \leq d_{\mathrm{match}}(M, M').
  \]
\end{theorem}
One may ask whether \cref{th:matching-distance_stability} still holds if one replaces $d_{\mathrm{match}}$ with the generalized bottleneck distance $d_B$ (Definition \ref{Def:bottleneck_general}).  The answer is no, even up to a constant factor \cite[Proposition 3.1]{botnan2022bottleneck}.   In contrast, the hook barcode, which we will introduce in \cref{Sec:Rel_Hom_Alg}, is Lipschitz stable with respect to $d_B$.  

\begin{figure}[tb]
\centering
\begin{subfigure}{.4\linewidth}
\centering
\begin{tikzpicture}[xscale=1.5,yscale=1.5]
\fill [gray,opacity=0,draw]
  (-.3,-.3) --
  (.3,-.3) {[rounded corners=10] --
  (.3,1.3)} --
  (-0.3,1.3) --
  cycle
  {};
  \fill [gray,opacity=0,draw]
  (.7,2.3) {[rounded corners=10] --
  (.7,.7)}  --
  (2.3,.7) --
  (2.3,2.3) --
  cycle
  {};
\path node[opacity=.2]   (v00) at (0,0)   {$\bullet$}   
      node[opacity=.2]    (v10) at (1,0) {$\bullet$} 
      node[opacity=.2]    (v20) at (2,0) {$\bullet$} 
      node[opacity=.2]    (v01) at (0,1)  {$\bullet$} 
       node[opacity=.2]    (v01b) at (0,.925) {} 
       node[opacity=.2]    (v01a) at (0,1.075) {} 
      node[opacity=.2]    (v11) at (1,1) {} 
       node[opacity=.2]    (v11b) at (1,.925) {} 
       node[opacity=.2]    (v11a) at (1,1.075) {} 
      node[opacity=.2]    (v21) at (2,1) {$\bullet$} 
     node[opacity=.2]    (v21b) at (2,.925) {} 
       node[opacity=.2]    (v21a) at (2,1.075) {} 
      node[opacity=.2]    (v02) at (0,2) {$\bullet$} 
      node[opacity=.2]    (v12) at (1,2) {$\bullet$} 
      node[opacity=.2]    (v22) at (2,2)  {$\bullet$} ;
     { 
        \draw[opacity=.2]   (v00) -- node[above]{}      (v10); 
     \draw[opacity=.2]   (v10) -- node[above]{}      (v20); 
     
     \draw[opacity=.2]   (v01) -- node[above]{}      (v11); 
     \draw[opacity=.2]   (v11) -- node[above]{}      (v21); 
     
     \draw[opacity=.2]   (v02) -- node[above]{}      (v12); 
     \draw[opacity=.2]   (v12) -- node[above]{}      (v22); 
     
      \draw[opacity=.2]   (v00) -- node[above]{}      (v01); 
     \draw[opacity=.2]   (v01) -- node[above]{}      (v02); 
     
     \draw[opacity=.2]   (v10) -- node[above]{}      (v11); 
     \draw[shorten <=2pt,opacity=.2]   (v11) -- node[above]{}      (v12); 
     
     \draw[opacity=.2]   (v20) --      (v21); 
     \draw[opacity=.2]   (v21) --    (v22); 
     
     \draw[line width=1pt,blue,opacity=.5]   (v01) --    (v12); 
     \draw[line width=1pt,blue,opacity=.5]   (v10) --    (v21); 
     \draw[line width=1pt,blue,opacity=.5]   (v01a) --    (v21a); 
     \draw[blue,opacity=.5,fill] (1,1) circle [radius=.05];
     \draw[line width=1pt,red,opacity=.5]   (v01b) --    (v11b);
     \draw[line width=1pt,red,opacity=.5]   (v11b) --    (v21b);
}     
\end{tikzpicture}
\caption{}
\end{subfigure}
\begin{subfigure}{.4\textwidth}
\centering
\begin{tikzpicture}[xscale=1.5,yscale=1.5]
\fill [gray,opacity=.2,draw]
  (-.3,-.3) --
  (.3,-.3) {[rounded corners=10] --
  (.3,1.3)} --
  (-0.3,1.3) --
  cycle
  {};
  \fill [gray,opacity=.2,draw]
  (.7,2.3) {[rounded corners=10] --
  (.7,.7)}  --
  (2.3,.7) --
  (2.3,2.3) --
  cycle
  {};
\path node[opacity=.2]   (v00) at (0,0)   {$\bullet$}   
      node[opacity=.2]    (v10) at (1,0) {$\bullet$} 
      node[opacity=.2]    (v20) at (2,0) {$\bullet$} 
      node[opacity=.2]    (v01) at (0,1)  {$\bullet$} 
       node[opacity=.2]    (v01b) at (0,.925) {} 
       node[opacity=.2]    (v01a) at (0,1.075) {} 
      node[opacity=.2]    (v11) at (1,1) {} 
       node[opacity=.2]    (v11b) at (1,.925) {} 
       node[opacity=.2]    (v11a) at (1,1.075) {} 
      node[opacity=.2]    (v21) at (2,1) {$\bullet$} 
     node[opacity=.2]    (v21b) at (2,.925) {} 
       node[opacity=.2]    (v21a) at (2,1.075) {} 
      node[opacity=.2]    (v02) at (0,2) {$\bullet$} 
      node[opacity=.2]    (v12) at (1,2) {$\bullet$} 
      node[opacity=.2]    (v22) at (2,2)  {$\bullet$}
      node[above right]  at (-.35,-.35) {\tiny $(0,1)^-$}
      node  at (1.5,1.5) {\tiny $(1,1)^+$}
       ;
     { 
      \draw[opacity=.2]   (v00) -- node[above]{}      (v10); 
     \draw[opacity=.2]   (v10) -- node[above]{}      (v20); 
     
     \draw[-{Latex[length=1.5mm, width=1.5mm]},thick,green!50!black]   (v01) -- node[above]{}      (v11); 
     \draw[opacity=.2]   (v11) -- node[above]{}      (v21); 
     
     \draw[opacity=.2]   (v02) -- node[above]{}      (v12); 
     \draw[opacity=.2]   (v12) -- node[above]{}      (v22); 
     
      \draw[opacity=.2]   (v00) -- node[above]{}      (v01); 
     \draw[opacity=.2]   (v01) -- node[above]{}      (v02); 
     
     \draw[opacity=.2]   (v10) -- node[above]{}      (v11); 
     \draw[shorten <=2pt,opacity=.2]   (v11) -- node[above]{}      (v12); 
     
     \draw[opacity=.2]   (v20) --      (v21); 
     \draw[opacity=.2]   (v21) --    (v22); 
     
     \draw[line width=1pt,blue,opacity=.5]   (v01) --    (v12); 
     \draw[line width=1pt,blue,opacity=.5]   (v10) --    (v21); 
     \draw[line width=1pt,blue,opacity=.5]   (v01a) --    (v21a); 
     \draw[blue,opacity=.5,fill] (1,1) circle [radius=.05];
     \draw[line width=1pt,red,opacity=.5]   (v01b) --    (v11b);
     \draw[line width=1pt,red,opacity=.5]   (v11b) --    (v21b);
}     
\end{tikzpicture}
\caption{}
\end{subfigure}  \caption{(a) Visualization of the signed barcode of \cref{eq:rank-decomp}.  Rectangles with positive and negative sign are represented using blue and red bars, respectively.  (b) The rank of the map $M_{(0,1)}\to M_{(1,1)}$ (represented by the green arrow) is 1;  This is the number of blue bars intersecting both gray regions minus the number of red bars intersecting both regions.}
  \label{fig:barcode_indec2_grid}
\end{figure}
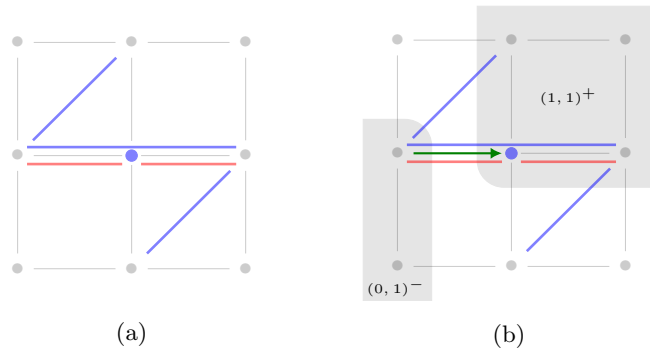

\paragraph{Visualization} Visualizing the rectangle barcode of a bipersistence module by directly plotting the rectangles can be messy for examples of realistic size, since many rectangles may overlap.  To obtain a cleaner visualization, \cite{botnan2021signed} proposes to instead represent each rectangle via a line segment (``bar'') connecting its infimum to its supremum, as illustrated in \cref{fig:barcode_indec2_grid}\,(a).  Examples involving larger persistence modules can be found in \cite{botnan2021signed}.  

One can read the rank invariant off of such a visualization: $\Rk(M_s\to M_t)$ is the number of bars in $\Rec$ intersecting both the sets $s^-=\{u\in P \colon u\leq s\}$ and $t^+=\{u\in P \colon u\geq t\},$
  minus the number of bars in $\Sec$ intersecting both sets; see \cref{fig:barcode_indec2_grid}\,(b).

\paragraph{Computation}\label{Sec:Compuitng_Signed_Barcode} The rectangle barcode of a p.f.d. persistence module $M$ indexed over a finite grid $\grid = \prod_{i=1}^d [1, n_i]$  can be computed by means of a simple inclusion-exclusion formula. Let $\alpha_{[s,t]}$ denote the multiplicity of the rectangle  \[[s,t] = \{u \in \R^d \mid s\leq u\leq t\},\] in the rectangle barcode $(\Rec, \Sec)$. That is, if $\alpha_{[s,t]}>0$ then $[s,t]$ appears with multiplicity $\alpha_{[s,t]}$ in $\Rec$, and if $\alpha_{[s,t]}<0$ then $[s,t]$ appears with multiplicity $-\alpha_{[s,t]}$ in $\Sec$. Then we have:
\begin{equation}\label{eq:incl_excl_rect}
\alpha_{[s,t]} = \sum_{\begin{smallmatrix}s'\leq s\\\|s'-s\|_\infty \leq 1\end{smallmatrix}} \sum_{\begin{smallmatrix}t'\geq t\\\|t'-t\|_\infty \leq 1\end{smallmatrix}} (-1)^{\|s'-s\|_1 + \|t'-t\|_1}\,\Rank M(s',t'). 
\end{equation}
\begin{remark}
The case $d=1$ gives the well-known inclusion-exclusion formula relating the persistence diagram of a one-parameter persistence module to its rank invariant~\cite{cohen2007stability}. The case $d=2$ specializes to the inclusion-exclusion formula for computing the multiplicities of summands in rectangle-decomposable 2-parameter persistence modules~\cite{botnan_et_al:LIPIcs:2020:12180}.
\end{remark}

A simple inspection of the formula in \cref{eq:incl_excl_rect} reveals that computation is bounded in time $O\left( 2^{2d}\, \prod_{i=1}^d n_i^2\right)$, assuming constant-time access to the ranks $\Rank M(s', t')$ and constant-time arithmetic operations. When the module $M$ comes from a simplicial filtration over the grid~$\grid$ with $n = \max_{i} n_i$ simplices in total, the rank invariant itself can be pre-computed and stored, e.g., by naively computing the ranks $\Rank M(s,t)$ for each pair $s\leq t\in\grid$ independently, which takes $O(n^{2d+\omega})$ time in total, where $\omega<2.373$ is the exponent for matrix multiplication~\cite{milosavljevic2011zigzag}. Adding in the computation time for the signed barcode yields a bound in~$O(n^{2d+\omega} + (2n)^{2d})$. In the special case where $d=2$, assuming the filtration is $1$-critical (i.e. each simplex has a unique minimal index of appearance in the filtration), there is an $O(n^4)$-time algorithm to compute the rank invariant~\cite{morozov2008homological, botnan_et_al:LIPIcs:2020:12180}, and computing its rectangle barcode also takes $O(n^4)$ time.

\subsection{Interval Barcodes via Decomposition of the Generalized Rank Invariant}
While a smaller signed barcode using intervals of other shapes than rectangles is sometimes possible, it comes at the expense of uniqueness, as illustrated by \cref{eq:rank-decomp-first} and \cref{eq:rank-decomp}. It is however possible to retain uniqueness if one instead decomposes the generalized rank invariant (Definition \ref{def:generalized-rank_invariant}). 

Let $P$ be an arbitrary poset. We say a collection $\Int$ of intervals in $P$ is \emph{locally finite} if for any $I,K\in\Int$ the set $\{J\in\Int : I\subseteq J\subseteq K\}$ is finite.

\begin{theorem}[\cite{botnan2021signed}]
\label{thm:genarlized-rank_invariant}
Let $\Int$ be a locally finite collection of intervals in $P$, and let $M$ be a p.f.d. $P$-module.  Then there exists a pair $(\Rec, \Sec)$, where $\Rec$ and $\Sec$ are finite multi-sets of intervals in $\Int$, such that 
\begin{equation} \Rank_\Int M = \Rank_\Int\left(\bigoplus_{I\in \Rec} k_I\right)-\Rank_\Int\left(\bigoplus_{J\in \Sec} k_J\right),\label{eq:RS2}\end{equation}
and if $(\Rec', \Sec')$ is any other such pair satisfying \cref{eq:RS2}, then $\Rec\subseteq \Rec'$ and $\Sec\subseteq \Sec'$. 
\end{theorem}

\begin{remark}
When $\Int$ consists of all intervals in an essentially finite poset $P$, the \emph{generalized persistence diagram} of \cite{kim2018generalized}, constructed via M\"obius inversion of the generalized rank invariant, is essentially the same as the decomposition $(\Rec, \Sec)$ of \cref{thm:genarlized-rank_invariant}.
\end{remark}

\subsection{Hook Barcodes via Relative Homological Algebra}\label{Sec:Rel_Hom_Alg}
Several recent works use \emph{relative homological algebra} to study generalized persistence.  In brief, the idea is to identify an \emph{exact structure} for a given invariant, i.e., a distinguished class of exact sequences which respects the structure of the invariant.  There is a notion of \emph{relative minimal projective resolution} associated to any exact structure, which is used to define a signed barcode. This approach to generalized persistence was initiated for the rank invariant in \cite{botnan2021signed}, and in parallel for other invariants in \cite{blanchette_brustle_hanson_2022}; see also \cite{oudot2021stability,asashiba2022approximation, 
 chacholski2022effective,botnan2022bottleneck} for other recent work, and \cite{buhler2010exact} for an excellent introduction to exact categories. 

In the special case of the rank invariant, this approach yields the \emph{hook barcode}, a signed barcode which is bottleneck-stable with respect to the interleaving distance \cite{botnan2022bottleneck}.  We now explain this case in more detail.

\begin{definition}~
\begin{itemize}
\item[(i)] A short exact sequence \[0\to M\to M'\to M''\to 0\] of $\R^n$-persistence modules is \emph{rank exact} if $\Rk M' = \Rk M+\Rk M''$,
\item[(ii)] A long exact sequence \[\cdots \xrightarrow{f_{i+1}}  X_i \xrightarrow{f_i} X_{i-1} \xrightarrow{f_{i-1}}\cdots \] is \emph{rank exact} if $\im f_i \to X_{i-1} \to \im f_{i-1}$ is rank exact for all $i$,
\item[(iii)] An $\R^n$-persistence module $M$ is \emph{rank projective} if ${\rm Hom}(M,-)$ takes each rank exact sequence to an exact sequence of vector spaces.
\end{itemize}
\end{definition}

\begin{definition}\label{Def:Hook}
A \emph{hook} is an interval in $\R^n$ of the form $\{ p\in \R^n : r\leq p\not\geq s\}$ for $r<s\in \R^n\cup \{\infty\}$.  
\end{definition}

One can show that if $M\colon \R^n\to \Vect$ is finitely presented and rank projective, then $M$ is \emph{hook-decomposable}, i.e., $M$ is interval-decomposable and each interval in $\B{M}$ is a hook.

\begin{theorem}[\cite{botnan2021signed,botnan2022bottleneck}]\label{Thm:Rank_Exact_Res}
A finitely presented $\R^n$-persistence module $M$ admits a minimal rank projective resolution \[0\to X_k\to X_{k-1} \to \cdots \to X_1\to X_0\to M,\]
with each $X_i$ finitely presented (and hook-decomposable) and $k\leq 2n-2$.  
\end{theorem}
It follows from the theory of exact categories that a minimal rank projective resolution is unique up to isomorphism. 
Given a finitely presented $\R^n$-persistence module $M$, let $\beta_i(M) = \B{X_i},$ where $X_i$ is as in \cref{Thm:Rank_Exact_Res}.  

\begin{definition}\label{Def:Hook_Signed_Barcode}
We define the \emph{hook barcode} of $M$ to be the pair 
\[ \left(\beta_{2\mathbb{N}}(M), \beta_{2\mathbb{N}+1}(M)\right) := \left(\bigcup_{i \text{ even}}\beta_i(M), \bigcup_{i \text{ odd}} \beta_i(M)\right).\]
\end{definition}

\begin{remarks}~
\begin{itemize}
\item[(i)] By virtue of the exact structure, one has that \[\Rk M = \sum_{i=0}^k (-1)^i \Rk X_i = \sum_{i \text{ even}} \Rk X_i - \sum_{i \text{ odd}} \Rk X_i .\] This implies that we can interpret the hook barcode of $M$ as a decomposition of $\Rk M$, as for the rectangle barcode.
\item [(ii)]  In contrast to the rectangle barcode, the two components of the hook barcode may have intervals in common. Removing common hooks from the decomposition in (i) yields a minimal decomposition of the rank function in terms of hook modules.  However, the main stability property of the hook barcode (\cref{luis-stability} below) does not hold for this minimal decomposition \cite[Remark 3.2]{botnan2022bottleneck}.  
\end{itemize}
\end{remarks}
Like the rectangle barcode, the hook barcode can be visualized by a collection of line segments, as the following rank projective resolution illustrates. 
\begin{center}
\includegraphics[scale=0.28]{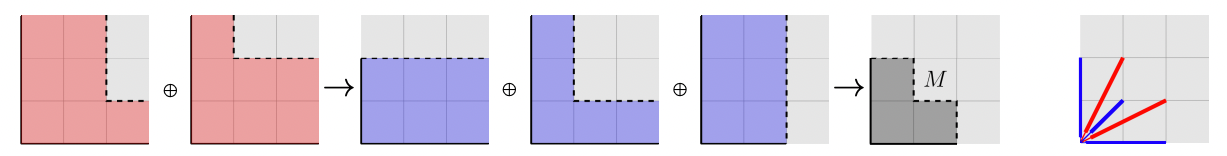}
\end{center}

The hook barcode is stable in the following sense:

\begin{theorem}[\cite{botnan2022bottleneck}]
\label{luis-stability}
If $M$ and $N$ are finitely presented $\R^n$-persistence modules, then
\[d_B\left(\beta_{2\mathbb{N}}(M)\cup \beta_{2\mathbb{N}+1}(N), \beta_{2\mathbb{N}}(N)\cup \beta_{2\mathbb{N}+1}(M)\right) \leq (2n-1)^2d_I(M,N).\]
\end{theorem}

The proof is an application of the algebraic stability of hook-decomposable modules (\cref{Thm:Generalized:Algebraic_Stability}\,(iv)).

\begin{remark}
Computation of minimal resolutions and \emph{relative Betti numbers} in exact categories is an area of active research. Recently, Chachólski et al. \cite{chacholski2022effective} provided a framework for computing relative Betti numbers of exact structures on upper semi-lattices. 
\end{remark}

\section{Local Conditions for Interval Decomposability}
\label{sec:block}
In this section, we introduce local conditions on a bipersistence module $M$ which force $M$ to be interval-decomposable, with each interval in $\B{M}$ of a particularly simple form.  We will see that such conditions arise naturally in the study of the \emph{interlevel  persistence} of $\R$-valued functions.

\subsection{Decomposition into Blocks and Rectangles}
Let $S$ and $T$ be totally ordered sets.
 
\begin{definition}\label{Def:Exactness_Conditions}
Given $P\subseteq S\times T$ and a persistence module $M\colon P\to \kVect$, 
\begin{itemize}
\item[(i)]
$M$ is \emph{middle-exact} \cite{botnan2020decomposition,cochoy2020decomposition,carlsson2009zigzag,carlsson2010zigzag} if 
\begin{equation}
 M_a \xrightarrow{M_{a,b}\oplus M_{a,c}} M_b\oplus M_c \xrightarrow{(M_{b,d} -M_{c,d})} M_d
\label{eq:middlex}
\end{equation}
is exact (i.e., the image of the first map equals the kernel of the second map) for all $a,b,c,d\in P$ of the form 
\[a=(x,y),\quad b=(x,y'),\quad c=(x',y),\quad d=(x', y').\] 
\item[(ii)]
$M$ is \emph{weakly-exact} \cite{botnan2020local,botnan_et_al:LIPIcs:2020:12180} if for all $a,b,c,d$ as above,
\begin{align*}
\Im(M_a \to M_d) &= \Im(M_b\to M_d)\cap \Im (M_c\to M_d),\\
\Ker(M_a\to M_d) &= \Ker(M_a\to M_b) + \ker(M_a\to M_c).
\end{align*} 
\end{itemize}
\end{definition}
One can check that if $M$ is middle-exact then it is weakly-exact.

Recall from Example \ref{ex:square} that a persistence module over a commutative square is interval-decomposable. A simple case-by-case inspection of the intervals gives the following characterization of middle-exact and weakly-exact persistence modules. 

\begin{proposition}~
\begin{itemize}
\item[(i)] $M$ is middle-exact if and only if its restriction to any square 
\[\begin{tikzcd}
M_c \ar[r] & M_d\\
M_a \ar[u]\ar[r] & M_b\ar[u]
\end{tikzcd}
\]
decomposes into interval modules supported on the set \[B:=\{\{a\},\{d\},\{a,b\},\{a,c\},\{b,d\},\{c,d\},\{a,b,c,d\}\}.\] 
\item[(i)] $M$ is weakly-exact if and only if its restriction to any square as above
decomposes into interval modules supported on the set \[B\cup \{\{b\},\{c\}\}.\] 
\end{itemize}
\end{proposition}

\begin{definition}\label{Def:Block}
An interval $I\subseteq S\times T$ is called a \emph{block} if $I$ can be written in one of the following ways:
\begin{enumerate}
\item $I=J_S\times J_T$ for downsets $J_S$ and $J_T$, 
\item $I=J_S\times J_T$ for upsets $J_S$ and $J_T$,
\item $I=J_S \times T$ for an interval $J_S$, 
\item $I=S\times J_T$ for an interval $J_T$. 
\end{enumerate}

\end{definition}
More generally, recall from Definition \ref{Def:Rectangle} that a rectangle in $S\times T$ is an interval of the form $I\times J$, where $I$ and $J$ are intervals in $S$ and $T$.  
We say that $M$ is \emph{block-decomposable} (respectively, \emph{rectangle-decomposable}) if $M$ is interval-decomposable and each interval of $\B{M}$ is a block (respectively, a rectangle). 

The following theorems tell us that a module is block-decomposable or rectangle-decomposable if and only if its restriction to each square is.

\begin{theorem}[\cite{botnan2020decomposition}]
A p.f.d. $S\times T$-module is middle-exact if and only if it 
 is block-decomposable.
\label{thm:block}
\end{theorem}

Given a totally ordered set $T$, a subset $T'\subset T$ is said to be \emph{coinitial} if for all $t\in T$ there exists $t'\in T'$ with $t'\leq t$.  For example, $\Z$ is a coinitial subset of $\R$.  

\begin{theorem}[\cite{botnan2020local}]
Suppose that $S$ and $T$ are totally ordered sets such that every interval in $S$ or $T$ admits a countable coinitial subset. Then a p.f.d. $S\times T$-module is weakly-exact if and only if it is rectangle-decomposable.
\label{thm:rectangle-decomp}
\end{theorem}

\begin{remark}
In the case that $S$ and $T$ satisfy the conditions of \cref{thm:rectangle-decomp}, \cref{thm:block} follows easily from \cref{thm:rectangle-decomp}.
\end{remark}

\begin{remark}
Proofs of \cref{thm:block} for the special cases $S\times T$ finite and $S\times T=\R^2$ first appeared in \cite{bendich2013homology} and \cite{cochoy2020decomposition}, respectively; see also \cite{carlsson2009zigzag,carlsson2010zigzag}, which had a strong influence on these results. For the case $S\times T$ finite, an elementary proof of \cref{thm:rectangle-decomp} appeared in \cite{botnan_et_al:LIPIcs:2020:12180}.
\end{remark}

As in \cref{Sec:Interlevel_Rips}, let \[\overline{U} := \{ (a,b)\in\R^{\text{op}}\times \R \colon a \leq b\},\]
and let \[U:= \{(a,b)\in\R^{\text{op}}\times \R  \colon a <  b\}.\]
As a (non-trivial) application of \cref{thm:block}, \cite{botnan2020decomposition} establishes the following result.

\begin{theorem}[\cite{botnan2020decomposition}]
Let $P$ be $U$ or $\overline{U}$. Any p.f.d. middle-exact $P$-module $M$ is interval-decomposable, and each interval of $\B{M}$ is of the form $P\cap I$, where $I$ is a block in $\R^{\op}\times \R$. 
\label{thm:upperT}
\end{theorem} 
We shall also refer to the intervals $P\cap I$ arising in \cref{thm:upperT} as \emph{blocks (in $P$)}. 

\subsection{Interlevel Persistent Homology}\label{Sec:Interlevel_Sets}
One can use \cref{thm:upperT} to define \emph{interlevel barcodes}, fundamental invariants of $\R$-valued functions which refine the standard sublevel barcodes.  Moreover, these invariants admit extensions given using relative homology, and the invariants in different homological degrees assemble nicely into a single object, the \emph{Mayer-Vietoris strip}.  We will now explain this. 

To formulate the definitions, one has to choose between using open intervals or closed intervals in $\R$.  For consistency with our treatment of sublevel filtrations (Definition \ref{Def:Sublevel}), we will use closed intervals here, though using open intervals is arguably more convenient.  We briefly discuss the approach via open intervals in Remark \ref{Rem:Open_Interlevel} below.

Given a continuous function $\gamma\colon W\to \R$, the \emph{interlevel bifiltration of $\gamma$} is the functor 
\[\FI{\gamma}\colon \overline{U}\to \Top, \quad  \FI{\gamma}(s,t) = \{x\in X : s \leq \gamma(x) \leq t\}.\]
It follows from the Mayer--Vietoris sequence that for all $i\geq 0$, the restriction of $\HH_i\FI{\gamma}$ to $U$
is middle-exact.   Thus, if $\HH_i\FI{\gamma}$ is p.f.d., then \cref{thm:upperT} provides a barcode of the restriction of $\HH_i\FI{\gamma}$ to $U$.
  
Furthermore, assuming that $\gamma$ is a Morse function or, more generally, \emph{of Morse type} \cite[Section 2]{carlsson2009zigzag}, $\HH_i\FI{\gamma}$ is 
itself middle-exact, and \cref{thm:upperT} then provides a barcode $\B{\HH_i\FI{\gamma}}$ of $\HH_i\FI{\gamma}$ without having to restrict to $U$.

Taking the intersection of each interval of $\B{\HH_i\FI{\gamma}}$ with $\partial U=\{(x,x)\mid x\in \R\}$ and projecting onto the first coordinate, we recover the $i^{\mathrm{th}}$ \emph{levelset barcode} of $\gamma$, as introduced in \cite{carlsson2009zigzag}; we denote this as $\mathcal L_i(\gamma)$.   The intervals in this barcode encode how the homology of the fibers $\gamma^{-1}(r)$ changes as $r$ varies.

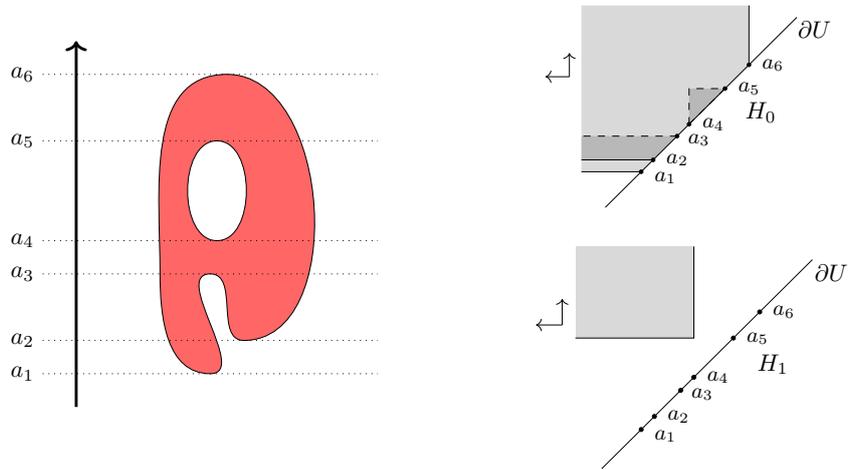
\begin{figure}
\centering
\scalebox{0.9}{
\begin{tikzpicture}[scale=.35]
\begin{scope}

\begin{scope}[rotate=90, xshift=-7cm, scale=1.4]
\fill[red!60,draw=black,even odd rule]
(0,0) to [out=90,in=180] (3,1.5) to [out=0,in=90] (9,-.5) to [out=270,in=270] (1,-1) to [out=90,in=270] (3,0) to [out=90,in=270] (0,0)
(4,-.2) to [out=90, in=90] (7,-.2) to [out=270, in=270] (4,-.2);
\draw[very thick, ->] (-1,4) -- (10,4);
\def\x{-5}
\def\y{5}
\draw[dotted] (0,\y) -- (0,\x);
\draw[dotted] (1,\y) -- (1,\x);
\draw[dotted] (3,\y) -- (3,\x);
\draw[dotted] (4,\y) -- (4,\x);

\draw[dotted] (7,\y) -- (7,\x);
\draw[dotted] (9,\y) -- (9,\x);

\node[left] at (0,\y) {$a_1$};
\node[left] at (1,\y) {$a_2$};
\node[left] at (3,\y) {$a_3$};
\node[left] at (4,\y) {$a_4$};
\node[left] at (7,\y) {$a_5$};
\node[left] at (9,\y) {$a_6$};

\end{scope}

\begin{scope}[xshift=18cm, scale=0.5, yshift=3cm]
\node at (10,5) {$H_0$};

\draw (-3,-3) -- (13,13); 

\node[right=2pt] (x1) at (0,-0.5) {\small $a_1$};
\node[right=2pt] (x2) at (1,1) {\small $a_2$};
\node[right=2pt] (x3) at (2.8,2.7) {\small $a_3$};
\node[right=2pt] (x4) at (4,4) {\small $a_4$};
\node[right=2pt] (x5) at (7,7) {\small $a_5$};
\node[right=2pt] (x6) at (9,9) {\small $a_6$};
\node[right=2pt] (x6) at (12,12) {$\partial U$};
\draw[draw=none, fill=black, fill opacity=0.15] (-5,0) -- (0,0) -- (9,9) -- (9,14) -- (-5,14) -- cycle;
\draw (-5,0) -- (0,0);
\draw (9,9) -- (9,14);

\draw[draw=none, fill=black, fill opacity=0.15]  (-5,1) -- (1,1) -- (3,3) -- (-5,3) -- cycle;
\draw (-5,1) -- (1,1);
\draw[dashed] (3,3) -- (-5,3);
\draw[->] (-6,8) -- (-6,10);
\draw[->] (-6,8) -- (-8,8);
\draw[dashed, fill=black, fill opacity=0.15]  (4,4) -- (7,7) -- (4,7) -- cycle;

\coordinate (y1) at (0,0);
\coordinate (y2) at (1,1);
\coordinate (y3) at (3,3);
\coordinate (y4) at (4,4);
\coordinate (y5) at (7,7);
\coordinate (y6) at (9,9);

\foreach \a in {y1, y2, y3, y4, y5, y6}
	\draw[fill=black] (\a) circle (1ex); 
\end{scope}

\begin{scope}[xshift=18cm, scale=0.55, yshift=-17cm]

\draw (-3,-3) -- (13,13); 

\node[right=2pt] (x1) at (0,-0.5) {\small $a_1$};
\node[right=2pt] (x2) at (1,1) {\small $a_2$};
\node[right=2pt] (x3) at (2.8,2.7) {\small $a_3$};
\node[right=2pt] (x4) at (4,4) {\small $a_4$};
\node[right=2pt] (x5) at (7,7) {\small $a_5$};
\node[right=2pt] (x6) at (9,9) {\small $a_6$};
\draw[->] (-6,8) -- (-6,10);
\draw[->] (-6,8) -- (-8,8);

\node[right=3pt] (x6) at (12,12) {$\partial U$};

\draw[draw=none, fill=black, fill opacity=0.15] (4,7) -- (4,14) -- (-5,14) -- (-5, 7) -- cycle;
\draw (-5,7) -- (4,7) -- (4,14);
\node at (10,5) {$H_1$};
\coordinate (y1) at (0,0);
\coordinate (y2) at (1,1);
\coordinate (y3) at (3,3);
\coordinate (y4) at (4,4);
\coordinate (y5) at (7,7);
\coordinate (y6) at (9,9);

\foreach \a in {y1, y2, y3, y4, y5, y6}
	\draw[fill=black] (\a) circle (1ex); 
\end{scope}

\end{scope}
\end{tikzpicture}}
\caption{The interlevel homology of the function $\gamma\colon W\to \R$ given by $\gamma(x,y)=y$, where $W\subset \R^2$ is as shown in red.  }
\label{fig:lzz-barcode}

\end{figure}

\begin{example}\label{Ex:Levelset}
\cref{fig:lzz-barcode} illustrates a function $\gamma\colon W\to \R$ and the barcodes 
$\B{\HH_i \FI{\gamma}}$.  By considering how the intervals in these barcodes intersect $\partial U$, we see that 
\[\mathcal L_0(\gamma)=\{[a_1, a_6], [a_2, a_3), (a_4, a_5)\}\] 
and $\mathcal L_i(\gamma)=\emptyset$ for $i> 0$.  
\label{ex:lzz-open}
\end{example}

\begin{figure}
\centering
 \scalebox{0.33}{
\begin{tikzpicture}[scale=0.6]
\begin{scope}[scale=1, xshift=0cm ]

\draw (-3,-3) -- (13,13) -- (-3,29) -- (-19,13) -- cycle;
\draw (-3,-3) -- (-3,29);
\draw (-19,13) -- (13,13);

\coordinate (y1) at (0,0);
\coordinate (y2) at (1,1);
\coordinate (y3) at (3,3);
\coordinate (y4) at (4,4);
\coordinate (y5) at (7,7);
\coordinate (y6) at (9,9);

\coordinate (hy1) at (0,26);
\coordinate (hy2) at (1,25);
\coordinate (hy3) at (3,23);
\coordinate (hy4) at (4,22);
\coordinate (hy5) at (7,19);
\coordinate (hy6) at (9,17);

\coordinate (ly1) at (-6,0);
\coordinate (ly2) at (-7,1);
\coordinate (ly3) at (-9,3);
\coordinate (ly4) at (-10,4);
\coordinate (ly5) at (-13,7);
\coordinate (ly6) at (-15,9);

\coordinate (xx) at (2,2);
\coordinate (yy) at = (8,8);

\coordinate (z1) at (-15,17);
\coordinate (z2) at (-13,19);
\coordinate (z3) at (-10,22);
\coordinate (z4) at (-9,23);
\coordinate (z5) at (-7,25);
\coordinate (z6) at (-6,26);

\foreach \a in {z1, z2,z3,z4,z5,z6, y1,y2,y3,y4,y5,y6,ly1,ly2,ly3,ly4,ly5,ly6, hy1,hy2,hy3,hy4,hy5,hy6}
	\draw[fill=black] (\a) circle (1ex);

\node[above left=-1pt] at (z1) {\huge $a_6$};
\node[above left=-1pt] at (z2) {\huge $a_5$};
\node[above left=-1pt] at (z3)  {\huge $a_4$};
\node[above left=-1pt] at (z4) {\huge $a_3$};
\node[above left=-1pt] at (z5) {\huge $a_2$};
\node[above left=-1pt] at (z6) {\huge $a_1$};

\node[below right=-1pt] at (y1) {\huge $a_1$};
\node[below right=-1pt] at (y2) {\huge $a_2$};
\node[below right=-1pt] at (y3)  {\huge $a_3$};
\node[below right=-1pt] at (y4) {\huge $a_4$};
\node[below right=-1pt] at (y5) {\huge $a_5$};
\node[below right=-1pt] at (y6) {\huge $a_6$};

\node[below left=-1pt] at (ly1) {\huge $a_1$};
\node[below left=-1pt] at (ly2) {\huge $a_2$};
\node[below left=-1pt] at (ly3)  {\huge  $a_3$};
\node[below left=-1pt] at (ly4) {\huge $a_4$};
\node[below left=-1pt] at (ly5) {\huge $a_5$};
\node[below left=-1pt] at (ly6) {\huge $a_6$};

\node[above right=-1pt] at (hy1) {\huge $a_1$};
\node[above right=-1pt] at (hy2) {\huge $a_2$};
\node[above right=-1pt] at (hy3)  {\huge $a_3$};
\node[above right=-1pt] at (hy4) {\huge $a_4$};
\node[above right=-1pt] at (hy5) {\huge $a_5$};
\node[above right=-1pt] at (hy6) {\huge $a_6$};

\draw[fill=black] (2,8) circle (1ex);
\draw[dashed] (xx) -- (2,8) -- (yy);

\node[left] at (1,10) {\Huge $Q_1$};

\draw[fill=black] (-8,8) circle (1ex);
\draw[dashed] (-8,2) -- (-8,8) -- (-14,8);

\node[right] at (-7,10) {\Huge $Q_3$};

\draw[fill=black] (-8,18) circle (1ex);
\draw[dashed] (-8,24) -- (-8,18) -- (-14,18);

\node[right] at (-7,16) {\Huge $Q_4$};

\draw[fill=black] (2,18) circle (1ex);
\draw[dashed] (8,18) -- (2,18) -- (2,24);
\node[left] at (1,16) {\Huge $Q_2$};

\foreach \a in {(2,24), (8,18),(xx), (yy), (-8,24), (-14,18), (-8,2), (-14,8)}
	\draw[fill=red] \a circle (2ex); 
	
\foreach \a in {(1.7,24),(1.7,2)}
	\node[left] at \a {\huge $x$};

\node[above] at (7.7,8.3) {\huge $y$};
\node[above] at (-13.7,8.3) {\huge $y$};

\node[below] at (-13.7,17.7) {\huge $y$};
\node[below] at (7.3,17.7) {\huge $y$};

\node[right] at (-7.7,24) {\huge $x$};
\node[right] at (-7.7,2) {\huge $x$};

\draw[->,ultra thick] (-4,14) -- (-4,16);
\draw[->,ultra thick] (-4,14) -- (-6,14);

\end{scope}
\end{tikzpicture}}
\caption{Visualization of the poset $Q$ indexing relative interlevel homology.}
\label{fig:pyramid}
\end{figure}
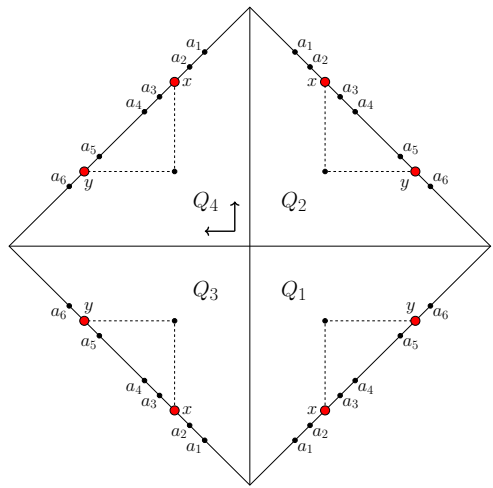

\paragraph{Relative Interlevel Homology}
For each $i\geq 0$, we can extend the above construction  
 to a relative version as follows: consider the set \[Q=Q_1\sqcup Q_2 \sqcup Q_3\sqcup Q_4\] of pairs of subsets of $\R$, where
\begin{align*} 
Q_1 &= \{([x,y],\emptyset )\colon x\leq y\} & Q_2 &= \{([x, \infty),[y, \infty))\colon x\leq y\}),\\
Q_3 &= \{((-\infty,y],(-\infty,x])\colon x\leq y\}, & Q_4 &= \{(\R,(-\infty, x]\cup [y, \infty))\colon x\leq y\}.
\end{align*}
We regard $Q$ as poset by taking $(A,B)\leq (A',B')$ if and only if $A\subseteq A'$ and $B\subseteq B'$.  Let $2\R$ denote the totally ordered set consisting of two disjoint copies of $\R$, with all elements in the first copy less than all elements in the second copy.  We have a natural identification of $Q$ with an interval in $(2\R^{\op})\times 2\R$.    
  
Given a continuous function $\gamma\colon W\to \R$, we define a functor $\mathcal{S}^{\bullet}(\gamma)$ from $Q$ to the category of pairs of topological spaces by \[\mathcal{S}(\gamma)^{\bullet}(A,B)=(\gamma^{-1}(A),\gamma^{-1}(B)).\]  If $\gamma$ is of Morse type, then $H_i\mathcal{S}^{\bullet}(\gamma)$ is interval-decomposable, where each interval is the intersection of a rectangle in $2\R^{\op}\times 2\R$ with $Q$ \cite{bendich2013homology,botnan_et_al:LIPIcs:2020:12180}. 
It is helpful to visualize $Q$ schematically as the decomposition of a tilted square into four triangles, as shown in \cref{fig:pyramid}; we can make formal sense of this schematic either by 
reparameterizing $Q$ as in \cite{bauer2021structure}, or by assuming that $\gamma$ is bounded above and below, in which case we may replace $Q$ with a subposet where $x$ and $y$ are correspondingly bounded, as in \cite{carlsson2009zigzag}.

\begin{figure}
\centering\centering
 \scalebox{1}{
\begin{tikzpicture}[scale=0.15, rotate=0]
\begin{scope}

\begin{scope}[scale=0.8]
\def\x{25}
\def\y{5}

\coordinate (y1) at (0,0);
\coordinate (y2) at (1,1);
\coordinate (y3) at (3,3);
\coordinate (y4) at (4,4);
\coordinate (y5) at (7,7);
\coordinate (y6) at (9,9);

\coordinate (w1) at ($ (y1) + (\x,-\x)$);
\coordinate (w2) at ($ (y2) + (\x,-\x)$);
\coordinate (w3) at ($ (y3) + (\x,-\x)$);
\coordinate (w4) at ($ (y4) + (\x,-\x)$);
\coordinate (w5) at ($ (y5) + (\x,-\x)$);
\coordinate (w6) at ($ (y6) + (\x,-\x)$);

\coordinate (x1) at (-15+\x,17-\x);
\coordinate (x2) at (-13+\x,19-\x);
\coordinate (x3) at (-10+\x,22-\x);
\coordinate (x4) at (-9+\x,23-\x);
\coordinate (x5) at (-7+\x,25-\x);
\coordinate (x6) at (-6+\x,26-\x);

\coordinate (ww1) at ($ (x1) + (\x,-\x)$);
\coordinate (ww6) at ($ (x6) + (\x,-\x)$);

\coordinate (z1) at (-15,17);
\coordinate (z2) at (-13,19);
\coordinate (z3) at (-10,22);
\coordinate (z4) at (-9,23);
\coordinate (z5) at (-7,25);
\coordinate (z6) at (-6,26);

\foreach \a in {z1, z2,z3,z4,z5,z6, y1,y2,y3,y4,y5,y6, w1,w2,w3,w4,w5,w6,x1,x2,x3,x4,x5,x6}
	\draw[fill=black] (\a) circle (1ex);

\node[above left=-2pt] at (x1) {\tiny $a_6$};
\node[above left=-2pt] at (x2) {\tiny $a_5$};
\node[above left=-2pt] at (x3)  {\tiny $a_4$};
\node[above left=-2pt] at (x4) {\tiny $a_3$};
\node[above left=-2pt] at (x5) {\tiny $a_2$};
\node[above left=-2pt] at (x6) {\tiny $a_1$};

\node[above left=-2pt]  at (z1) {\tiny $a_6$};
\node[above left=-2pt]  at (z2) {\tiny $a_5$};
\node[above left=-2pt]  at (z3)  {\tiny $a_4$};
\node[above left=-2pt]  at (z4) {\tiny $a_3$};
\node[above left=-2pt]  at (z5) {\tiny $a_2$};
\node[above left=-2pt]  at (z6) {\tiny $a_1$};

\node[below right=-2pt] at (y1) {\tiny $a_1$};
\node[below right=-2pt]  at (y2) {\tiny $a_2$};
\node[below right=-2pt]  at (y3)  {\tiny $a_3$};
\node[below right=-2pt]  at (y4) {\tiny $a_4$};
\node[below right=-2pt]  at (y5) {\tiny $a_5$};
\node[below right=-2pt]  at (y6) {\tiny $a_6$};

\node[below right=-2pt]  at (w1) {\tiny $a_1$};
\node[below right=-2pt]  at (w2) {\tiny $a_2$};
\node[below right=-2pt]  at (w3)  {\tiny $a_3$};
\node[below right=-2pt]  at (w4) {\tiny $a_4$};
\node[below right=-2pt]  at (w5) {\tiny $a_5$};
\node[below right=-2pt]  at (w6) {\tiny $a_6$};

\draw (-3,-3) -- (13,13) -- (-3,29) -- (-19,13) -- cycle;
\draw (-3+\x,-3-\x) -- (13+\x,13-\x) -- (-3+\x,29-\x) -- (-19+\x,13-\x) -- cycle;

\draw[dashed, draw=none, fill=black, fill opacity=0.15] (-6,0) -- (0,0) -- (9,9) -- (9,17) -- (-6,17) -- cycle;
\draw (-6,0) -- (0,0) -- (9,9) -- (9,17);
\draw[dashed] (9,17) -- (-6,17) -- (-6,0);

\draw[draw=none, fill=black, fill opacity=0.15]  (-7,1) -- (1,1) -- (3,3) -- (-7,3) -- cycle;
\draw (-7,1) -- (1,1);
\draw[dashed] (3,3) -- (-7,3) -- (-7,1);

\draw[dashed, fill=black, fill opacity=0.15]  (4,4) -- (7,7) -- (4,7) -- cycle;

\draw[draw=none, fill=black, fill opacity=0.2] (4+\x,7-\x) -- (4+\x,22-\x) -- (-10+\x,22-\x) -- (-13+\x,19-\x) -- (-13+\x,7-\x)--  cycle;
\draw (-13+\x,7-\x) -- (4+\x,7-\x) -- (4+\x,22-\x);
\draw[dashed] (4+\x,22-\x) -- (-10+\x,22-\x);
\draw[dashed] (-13+\x,7-\x) -- (-13+\x,19-\x);

\draw (-7+\x,25-\x) -- (-7+\x,3-\x)--(-9+\x,3-\x);
\draw[dashed] (-9+\x,23-\x) -- (-9+\x,3-\x);
\draw[draw=none, fill=black, fill opacity=0.2] (-7+\x,25-\x) -- (-7+\x,3-\x) -- (-9+\x,3-\x) -- (-9+\x,23-\x) -- cycle;

\draw (-6+\x,26-\x) -- (-6+\x,17-\x) -- (-15+\x, 17-\x);
\draw[dashed, draw=none, fill=black, fill opacity=0.2] (-6+\x,26-\x) -- (-6+\x,17-\x) -- (-15+\x, 17-\x); -- cycle;

\draw (-3,-3) -- (-3,29);
\draw (-19,13) -- (13,13);

\draw (-3+\x,-3-\x) -- (-3+\x,29-\x);
\draw (-19+\x,13-\x) -- (13+\x,13-\x);

\node at (-10,0) {$H_0$};
\node at (-10+\x,-\x) {$H_1$};

\draw [thick,->] (x6) to [out=150,in=-120] (y1);
\draw [thick, ->] (x1) to [out=150,in=30] (y6);

\draw [dashed,thick,->] (ww6) to [out=150,in=-120] (w1);
\draw [dashed, thick, ->] (ww1) to [out=150,in=30] (w6);

\end{scope}

\begin{scope}[xshift=25cm,yshift=15cm, scale=0.8]

\node[rotate=-45] at (12,18) {Mayer--Vietoris strip};

\def\x{16}
\def\y{5}
\draw (-3,-3) -- (13,13) -- (-3,29) -- (-19,13) -- cycle;
\draw (-3+\x,-3-\x) -- (13+\x,13-\x) -- (-3+\x,29-\x) -- (-19+\x,13-\x) -- cycle;
\draw (-3+\x, -3-\x) -- (-3+1.5*\x, -3-1.5*\x);
\draw (13+\x, 13-\x) -- (13+1.5*\x, 13-1.5*\x);

\coordinate (y1) at (0,0);
\coordinate (y2) at (1,1);
\coordinate (y3) at (3,3);
\coordinate (y4) at (4,4);
\coordinate (y5) at (7,7);
\coordinate (y6) at (9,9);

\coordinate (x1) at (-15+2*\x,17-2*\x);
\coordinate (x2) at (-13+2*\x,19-2*\x);
\coordinate (x3) at (-10+2*\x,22-2*\x);
\coordinate (x4) at (-9+2*\x,23-2*\x);
\coordinate (x5) at (-7+2*\x,25-2*\x);
\coordinate (x6) at (-6+2*\x,26-2*\x);

\coordinate (z1) at (-15,17);
\coordinate (z2) at (-13,19);
\coordinate (z3) at (-10,22);
\coordinate (z4) at (-9,23);
\coordinate (z5) at (-7,25);
\coordinate (z6) at (-6,26);

\foreach \a in {z1, z2,z3,z4,z5,z6, y1,y2,y3,y4,y5,y6, x1,x2,x3,x4,x5,x6}
	\draw[fill=black] (\a) circle (1ex);

\node[below right=-2pt] at (x1) {\tiny $a_6$};
\node[below right=-2pt] at (x2) {\tiny  $a_5$};
\node[below right=-2pt] at (x3)  {\tiny $a_4$};
\node[below right=-2pt] at (x4) {\tiny $a_3$};
\node[below right=-2pt] at (x5) {\tiny $a_2$};
\node[below right=-2pt] at (x6) {\tiny $a_1$};

\node[above left=-2pt] at (z1) {\tiny $a_6$};
\node[above left=-2pt] at (z2) {\tiny $a_5$};
\node[above left=-2pt] at (z3)  {\tiny $a_4$};
\node[above left=-2pt] at (z4) {\tiny $a_3$};
\node[above left=-2pt] at (z5) {\tiny $a_2$};
\node[above left=-2pt] at (z6) {\tiny $a_1$};

\draw[dashed, draw=none, fill=black, fill opacity=0.15] (-6,0) -- (9,0) -- (9,17) -- (-6,17) -- cycle;
\draw (-6,0) -- (9,0) -- (9,17);
\draw[dashed] (9,17) -- (-6,17) -- (-6,0);

\draw[dashed, draw=none, fill=black, fill opacity=0.15]  (-7,1) -- (23,1) -- (23,3) -- (-7,3) -- cycle;
\draw (-7,1) -- (23,1) -- (23,3);
\draw[dashed]  (-7,1) -- (-7,3) -- (23,3);

\draw[draw=none, fill=black, fill opacity=0.15]  (19,7) -- (4,7) -- (4,-10) -- (19,-10);
\draw (4,-10) -- (19,-10) -- (19,7);
\draw[dashed] (19,7) -- (4,7) -- (4,-10);

\end{scope}

\end{scope}
\end{tikzpicture}}
\caption{Relative interlevel set modules $H_i\mathcal S^{\bullet}(\gamma)$ and the corresponding Mayer--Vietoris strip of the function $\gamma$ from Example~\ref{ex:lzz-open}; see Examples~\ref{ex:pyramid} and \ref{ex:strip} and the surrounding discussion.}
\label{fig:mv-strip}

\end{figure}
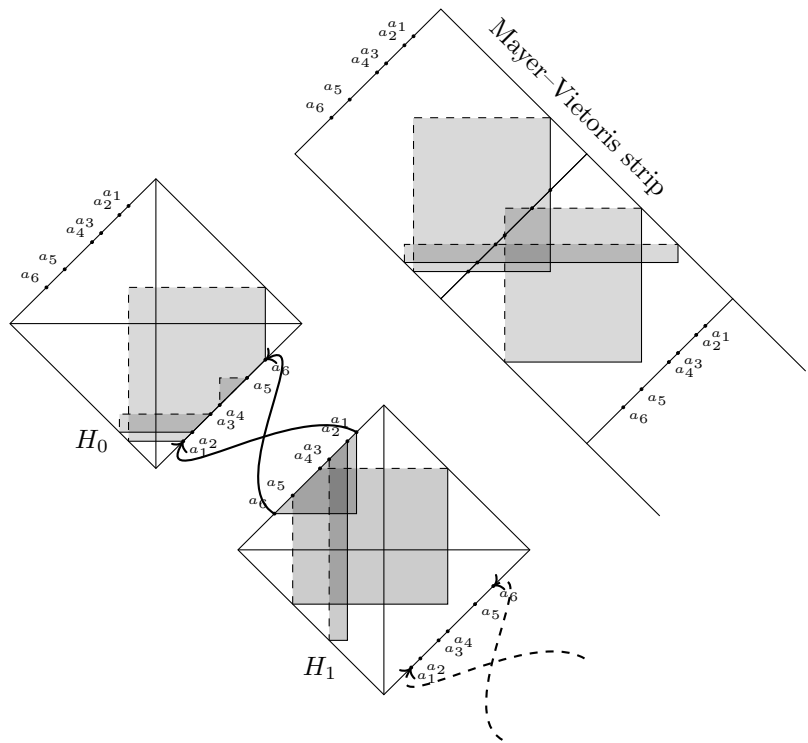

\begin{example}
\cref{fig:mv-strip} shows the decompositions of the modules $H_0\mathcal{S}^{\bullet}(\gamma)$ and $H_1\mathcal{S}^{\bullet}(\gamma)$ for $\gamma$ the function of~\cref{Ex:Levelset} and ~\cref{fig:pyramid}. Observe that the off-diagonal block in $\HH_1\FI{\gamma}$ now extends all the way to the top of $Q$. Moreover, if we reflect the block around the anti-diagonal passing through the center of $Q$, then we see that the boundary of the block aligns perfectly with the boundary of the triangle from $\HH_0\FI{\gamma}$. Geometrically, this makes sense, as the circle in $\gamma^{-1}[a_4, a_5]$ gives rise to two connected components in all pre-images $\gamma^{-1}[x,y]$ where $[x,y]\subseteq (a_4, a_5)$.
\label{ex:pyramid}
\end{example}

\paragraph{The Mayer--Vietoris Strip} In the previous example, there is a natural matching between $\B{H_0\mathcal{S}^{\bullet}(\gamma)}$ and $\B{H_1\mathcal{S}^{\bullet}(\gamma)}$.  It turns out that one always has such a matching.  To explain, using excision and the boundary map in the relative Mayer--Vietoris sequence, one obtains maps connecting the diagrams $H_i\mathcal{S}^{\bullet}(\gamma)$ and $H_{i+1}\mathcal{S}^{\bullet}(\gamma)$ for all $i$.  Using these maps and reflecting every other diagram $H_i\mathcal{S}^{\bullet}$ as shown in \cref{fig:mv-strip}, the diagrams $\{H_i\mathcal{S}^{\bullet}(\gamma)\}_{i\in \mathbb N}$ assemble into one big persistence module called the \emph{Mayer--Vietoris strip} of $\gamma$.  Let $Q'$ denote the indexing poset of this persistence module.  Like $Q$, the poset $Q'$ is naturally identified with an interval in a product of two totally ordered sets, so rectangles in $Q'$ are well defined.  

It turns out that the Mayer--Vietoris strip is middle-exact.  Using this, one can prove the following structure theorem.

\begin{theorem}[\cite{bauer2021structure,bendich2013homology,botnan2020local}]\label{Thm:MV_Structure}
If $\gamma$ is of Morse type, then the Mayer--Vietoris strip of $\gamma$ is interval-decomposable, where each interval is a maximal rectangle in $ Q'$. 
\end{theorem}

\begin{example}The Mayer--Vietoris strip of the function in~\cref{fig:pyramid}\,(a) is shown in \cref{fig:mv-strip}.
\label{ex:strip}
\end{example}

\begin{remark}[Computation]
The barcode of the Mayer--Vietoris strip of $\gamma$ is fully determined by the levelset barcodes of $\gamma$ in all homological degrees, and is also fully determined by the \emph{extended persistence} barcodes \cite{cohen2009extending} of $\gamma$ in all homological degrees; this is implicit in \cite{carlsson2009zigzag}.  It follows that in practice, one can obtain the barcode of the Mayer--Vietoris strip by either a (1-parameter) extended persistence computation or a zigzag persistence computation.  Conversely, the barcode of the Mayer-Vietoris strip determines both the levelset and extended persistence barcodes.
\end{remark}

\begin{remark}[Stability]
\label{rem:stability-strip}
The bottleneck distance (Definitions \ref{Def:bottleneck} and \ref{Def:bottleneck_general}) generalizes straightforwardly to barcodes consisting of rectangles in $Q'$, and it can be shown that the barcode of the Mayer--Vietoris strip is stable with respect to this and the $L^\infty$-distance on functions \cite{carlsson2009zigzag,bauer2021structure}. This is a direct consequence of the stability of extended persistent homology \cite{cohen2009extending}. 
\end{remark}
 
\begin{remark}[Interlevel Persistence with Open Intervals]\label{Rem:Open_Interlevel}
We now briefly discuss the corresponding interlevel persistence theory using open intervals.  
Given a continuous function, $\gamma:W\to \R$, we can define
\[\FI{\gamma}^o\colon U \to {\rm Top}, \qquad \FI{\gamma}^o(s,t) = \{x\in X \mid s< \gamma(x) < t\}.\]
By the Mayer--Vietoris sequence, $\FI{\gamma}^o$ is middle-exact.  Thus by \cref{thm:upperT}, if $H_i\FI{\gamma}^o$ is p.f.d., then it is interval-decomposable.  Analogously, we can define a variant of the Mayer--Vietoris strip using open intervals, and we have a structure theorem \cite{botnan2020local,bauer2021structure} which says that if this persistence module is p.f.d., then it is interval-decomposable where the intervals are maximal rectangles. As in the case of closed intervals, this decomposition is stable.
\end{remark}

\begin{remark}[Related Results]
See Bauer et al. \cite{bauer2021structure} for a careful construction of the Mayer--Vietoris strip for an arbitrary cohomology theory, and a corresponding structure theorem in the p.f.d. setting. A structure theorem for middle-exact p.f.d. modules over the Mayer--Vietoris strip not necessarily coming from functions can be found in \cite{botnan2020local}. A related construction and structure theorem in the context of derived sheaf theory can be found in \cite{berkouk2019level}.  
\end{remark}

\bibliographystyle{plain}
\bibliography{references}

\end{document}